\documentclass[12pt]{amsart}
\usepackage[]{amscd,amsfonts,amsmath,amsxtra,amssymb}
\usepackage{graphics}
\usepackage[all]{xy}
\newtheorem{sub}{}[section]
\newtheorem{subsub}{}[sub]
\font\ttx=cmcsc10

\font\tte=cmbsy10

\def\ov#1{\overline{#1}}
\def\codim{\mathop{\rm codim}\nolimits}
\def\coker{\mathop{\rm coker}\nolimits}
\def\Hom{\mathop{\rm Hom}\nolimits}
\def\HHom{\mathop{\mathcal Hom}\nolimits}
\def\Ext{\mathop{\rm Ext}\nolimits}
\def\EExt{\mathop{\mathcal Ext}\nolimits}

\def\Hilb{\mathop{\rm Hilb}\nolimits}
\def\Aut{\mathop{\rm Aut}\nolimits}
\def\Ad{\mathop{\rm Ad}\nolimits}
\def\End{\mathop{\rm End}\nolimits}
\def\imm{\mathop{\rm Im}\nolimits}
\def\lra{\longrightarrow}
\def\som{\mathop{\hbox{$\displaystyle\bigoplus$}}\limits}
\def\ssom{\mathop{\hbox{$\displaystyle\oplus$}}\limits}
\def\sigg{\mathop{\hbox{$\displaystyle\sum$}}\limits}

\def\supp{\mathop{\hbox{$\sup$}}\limits}

\def\paragra{{\tte \char120}}
\def\para{\paragra~\hskip -2pt}

\def\hfl#1#2{\smash{\mathop{\ \hbox to 12mm{\rightarrowfill}}
\limits^{\scriptstyle#1}_{\scriptstyle#2} \ }}
\def\hflb#1#2{\smash{\mathop{\hbox to 12mm{\leftarrowfill}}
\limits^{\scriptstyle#1}_{\scriptstyle#2}}}

\def\m#1{{\hbox{$#1$}}}

\def\ot{\otimes}
\def\og{\leavevmode\raise.3ex\hbox{$\scriptscriptstyle\langle\!\langle$}}
\def\fg{\leavevmode\raise.3ex\hbox{$\scriptscriptstyle\,\rangle\!\rangle$}}

\def\nsp{\lbrace 0\rbrace}

\def\dsp{\displaystyle}
\def\Ssect#1#2{\pagebreak[3]\begin{sub}\label{#2}{\sc\small\small  #1}\rm\medskip}
\def\SSsect#1#2{\pagebreak[2]\begin{subsub}\label{#2}{#1}\end{subsub}}
\def\sepsec{\vskip 2.5cm}
\def\sepsub{\vskip 1.5cm}
\def\sepsubsub{\vskip 1cm}
\def\sepprop{\vskip 0.8cm}
\def\sepx{\vskip 0.6cm}
\def\xmat#1{\[\xymatrix{#1}\]}
\def\Bdiag#1#2#3#4#5#6#7#8#9{\[{\xymatrix{
& 0\ar[d] & 0\ar[d] & 0\ar[d] \\
0\ar[r] & #1\ar[r]\ar[d] & #2\ar[d]\ar[r] & #3\ar[d]\ar[r] & 0 \\
0\ar[r] & #4\ar[r]\ar[d] & #5\ar[d]\ar[r] & #6\ar[d]\ar[r] & 0 \\
0\ar[r] & #7\ar[r]\ar[d] & #8\ar[d]\ar[r] & #9\ar[d]\ar[r] & 0 \\
& 0 & 0 & 0 & \\
}}\]}
\def\flinc{\ar@{^{(}->}}
\def\fleq{\ar@{=}}
\def\flon{\ar@{->>}}
\def\fmaps{\ar@{|-{>}}}
\def\wT{{\widetilde T}}
\newcommand{\N}{{\mathbb N}}
\newcommand{\M}{{\mathbb M}}
\newcommand{\A}{{\mathbb A}}
\newcommand{\B}{{\mathbb B}}
\newcommand{\Z}{{\mathbb Z}}
\newcommand{\R}{{\mathbb R}}
\newcommand{\C}{{\mathbb C}}
\newcommand{\Q}{{\mathbb Q}}
\renewcommand{\P}{{\mathbb P}}
\newcommand{\D}{{\mathbb D}}
\newcommand{\F}{{\mathbb F}}
\newcommand{\E}{{\mathbb E}}
\newcommand{\G}{{\mathbb G}}

\newcommand{\V}{{\mathbb V}}
\newcommand{\W}{{\mathbb W}}

\def\T{{\mathbb T}}
\newcommand{\ka}{{\mathcal A}}

\newcommand{\kd}{{\mathcal D}}
\newcommand{\ke}{{\mathcal E}}
\newcommand{\kf}{{\mathcal F}}
\newcommand{\kg}{{\mathcal G}}
\newcommand{\kh}{{\mathcal H}}
\newcommand{\ki}{{\mathcal I}}
\newcommand{\kj}{{\mathcal J}}
\newcommand{\kk}{{\mathcal K}}
\newcommand{\kl}{{\mathcal L}}

\newcommand{\kn}{{\mathcal N}}
\newcommand{\ko}{{\mathcal O}}
\newcommand{\kp}{{\mathcal P}}

\newcommand{\kt}{{\mathcal T}}
\newcommand{\ku}{{\mathcal U}}
\newcommand{\kv}{{\mathcal V}}
\newcommand{\kw}{{\mathcal W}}
\newcommand{\kx}{{\mathcal X}}
\newcommand{\ky}{{\mathcal Y}}
\newcommand{\kz}{{\mathcal Z}}
\begin{document}
\newtheorem{xprop}{Proposition}[section]
\newtheorem{xlemm}[xprop]{Lemme}
\newtheorem{xtheo}[xprop]{Th\'eor\`eme}
\newtheorem{xcoro}[xprop]{Corollaire}
\newtheorem{quest}{Question}
\newtheorem{defin}[xprop]{D\'efinition}
\newtheorem{remar}[xprop]{Remarque}
\newtheorem{remars}[xprop]{Remarques}

\def\refname{R\'ef\'erences}
\def\contentsname{Sommaire}
\def\proofname{D\'emonstration}
\def\abstractname{R\'esum\'e}

\author{Jean--Marc Dr\'{e}zet}
\address{
Institut de Math\'{e}matiques de Jussieu\\
UMR 7586 du CNRS\\ 
175 rue du Chevaleret\newline   
F-75013 Paris, France}
\email{drezet@math.jussieu.fr}

\title[{\tiny D\'eformations des extensions larges de faisceaux}]
{D\'eformations des extensions larges de faisceaux}
\maketitle
\tableofcontents


\section{Introduction}
\label{intro}
Soient $X$ une vari\'et\'e alg\'ebrique projective lisse irr\'eductible sur
$\C$, de dimension $d>0$ et $\ko_X(1)$ un fibr\'e en droites tr\`es ample sur
$X$. Si $\kf$ est un faisceau coh\'erent sur $X$ tel que $rg(\kf)>0$, on pose
$$\mu(\kf) \ = \frac{c_1(\kf).\ko_X(1)^{d-1}}{rg(\kf)} .$$
Il est bien connu que les faisceaux semi-stables relativement \`a $\ko_X(1)$,
de rang et de classes de Chern donn\'es sur $X$ admettent une vari\'et\'e de
modules grossiers, qui est dans certains cas une vari\'et\'e de modules fins
(cf. \cite{ma1}, \cite{ma2}, \cite{si}, \cite{yo}). D'autres ensembles de 
classes
d'isomorphisme de faisceaux coh\'erents sur $X$ peuvent aussi avoir une 
vari\'et\'e de modules fins (cf. \cite{dr}). Les cas trait\'es dans \cite{dr}
concernent des faisceaux qui ne sont pas loin d'\^etre semi-stables (par 
exemple
les faisceaux {\em prioritaires} sur $\P_2$ d\'efinis dans \cite{hi_la}).
On s'int\'eresse ici \`a des faisceaux {\em tr\`es instables}: on va 
\'etudier des
fibr\'es vectoriels $E$ poss\`edant un sous-faisceau $\Gamma$ tel que \
\m{\mu(\Gamma)\gg\mu(E)}. On s'int\'eresse aux d\'eformations de tels fibr\'es.

Rappelons qu'un faisceau coh\'erent $\Gamma$ sur $X$ poss\`ede une {\em
d\'eformation semi-universelle} $\kg$, qui est un faisceau 
analytique sur
$S\times X$, $(S,s_0)$ \'etant un germe de vari\'et\'e analytique, tel que \
\m{\kg_{s_0}\simeq\Gamma}. De plus, l'espace tangent $T_{s_0}S$ est
canoniquement isomorphe \`a $\Ext^1(\Gamma,\Gamma)$ (cf. \ref{S_T},
\cite{si_tr}). Le germe $(S,s_0)$ s'appelle la {\em base} de la d\'eformation
semi-universelle de $\Gamma$. Dans cet article on dira que $\Gamma$ est {\em
lisse} si $S$ est lisse en \m{s_0}. On dit que $\Gamma$ est {\em 2-lisse} si le
morphisme trace \m{\Ext^2(\Gamma,\Gamma)\to H^2(\ko_X)} est injectif (cf.
\ref{trace}). Un faisceau 2-lisse est lisse. 

\sepsub

\pagebreak[3]\begin{sub}{\sc\small\small Fibr\'es instables de rang 2 sur 
$\P_n$}
\label{intro-1}\rm

Les fibr\'es vectoriels instables de rang 2 sur $\P_2$ ou $\P_3$ ont d\'ej\`a
\'et\'e \'etudi\'es.
Soit \m{n\geq 2} un entier.
Pour \'etudier les fibr\'es vectoriels alg\'ebriques de rang 2 sur \m{\P_n} on
utilise la {\em construction de Serre} qui \'etablit un lien entre 
ces fibr\'es
et les sous-vari\'et\'es localement intersections compl\`etes de codimension 2
de \m{\P_n} (cf. \cite{O_S_S}, chap. I, \para 5, \cite{ho}, \cite{ba_vdv}, 
\cite{ha0}).
Soit $E$ un fibr\'e vectoriel alg\'ebrique de rang 2 et de classes de Chern 
$c_1=0$ ou $1$, $c_2$ sur $\P_n$. Soit $d$ le plus grand entier $\delta$
tel que \ $h^0(E(-\delta))>0$. Alors $E$ est non semi-stable si et seulement si
\m{d>0}. Si  c'est le cas on dit que $d$ est le 
{\em degr\'e d'instabilit\'e} de $E$. On a une suite exacte
$$(*) \ \ \ \ \ \ 0\lra\ko(d)\lra E\lra\ki_Z(c_1-d)\lra 0$$
$Z$ \'etant une sous-vari\'et\'e ferm\'ee localement intersection compl\`ete
de codimension 2 de \m{\P_n}, $\ki_Z$ d\'esignant son faisceau d'id\'eaux.

Le cas a de $\P_2$ \'et\'e trait\'e par S.A. Str\o mme dans \cite{st} 
($Z$ est alors un sous-sch\'ema de dimension 0, de longueur \ 
\m{m=c_2+d(d-c_1)}). La construction de Serre permet de construire une 
"vari\'et\'e de
modules" pour de tels fibr\'es : c'est la vari\'et\'e $M(d,c_1,c_2)$ des paires
$(Z,\sigma)$ avec \ \m{Z\in\Hilb^m(\P_2)}, $\sigma$ parcourant l'ouvert de
\m{\P(\Ext^1(\ki_Z(c_1-d),\ko(d)))} correspondant aux extensions $(*)$ telles
que $E$ soit localement libre. Str\o mme montre que les d\'eformations d'un
fibr\'e de rang 2 de degr\'e d'instabilit\'e $d$ g\'en\'erique sont aussi des 
fibr\'es de degr\'e d'instabilit\'e $d$. Cependant on a
$$\dim(\Ext^1(E,E)) \ > \ \dim(M(d,c_1,c_2)) ,$$
ce qui signifie que la base d'une d\'eformation universelle de $E$ n'est pas
r\'eduite. En particulier, $M(d,c_1,c_2)$ ne peut pas \^etre une vari\'et\'e de 
modules fins.

Le cas de $\P_3$ a \'et\'e trait\'e par C. B\u anic\u a dans \cite{ba} (voir
aussi \cite{ha2}, \para 1). La construction de Serre permet aussi dans ce cas de
construire une sorte de vari\'et\'e de modules pour les fibr\'es de degr\'e
d'instabilit\'e donn\'e.

\end{sub}

\sepsub

\pagebreak[3]\begin{sub}{\sc\small\small Extensions larges}\label{intro-2}\rm

Si $T$ est un faisceau pur sur $X$, on note
$$\widetilde{T} \ = \ T^\vee\ot\omega_X^{-1} .$$
On dit que $T$ est {\em parfait} si \ \m{\EExt^q(T,\ko_X)=0} \ pour
\m{q>\codim(T)} (cf. \ref{grat}).

Soit $F$ un faisceau coh\'erent non nul sur $X$. On dit que $F$ est {\em
r\'egulier} si $F$ est sans torsion, \m{F^{**}} simple, \m{F^{**}}
2-lisse, et si \m{F^{**}/F} est nul ou parfait de codimension 2 (cf. 
\ref{f_reg}).

Soient $\Gamma$ un faisceau localement libre, $F$ un faisceau r\'egulier sur
$X$, \m{T=F^{**}/F}. On consid\'ere une extension
$$0\lra\Gamma\lra\ke\lra F\lra O .$$
On en d\'eduit un morphisme \ \m{\phi : 
\Gamma^*\to\EExt^1(F,\ko_X)\simeq\widetilde{T}}. Il est ais\'e
de voir que $\ke$ est localement libre si et seulement si $\phi$ est surjectif.
Supposons que ce soit le cas et que \m{G=\ker(\phi)} soit r\'egulier. On obtient
alors deux extensions
\[(L) \ \ \ \ \ \ \ \ 0\lra G^*\lra\ke\lra F\lra O ,\]
\[(L^*) \ \ \ \ \ \ \ \ 0\lra F^*\lra\ke^*\lra G\lra O .\]
On dit que \m{(L)} (ou $\ke$) est une {\em extension large} si on a :
\[\Ext^i(F^{**},G^*) \ = \ \nsp \ \ \ \ {\rm si} \ \ i\geq 1 , \]
\[\Ext^i(F,G^*) \ = \ \Ext^i(G,F^*) \ = \ \nsp \ \ \ \ 
{\rm si} \ \ i\geq 2 , \]
\[\Ext^i(G^*,F^{**}) \ = \ \nsp \ \ \ \ {\rm si} \ \ i<\dim(X) . \]

Remarquons que ces conditions impliquent que $T$ est non trivial et qu'on a des
isomorphismes canoniques
\[\Ext^1(F,G^*)\simeq\Hom(G^{**},\wT) , \ \ \ \
\Ext^1(G,F^*)\simeq\Hom(F^{**},T) .\]
Le morphisme \m{G^{**}\to\wT} correspondant \`a \m{(L)} n'est autre que $\phi$,
et le morphisme \m{F^{**}\to T} correspondant \`a \m{(L')} est le morphisme
quotient. 

Les fibr\'es instables de rang 2 sur \m{\P_2} sont des extensions larges, ainsi
que certains fibr\'es instables de rang 2 sur \m{\P_3}.

Le but de cet article est de commencer l'\'etude des extensions larges, et plus
particuli\`erement des espaces $\Ext^1(\ke,\ke)$. On verra que $\Ext^1(\ke,\ke)$
poss\`ede des filtrations canoniques induites par $(L)$ et $(L^*)$, et
on \'etudiera les interactions entre ces filtrations.

\end{sub}

\sepsub

\pagebreak[3]\begin{sub}{\sc\small\small \'Etude de \m{\Ext^1(\ke,\ke)}}
\label{intro-3}\rm

Si on utilise la suite exacte $(L)$ on obtient le diagramme commutatif avec 
lignes et colonnes exactes :
\xmat{
& 0\ar[d] & 0\ar[d] & 0\ar[d] \\
0\ar[r] & \Ext^1(F,G^*)/\C\rho\ar[r]\ar[d] & M\ar[d]\ar[r] & 
\Ext^1(G^*,G^*)\ar[d]\ar[r] & 0 \\
0\ar[r] & N\ar[r]\ar[d] & \Ext^1(\ke,\ke)\ar[d]\ar[r] & 
K\ar[d]\ar[r] & 0 \\
0\ar[r] & \Ext^1(F,F)\ar[r]\ar[d] & \Ext^1(\ke,F)\ar[d]\ar[r] & 
A_2(\rho)\ar[d]\ar[r] & 0 \\
& 0 & 0 & 0 & }
avec \ \m{M=\Ext^1(\ke,G^*)}, \m{N_*=\Ext^1(F,\ke)}, $K$ \'etant un sous-espace
vectoriel de \m{\Ext^1(G^*,\ke)},
$A_2(\rho)$ d\'esignant le noyau de la multiplication par $\rho$
$$\Ext^1(G^*,F)\lra\Ext^2(F,F) .$$
On obtient de m\^eme en utilisant $(L^*)$ le diagramme commutatif avec 
lignes et colonnes \hbox{exactes :}
\xmat{
& 0\ar[d] & 0\ar[d] & 0\ar[d] \\
0\ar[r] & \Ext^1(G,F^*)/\C\pi\ar[r]\ar[d] & M_*\ar[d]\ar[r] & 
\Ext^1(F^*,F^*)\ar[d]\ar[r] & 0 \\
0\ar[r] & N\ar[r]\ar[d] & \Ext^1(\ke,\ke)\ar[d]\ar[r] & 
K_*\ar[d]\ar[r] & 0 \\
0\ar[r] & \Ext^1(G,G)\ar[r]\ar[d] & \Ext^1(\ke^*,G)\ar[d]\ar[r] & 
A_2(\pi)\ar[d]\ar[r] & 0 \\
& 0 & 0 & 0 & }
avec \ \m{M_*=\Ext^1(\ke^*,F^*)}, \m{N=\Ext^1(G,\ke^*)}, \m{K_*} \'etant un 
sous-espace vectoriel de \m{\Ext^1(F^*,\ke^*)},
$A_2(\pi)$ d\'esignant le noyau de la multiplication par $\pi$
$$\Ext^1(F^*,G)\lra\Ext^2(G,G) .$$
Remarquons qu'on a \m{M\subset N} et \m{M_*\subset N_*}. 
Les deux filtrations \ \m{0\subset M\subset N\subset\Ext^1(\ke,\ke)}, 
\m{0\subset M_*\subset N_*\subset\Ext^1(\ke,\ke)} \ et leur interaction
sont utilis\'ees dans \ref{Larg00} pour \'etudier \m{\Ext^1(\ke,\ke)}. 
On montre dans \ref{s_ext_6} que le sous-espace vectoriel
\[{\bf T} \ = \ (N\cap N_*)+M+M_* \ = \ (N+M_*)\cap(N_*+M)\]
de \m{\Ext^1(\ke,\ke)} est le {\em tangent \`a l'espace des extensions},
c'est-\`a-dire qu'il correspond aux d\'eformations de $\ke$  obtenues en
d\'eformant \m{F^{**}}, \m{G^{**}}, $T$, $\pi$ et $\rho$.

Si $X$ est une surface on a toujours \ \m{{\bf T}\not=\Ext^1(\ke,\ke)}, mais
si \m{\dim(X)>2} on peut avoir l'\'egalit\'e (cf. \ref{intro-5}, \ref{intro-6}).

La double d\'efinition de $\ke$, au moyen des suites exactes \m{(L)} et 
\m{(L^*)}, produit deux d\'efinitions de certains morphismes canoniques, et il
n'est pas \'evident que ces deux d\'efinitions co\"\i ncident (c'est m\^eme
parfois faux).
Par exemple on d\'efinit dans \ref{Larg00} deux morphismes canoniques
$$\xi_0, \xi_0^* : \End(T)\lra\Ext^1(\ke,\ke)$$
et on montre (proposition \ref{id1}) que \ \m{\xi_0=\xi_0^*}. De m\^eme on
a deux morphismes canoniques
$$\xi_2, \xi_2^* : \Ext^1(\ke,\ke)\lra\Ext^2(T,T)$$
mais ici on a \ \m{\xi_2=-\xi_2^*} (cf. proposition \ref{id4b}).

On d\'ecrit et on \'etudie en \ref{s_ext4} plusieurs sous-espaces vectoriels
canoniques de \m{\Ext^1(\ke,\ke)} :
\begin{enumerate}
\item[--] Le sous-espace vectoriel $Q$ correspondant aux d\'eformations de $\ke$
obtenues en d\'eformant $G^{**}$, $\pi$ et $\rho$, et le sous-espace vectoriel
analogue $Q_*$ correspondant aux d\'eformations de $\ke$ obtenues en 
d\'eformant $F^{**}$, $\pi$ et $\rho$.
\item[--] Le sous-espace vectoriel \ $Q+M_*=Q_*+M=M+M_*$ \ correspondant aux
d\'eformations de $\ke$ obtenues en d\'eformant $F^{**}$, $G^{**}$, $\pi$ et 
$\rho$.
\item[--] Le sous-espace vectoriel correspondant aux d\'eformations de $\ke$
obtenues en d\'eformant $T$, $\pi$ et $\rho$. On montre que c'est $N\cap N_*$.
\end{enumerate}
\end{sub}

\sepsub

\pagebreak[3]\begin{sub}{\sc\small\small Vari\'et\'es de modules d'extensions
larges}\label{intro-4}\rm

Soient $\kx$, $\ky$, $\kz$ des ensembles ouverts de faisceaux coh\'erents sur 
$X$, admettant des vari\'et\'es de modules fins {\bf M}, {\bf N}, {\bf Z} 
respectivement (cf. \ref{modfin}), les faisceaux de $\kz$ \'etant parfaits de
codimension 2. On suppose que les extensions non triviales de type \m{(L)}, avec
\m{F^{**}} dans $\kx$, \m{G^{**}} dans $\ky$, $T$ dans $\kz$ et $\ke$ localement
libre, sont des extensions larges. Soit \m{{\rm Larg}(\kx,\ky,\kz)} l'ensemble
des classes d'isomorphisme des fibr\'es $\ke$ obtenus.
On peut construire sous certaines hypoth\`eses une {\em 
vari\'et\'e de modules} \m{\M(\kx,\ky,\kz)}
pour les extensions larges de \m{{\rm Larg(\kx,\ky,\kz)}}. On obtient m\^eme
un {\em fibr\'e universel} sur \ \m{\M(\kx,\ky,\kz)\times X} (d\'efini
localement, cf. \ref{modfin}). Pour obtenir ces r\'esultats on munit 
\m{{\rm Larg(\kx,\ky,\kz)}} d'une structure de vari\'et\'e alg\'ebrique
int\`egre naturelle. Au point $q$ de \m{\M(\kx,\ky,\kz)} correspondant \`a
l'extension large $\ke$, le morphisme de d\'eformation  infinit\'esimale de 
Koda\"\i ra-Spencer
\[\omega_q : T\M(\kx,\ky,\kz)_q\lra\Ext^1(\ke,\ke)\]
(cf. \ref{K-S}) est injectif et son image est le sous-espace vectoriel $\bf T$
d\'efini pr\'ec\'edemment.

On d\'efinit en \ref{vmod01} les {\em familles
pures} d'extensions larges du type pr\'ec\'edent (cette d\'efinition est une
g\'en\'eralisation de celle de \cite{st}), et la vari\'et\'e \m{\M(\kx,\ky,\kz)}
repr\'esente le foncteur $F$ de la cat\'egorie des vari\'et\'es alg\'ebriques
complexes dans celle des ensembles d\'efini par
\[F(S) \ = \ \biggl\lbrace\text{familles pures de fibr\'es de }{\rm Larg}
(\kx,\ky,\kz)\text{ param\'etr\'ees par }S\biggr\rbrace\bigg/\sim\]
o\`u $\sim$ est la relation d\'equivalence suivante : deux familles pures
d'extensions larges $\kf_0$, $\kf_1$ param\'etr\'ees par $S$ sont \'equivalentes
si et seulement si il existe un recouvrement ouvert \m{(U_i)_{i\in I}} de $S$ et
 tel que pour tout \m{i\in I} on ait \
\m{\kf_{1\mid U_i\times X}\simeq\kf_{0\mid U_i\times X}}.

Si $X$ est une surface \m{\M(\kx,\ky,\kz)} n'est pas une vari\'et\'e de modules
fins (c'est-\`a-dire que les fibr\'es universels locaux ne sont pas des
d\'eformations compl\`etes). Dans ce cas on a en effet \ \m{{\bf
T}\not=\Ext^1(\ke,\ke)}. Cette in\'egalit\'e peut avoir plusieurs causes (cf.
\ref{intro-5}).

Si \m{\dim(X)>2}, \m{\M(\kx,\ky,\kz)} est une vari\'et\'e de modules fins si les
conditions suivantes sont r\'ealis\'ees : pour tous faisceaux $\F$ de $\kx$,
$\G$ de $\ky$, $\T$ de $\kz$ on a
\[\Ext^i(\F,\T) \ = \ \Ext^i(\G,\widetilde{T}) \ = \ \nsp \quad \text{ si }
i\geq 1 ,\]
\[H^0(\F\ot\T) \ = \ H^0(\G\ot\widetilde{T}) \ = \ \nsp .\]
Dans ce cas on a en effet \ \m{{\bf T}=\Ext^1(\ke,\ke)}. On obtient ainsi de
nouvelles vari\'et\'es de modules fins constitu\'ees de fibr\'es vectoriels non
simples (cf. \cite{dr}). Par exemple supposons que \m{X=\P_3}. Soit \m{n>4} un
entier. On prend pour $\bf Z$ la 
grassmannienne des droites de \m{\P_3} ($\kz$ est constitu\'e des faisceaux 
structuraux de ces droites),  $\bf M$ la vari\'et\'e de modules des fibr\'es de
corr\'elation nulle ($\kx$ est constitu\'e des \m{E(n)}, $E$ \'etant un fibr\'e
de corr\'elation nulle) et pour $\bf N$ la vari\'et\'e r\'eduite \`a un point
correspondant au fibr\'e \m{\ko_{\P_3}}. On obtient des extensions larges du
type
\[0\lra E(n)\lra\ke\lra\ki_\ell\lra 0 ,\]
o\`u $E$ est un fibr\'e de corr\'elation nulle, $\ell$ une droite de \m{\P_3} et
\m{\ki_\ell} son faisceau d'id\'eaux. Le fibr\'e $\ke$ est de rang 3 et de
classes de Chern \m{2n}, \m{n^2+2}, \m{2n+2}. Dans ce cas la vari\'et\'e
\m{\M(\kx,\ky,\kz)} est lisse et c'est une vari\'et\'e de modules fins. On a
\[\dim(\M(\kx,\ky,\kz)) \ = \ \dim(\Ext^1(\ke,\ke)) \ = \ 2n+14 ,\]
\[\dim(\End(\ke))= \frac{n(n+2)(n+4)}{3}+1 ,\quad
\dim(\Ext^2(\ke,\ke))=2n+14 .\]
\end{sub}

\sepsub

\pagebreak[3]\begin{sub}{\sc\small\small Extensions larges sur les surfaces}
\label{intro-5}\rm

Supposons que $X$ soit une surface. Si on consid\`ere une extension large
\m{(L)}, on a toujours \m{{\bf T}\not=\Ext^1(\ke,\ke)}. 

Dans le cas des fibr\'es
instables de rang 2 sur $\P_2$ cela est d\^u au fait que la base d'une
d\'eformation semi-universelle de $\ke$ n'est pas r\'eduite. Dans ce cas $\bf T$
est l'espace tangent au point correspondant \`a $\ke$ de la vari\'et\'e de
modules des fibr\'es instables correspondante. Mais c'est un espace de modules
grossier uniquement pour les familles {\em pures} de fibr\'es instables. 

Dans le cas g\'en\'eral, on a une situation analogue \`a la pr\'ec\'edente dans
le cas o\`u \m{\mu(G^*)\gg\mu(F)}. 
On suppose que le groupe de Picard de $X$ est isomorphe \`a 
$\Z$, de g\'en\'erateur ample \m{\ko_X(1)}.
Pour \m{i=0,1} soit \m{{\bf M}_i} une vari\'et\'e de modules des faisceaux
semi-stables non vide sur $X$. On s'int\'eresse \`a des extensions larges du 
type
$$0\lra G^*(d)\lra\ke\lra F\lra 0 $$
o\`u \m{G^*} (resp. $F$) est un faisceau semi-stable de \m{{\bf M}_0} (resp. de
\m{{\bf M}_1}). Il existe toujours de telles extensions larges si \m{d\gg 0}. 
On d\'emontre en \ref{def_larg_0} le

\sepprop

{\bf Th\'eor\`eme : }{\em Si \m{d\gg 0} les d\'eformations de 
$\ke$ sont des extensions larges du m\^eme type.}

\sepprop

Il en d\'ecoule que si \m{d\gg0}, et si \m{(S,s_0)} est la base d'une
d\'eformation semi-universelle de $\ke$, alors \m{\bf T} est l'espace tangent
en \m{s_0} \`a la sous-vari\'et\'e r\'eduite induite de $S$.

On donne cependant en \ref{ex_larg_3} des exemple d'extensions larges \m{(L)}
sur \m{\P_2}, o\`u $\ke$ est instable mais {\em prioritaire} (c'est-\`a-dire que
\m{\Ext^2(\ke,\ke(-1))=\nsp}, cf. \cite{hi_la}), donc 2-lisse, et se d\'eforme 
en fibr\'es stables.

On calcule en \ref{Ext3} le produit canonique
\[\Ext^1(\ke,\ke)\times\Ext^1(\ke,\ke)\lra\Ext^2(\ke,\ke) .\]
On peut en d\'eduire des informations sur le {\em module formel} (cf.
\ref{Mod_Extens}) $A$ de $\ke$. On peut sous certaines hypoth\`eses obtenir une
description compl\`ete de \m{A/m_A^3} (prop. \ref{Ext3_4}) comme dans 
\cite{st2} pour les fibr\'es instables de rang 2 sur \m{\P_2}. 
\end{sub}

\sepsub

\pagebreak[3]\begin{sub}{\sc\small\small Extensions larges sur les vari\'et\'es
de dimension sup\'erieure \`a 2}
\label{intro-6}\rm

On suppose que \m{\dim(X)>2}. Contrairement \`a ce qui se passe sur les surfaces
il se peut que les faisceaux r\'eguliers ne se d\'eforment pas en faisceaux
localement libres. Soit $F$ un faisceau r\'egulier sur $X$ et \m{T=F^{**}/F}.
On montre dans \ref{reg_big} que si 
\[\Ext^1(F^{**},T)=\Ext^2(F^{**},T)=\Hom(F^{**},\wT)=\nsp\]
toutes les d\'eformations de $F$ s'obtiennent en d\'eformant \m{F^{**}}, $T$ et
le morphisme surjectif \m{F^{**}\to T}. De plus $F$ est lisse, mais pas 2-lisse
en g\'en\'eral.

Soit maintenant une extension large $(L)$ telle que
\[\Ext^1(F^{**},T)=\Ext^2(F^{**},T)=\Hom(F^{**},\wT)=\nsp ,\]
\[\Ext^1(G^{**},\wT)=\Ext^2(G^{**},\wT)=\Hom(G^{**},T)=\nsp .\]
On montre dans \ref{s_ext_6} que les d\'eformations de $\ke$ sont encore des
extensions larges du m\^eme type, c'est-\`a-dire s'obtiennent en d\'eformant
\m{F^{**}}, \m{G^{**}}, $T$ et les morphismes surjectifs \m{F^{**}\to T},
\m{G^{**}\to\wT}. Si $T$ est lisse alors on montre que $\ke$ l'est aussi.
\end{sub}

\sepsub

\pagebreak[3]\begin{sub}{\sc\small\small Questions}
\label{intro-7}\rm

Les extensions larges sur les vari\'et\'es de dimension sup\'erieure \`a 2 sont
relativement simples \`a \'etudier, on trouve de nombreux cas o\`u les fibr\'es
obtenus sont lisses et o\`u leurs d\'eformations sont des extensions larges du
m\^eme type. Il faut bien entendu faire des hypoth\`eses sur les constituants de
ces extensions larges, en particulier sur les faisceaux parfaits de codimension 
2. Mais cela est in\'evitable.

Le situation est beaucoup plus compliqu\'ee sur les surfaces. Consid\'erons
l'extension large $(L)$ sur une surface. On d\'etermine au chapitre 9 le produit
canonique
\[\mu_0 : \Ext^1(\ke,\ke)\times\Ext^1(\ke,\ke)\lra\Ext^2(\ke,\ke)\]
et le lieu \m{\sigma\in\Ext^1(\ke,\ke)} tels que \m{\mu_0(\sigma,\sigma)=0} est
une r\'eunion finie de sous-espaces vectoriels contenant $\bf T$. Rappelons
qu'on a \m{\mu_0(\sigma,\sigma)=0} si et seulement si la d\'eformation
infinit\'esimale double de $\ke$ d\'efinie par $\sigma$ s'\'etend en une
d\'eformation triple (cf. \ref{triple}). On peut aussi \'etudier les
d\'eformations $n$-tuples de $\ke$. Ce sont par d\'efinition les fibr\'es
vectoriels $\kf$ sur \ \m{X\times{\rm Spec}(\C[t]/(t^n))} \ dont la restriction
\`a \ \m{X\times \lbrace*\rbrace} \ ($*$ d\'esignant le point ferm\'e de 
\m{{\rm Spec}(\C[t]/(t^n))}) est isomorphe \`a $\ke$.

\sepprop

{\bf Question : }{\em Quel est le lieu $V_n$ des \m{\sigma\in\Ext^1(\ke,\ke)} 
tels que la d\'eformation infinit\'esimale double de $\ke$ d\'efinie par 
$\sigma$ s'\'etende en une d\'eformation $n$-tuple ?}
\end{sub}

\sepprop

Les d\'eformations $n$-tuples de $\ke$ correspondent aux {\em extensions 
multiples d'ordre} $n$ de $\ke$. Ce sont des fibr\'es vectoriels $\kf$ sur $X$
munis d'une filtration
\[0=\kf_0\subset\kf_1\subset\cdots\subset\kf_n=\kf \]
telle que \ \m{\kf_i/\kf_{i-1}\simeq\ke} \ pour \m{1\leq i\leq n}, et
poss\`edant des propri\'et\'es suppl\'ementaires.
Le cas le plus simple est celui des d\'eformations triples. Soit \m{\sigma\in
\Ext^1(\ke,\ke)} tel que \ \m{\mu_0(\sigma,\sigma)=0}. Soit
\[(*) \quad\quad 0\lra\ke\lra\ke_\sigma\lra\ke\lra 0\]
l'extension correspondante. Une d\'eformation triple de $\ke$ \'etendant la
d\'eformation double correspondant \`a $\sigma$ correspond \`a un fibr\'e
vectoriel $\kf$ s'ins\'erant dans un diagramme commutatif avec lignes et
colonnes exactes
\xmat{
 & & 0\ar[d] & 0\ar[d] \\
0\ar[r] & \ke\ar[r]\fleq[d] & \ke_\sigma\ar[r]\ar[d] & \ke\ar[r]\ar[d] & 0 \\
0\ar[r] & \ke\ar[r] & \kf\ar[r]\ar[d] & \ke_\sigma\ar[r]\ar[d] & 0 \\
 & & \ke\ar[d]\fleq[r] & \ke\ar[d]\\
 & & 0 & 0
}
o\`u la suite exacte horizontale du haut et la suite exacte verticale de droite
sont identiques \`a \m{(*)}. Soient $\alpha$ l'\'el\'ement de 
\m{\Ext^1(\ke,\ke_\sigma)} correspondant \`a la suite exacte verticale du milieu
et $\beta$ l'\'el\'ement de \m{\Ext^1(\ke_\sigma,\ke)} correspondant \`a la 
suite exacte horizontale du milieu. Alors la d\'eformation triple de $\ke$
correspondant \`a $\kf$ s'\'etend en une d\'eformation quadruple si et seulement
si le produit de $\alpha$ et $\beta$ dans \m{\Ext^2(\ke,\ke)} est nul. On peut
donc ramener le probl\`eme de l'extension des d\'eformations triples en
extensions quadruples \`a l'\'etude des espaces vectoriels 
\m{\Ext^1(\ke,\ke_\sigma)}, \m{\Ext^1(\ke_\sigma,\ke)} et du produit
\[\Ext^1(\ke,\ke_\sigma)\times\Ext^1(\ke_\sigma,\ke)\lra\Ext^2(\ke,\ke) .\]

Dans le cas o\`u les hypoth\`eses du th\'eor\`eme 1.5 sont v\'erifi\'ees on
peut s'attendre \`a ce que \ \m{\cap_{n\geq 2}V_n={\bf T}}. 

\sepprop

{\bf Question : }{\em Quel est le plus petit entier $n$ tel que \m{V_n=V_{n+p}}
pour \m{p>0} ?}

\sepsub

\pagebreak[3]\begin{sub}{\sc\small\small Plan des chapitres suivants}
\label{intro-8}\rm

{\em Chapitre 2}

Le chapitre 2 est consacr\'e \`a des rappels de notations et de notions
utilis\'ees dans cet article. 

On donne en \ref{trace} la d\'efinition et quelques propri\'et\'es du
morphisme trace. 

On rappelle en \ref{Ext_cst} la construction des Ext de
faisceaux coh\'erents au moyen de r\'esolutions localement libres, qui
est utilis\'ee au chapitre 7.

On donne \ref{polHN} la d\'efinition et les propri\'et\'es du polyg\^one de
Harder-Narasimhan d'un faisceau coh\'erent sur une surface et au \ref{modfin} 
celles des vari\'et\'es de modules fins de faisceaux coh\'erents. 

Au \ref{grat} on d\'efinit les faisceaux parfaits, apr\`es quelques rappels sur
les faisceaux purs et r\'eflexifs. 

On donne en  \ref{except} la d\'efinition des
fibr\'es exceptionnels et de la fonction d'existence des fibr\'es stables sur
\m{\P_2}. Ces r\'esultats sont utilis\'es en \ref{ex_larg_3} pour produire des 
exemples d'extensions larges se d\'eformant en fibr\'es stables.

En \ref{quotalg} on \'etudie des quotients alg\'ebriques par un 
groupe r\'eductif qui seront utilis\'es en \ref{vmod} pour construire les
vari\'et\'es de modules d'extensions larges sur les surfaces. 

Au \ref{concav}
on donne un r\'esultat sur les fonctions r\'eelles concaves qu'on utilise pour
prouver le th\'eor\`eme \'enonc\'e en \ref{intro-5}.

\sepx

{\em Chapitre 3}

Le chapitre 3 est consacr\'e aux d\'eformations des faisceaux coh\'erents sur
une vari\'et\'e alg\'ebrique projective lisse. On rappelle en particulier ce 
que sont une d\'eformation semi-universelle et le module formel d'un faisceau
coh\'erent, ainsi le morphisme de d\'eformation infinit\'esimale de 
Koda\"\i ra-Spencer d'une famille de faisceaux. On termine le chapitre par une
\'etude des d\'eformations d'une extension de faisceaux coh\'erents.

\sepx

\pagebreak[2]{\em Chapitre 4}

Le chapitre 4 est une \'etude d\'etaill\'ee des extensions de faisceaux 
coh\'erents. Je pense que certains r\'esultats (comme ceux sur les morphismes
d'extensions ou les extensions duales) sont bien connus, mais j'ai inclus les
\'enonc\'es pr\'ecis et les d\'emonstrations faute d'avoir trouv\'e une
r\'ef\'erence. Les r\'esultats de \ref{ext_const} sont utilis\'es en
\ref{Def_reg} dans l'\'etude des d\'eformations des faisceaux r\'eguliers.

\sepx

{\em Chapitre 5}

Le chapitre 5 est consacr\'e \`a l'\'etude des faisceaux r\'eguliers. On
\'etudie leurs d\'eformations en \ref{Def_reg} : si $F$ est un faisceau
r\'egulier et \m{T=F^{**}/F}, on d\'etermine le sous-espace vectoriel de
\m{\Ext^1(F,F)} correspondant aux d\'eformations obtenues en d\'eformant
\m{F^{**}}, $T$ et le morphisme surjectif \m{F^{**}\to T}. 

En \ref{reg_surf} on \'etudie les faisceaux r\'eguliers sur une surface et on 
construit dans certains cas des vari\'et\'es de modules de tels faisceaux.

En \ref{reg_big} on \'etudie les faisceaux r\'eguliers sur une vari\'et\'e de
dimension sup\'erieure \`a 2. On donne des conditions suffisantes pour que si
$F$ est un faisceau r\'egulier et \m{T=F^{**}/F}, les seules d\'eformations de
$F$ s'obtiennent  en d\'eformant \m{F^{**}}, $T$ et le morphisme surjectif 
\m{F^{**}\to T}. Ces conditions impliquent aussi la lissit\'e de $F$. Pour
cela il faut partir de bons faisceaux parfaits $T$ (ils doivent \^etre lisses).
Sur \m{\P_3} on peut obtenir de tels faisceaux en partant d'ouverts lisses de 
sch\'emas de Hilbert de courbes lisses (il en existe beaucoup m\^eme si
d'apr\`es \cite{ma_pe} les sch\'emas de Hilbert de courbes ne sont pas 
en g\'en\'eral r\'eduits). On consid\`ere ensuite la jacobienne relative au 
dessus de ces 
ouverts, qui est une sorte de vari\'et\'e de modules de faisceaux parfaits.
En g\'en\'eral cependant les d\'eformations de $F$ porteront la trace des
pathologies des d\'eformations de $T$.

\sepx

{\em Chapitre 6}

Dans le chapitre 6 on d\'efinit les extensions larges, et plusieurs exemples en
sont donn\'es en \ref{Larg_exemp}. On donne en \ref{Larg_constr} une
description des extensions larges utilisant des r\'esolutions localement libres
de faisceaux qui est utilis\'ee dans le chapitre 7.

\sepx

{\em Chapitre 7}

Le chapitre 7 est consacr\'e \`a l'\'etude de \m{\Ext^1(\ke,\ke)}, o\`u $\ke$
est une extension large correspondant \`a $(L)$. 

Dans \ref{end_inc} on d\'eduit de \m{(L)} un morphisme canonique \ 
\m{\xi_0:\End(T)\to
\Ext^1(\ke,\ke)} \ et on mon-\break
tre qu'il est \'egal au morphisme analogue provenant
de \m{(L^*)} (compte tenu de l'identifi-\break
cation \m{\End(T)=\End(\wT)}). 
De m\^eme dans \ref{on_ext} on d\'eduit de \m{(L)} un morphisme canoni-\break 
que 
\m{\xi_2:\Ext^1(\ke,\ke)\to\Ext^2(T,T)} \ et on montre que le morphisme
analogue \hfil\break
\m{\xi_2^*:\Ext^1(\ke^*,\ke^*)\to\Ext^2(\wT,\wT)} d\'eduit de \m{(L^*)}
est \'egal \`a \m{-\xi_2}.

Dans \ref{s_ext4b} et \ref{s_ext_5}
on consid\`ere d'autres diagrammes canoniques, certains
commutatifs, et d'autres anticommutatifs. La plupart des d\'emonstrations sont
omises, elles sont analogues \`a celles de \ref{end_inc} et \ref{on_ext}. 

Les diagrammes commutatifs \'evoqu\'es en \ref{intro-3} sont construits en
\ref{s_ext_0}. On les utilise, ansi que certains r\'esultats des sections
pr\'ec\'edentes pour \'etudier dans \ref{s_ext4}
des sous-espaces vectoriels canoniques de
\m{\Ext^1(\ke,\ke)}, comme ceux \'evoqu\'es \`a la fin de \ref{intro-3}.

On \'etudie en \ref{s_ext_6} le sous-espace vectoriel $\bf T$ de 
\m{\Ext^1(\ke,\ke)}. On montre en particulier que si \m{\dim(X)>2} on peut
obtenir \ \m{{\bf T}=\Ext^1(\ke,\ke)} \ faisant des hypoth\`eses convenables.

On d\'etermine en \ref{s_ext_4} les produits canoniques
\[\mu_G : \End(\ke)\times\Ext^1(\ke,\ke)\lra\Ext^1(\ke,\ke) , \quad
\mu_D : \Ext^1(\ke,\ke)\times\End(\ke)\lra\Ext^1(\ke,\ke) .\]
En particulier il en d\'ecoule que \m{\Aut(\ke)} agit trivialement sur $\bf T$.

\sepx

{\em Chapitre 8}

Le chapitre 8 est consacr\'e \` la constrution des vari\'et\'es de modules
d'extensions larges.

On d\'efinit en \ref{vmod01} les familles {\em pures}
d'extensions larges. Les vari\'et\'es de modules d'extensions larges sont
en fait des espaces de modules grossiers pour les familles pures d'extensions
larges.

On montre en \ref{mod_univ} qu'il existe des {\em fibr\'es universels d\'efinis
localement} pour les vari\'et\'es de modules d'extensions larges. Quand
\m{\dim(X)>2} il peut m\^eme se produire que ces derni\`eres soient des
vari\'et\'es de modules fins. On en donne deux exemples en \ref{mod_ex_p3}
sur \m{\P_3}.

\sepx

{\em Chapitre 9}

Le chapitre 9 est consacr\'e aux extensions larges sur les surfaces. 

On d\'etermine en \ref{Ext3} le produit canonique
\[\Ext^1(\ke,\ke)\times\Ext^1(\ke,\ke)\lra\Ext^2(\ke,\ke) .\]

Le th\'eor\`eme \'enonc\'e en \ref{intro-5} est d\'emontr\'e en
\ref{def_larg_0}.

\end{sub}

\sepsec


\pagebreak[4]\section{Pr\'eliminaires}

\medskip

\Ssect{Notations}{notations}

\pagebreak[2]\begin{subsub}{Produits dans les }\rm Ext.
Soient $E$, $F$, $G$ des faisceaux coh\'erents sur une vari\'et\'e 
alg\'ebrique. Soient $p$, $q$ des entiers, \m{u\in\Ext^p(E,F)},
\m{v\in\Ext^q(F,G)}. On notera
\[l^{p,q}_u : \Ext^q(F,G)\lra\Ext^{p+q}(E,G), \ \ \ \
r^{p,q}_v : \Ext^q(E,F)\lra\Ext^{p+q}(E,G)\]
les multiplications par $u$, $v$ respectivement. Si aucune confusion n'est
\`a craindre on notera \ $l_u=l^{p,q}_u$, $r_v=r^{p,q}_v$.
\end{subsub}

\sepprop

\pagebreak[2]\begin{subsub}{Invariants logarithmiques. }\rm
Soient $X$ une vari\'et\'e alg\'ebrique projective lisse et irr\'eductible
de dimension \m{n>0}, 
\m{\ko_X(1)} un fibr\'e en droites tr\`es ample sur $X$. Soit $E$ un faisceau
coh\'erent sur $X$, tel que \m{rg(E)>0}. Le nombre rationnel
$$\mu_0(E) \ = \ \frac{c_1(E)}{rg(E)} \ \in A^1(X)\ot\Q $$
s'appelle la {\em pente} de $E$ (relativement \`a \m{\ko_X(1)}), et
$$\Delta(E) \ = \ \frac{1}{r}(c_2-\frac{r-1}{2r}c_1(E)^2) \ \in A^2(X)\ot\Q $$
est le {\em discriminant} de $E$. Plus g\'en\'eralement on d\'efinit 
formellement les {\em invariants logarithmiques} 
\ \m{\Delta_i(E)\in A^i(E)\ot\Q} \ de $E$ par
\[{\rm log}(ch(E)) \ = \ {\rm log}(r)+
\sigg_{1\leq i\leq n}(-1)^{i+1}\Delta_i(E) \]
\m{ch(E)} d\'esignant le caract\`ere de Chern de $E$ (cf. \cite{dr3}). On a
\ \m{\Delta_1(E)=\mu(E)}, \m{\Delta_2(E)=\Delta(E)}, et
\[\Delta_3(E) \ = \ \frac{1}{r}(\frac{c_3(E)}{2} + c_1(E)c_2(E)(\frac{1}{r} - 
\frac{1}{2})
+ c_1(E)^3(\frac{1}{3r^2} - \frac{1}{2r} + \frac{1}{6})) .\]
Ces invariants ont les propri\'et\'es suivantes :
\begin{enumerate}
\item[(i)] Si $L$ est un fibr\'e en droites sur $X$, on a \ $\Delta_i(E)=0$ \
si $i\geq 2$.
\item[(ii)] Si $E$, $F$ sont des faisceaux coh\'erents de rang positif sur $X$,
on a \ $\Delta_i(E\ot F)=\Delta_i(E)+\Delta_i(F)$ \ pour $1\leq i\leq n$.
\item[(iii)] Si $E$ est un faisceau localement libre non nul sur $X$, on a \
$\Delta_i(E^*)=(-1)^i\Delta_i(E)$ \ pour $1\leq i\leq n$.
\end{enumerate}
Pour tous faisceaux coh\'erents $E$, $F$ sur $X$ on note
$$\chi(E,F) \ = \ \sigg_{0\leq i\leq n}(-1)^i\dim(\Ext^i(E,F)).$$
Il existe un polyn\^ome $P_X$ \`a $n$ ind\'etermin\'ees tel que pour tout 
faisceau coh\'erent $E$ de rang positif sur $X$ on ait
\[\chi(E) \ = \ rg(E)P_X(\Delta_1(E),\ldots,\Delta_n(E)) .\]
Si $E$ et $F$ sont des faisceaux coh\'erents de rang positif sur $X$ on a
\[\chi(E,F) \ = \ rg(E)rg(F)P_X(\Delta_1(F)-\Delta_1(E),\ldots,
\Delta_n(F)-\Delta_n(E)) .
\]
\end{subsub}

\sepprop

\pagebreak[2]\begin{subsub}{Th\'eor\`eme de Riemann-Roch sur les surfaces. }\rm
On suppose que \m{\dim(X)=2}. On a alors
\[P_X(\alpha_1,\alpha_2) \ = \ \frac{\alpha_1(\alpha_1-\omega_X)}{2}+
\chi(\ko_X)-\alpha_2 .\]
\m{\omega_X} d\'esignant le fibr\'e canonique sur $X$.
Si $E$, $F$ sont des faisceaux coh\'erents sur $X$ tels que \ \m{rg(E)>0} \ et
\ \m{rg(F)>0}, on a donc
$$\chi(E,F) \ = \  rg(E).rg(F).\bigl(\frac{(\mu(F)-\mu(E))(\mu(F)-\mu(E)-
\omega_X)}{2}+\chi(\ko_X)-\Delta(E)-\Delta(F)\bigr).$$
Si $E$ et $F$ ne sont plus n\'ecessairement de rang positif, on a
$$\chi(E,F) \ = \ -c_1(E)c_1(F) 
-rg(E)c_2(F)-rg(F)c_2(E)+rg(E)rg(F)\chi(\ko_X)$$
$$ \ \ \ \ \ \ \ \ \ + 
\frac{1}{2}(rg(F).\omega_Xc_1(E)-rg(E).\omega_Xc_1(F)+
rg(E)c_1(F)^2+rg(F)c_1(E)^2).$$
\end{subsub}

\sepprop

\pagebreak[2]\begin{subsub}{Th\'eor\`eme de Riemann-Roch sur $\P_3$. }
\label{RR_P3}\rm
On a 
\[\P_{\P_3}(\alpha_1,\alpha_2,\alpha_3) \ = \
\alpha_3-\alpha_1\alpha_2-2\alpha_2+
\frac{(\alpha_1+3)(\alpha_1+2)(\alpha_1+1)}{6} .\]
Si $E$ un faisceau coh\'erent sur \m{\P_3} de rang \m{r>0}, on a donc
\[\chi(E) \ = \ r\bigl(\Delta_3-\Delta_1\Delta_2-2\Delta_2
+\frac{(\Delta_1+3)(\Delta_1+2)(\Delta_1+1)}{6}\bigr) \]
avec \m{\Delta_i=\Delta_i(E)} pour \m{i=1,2,3}.
\end{subsub}

\end{sub}

\sepsub

\Ssect{Diagrammes 3x3}{diag3x3}

On appelle {\em diagramme 3x3} (dans une cat\'egorie ab\'elienne) un 
diagramme commutatif avec lignes et colonnes exactes du type
\xmat{
& 0\ar[d] & 0\ar[d] & 0\ar[d] \\
0\ar[r] & A\ar[r]\ar[d] & M\ar[d]\ar[r] & C\ar[d]\ar[r] & 0 \\
0\ar[r] & N\ar[r]\ar[d] & E\ar[d]\ar[r] & F\ar[d]\ar[r] & 0 \\
0\ar[r] & G\ar[r]\ar[d] & H\ar[d]\ar[r] & K\ar[d]\ar[r] & 0 \\
& 0 & 0 & 0 & }
Ce diagramme est isomorphe au diagramme
\xmat{
& 0\ar[d] & 0\ar[d] & 0\ar[d] \\
0\ar[r] & M\cap N\ar[r]\ar[d] & M\ar[d]\ar[r] & M/(M\cap N)\ar[d]\ar[r] & 0 \\
0\ar[r] & N\ar[r]\ar[d] & E\ar[d]\ar[r] & E/N\ar[d]\ar[r] & 0 \\
0\ar[r] & N/(M\cap N)\ar[r]\ar[d] & E/M\ar[d]\ar[r] & E/(M+N)\ar[d]\ar[r] & 0
 \\
& 0 & 0 & 0 & }

\end{sub}

\sepsub

\Ssect{Morphisme trace}{trace}

Soit $\ke$ un faisceau coh\'erent sur une vari\'et\'e alg\'ebrique projective
lisse connexe $X$ de dimension $n$. Pour tout entier $i$ on note
\[tr_i(\ke) : \Ext^i(\ke,\ke)\lra H^i(\ko_X)\]
le morphisme trace , ou $tr_i$ et m\^eme $tr$ s'il n'y a pas 
d'ambigu\"\i t\'e. On note\hfil\break 
\m{\Ad^i(\ke) \ = \ \ker(tr_i(\ke)).}

Pour tous entiers $p, q\leq0$, on note
\[\mu_{pq} : \Ext^p(\ke,\ke)\times\Ext^q(\ke,\ke)\lra
\Ext^{p+q}(\ke,\ke)\]
\[i_{pq} : \Ext^p(\ke,\ke)\times\Ext^q(\ke,\ke)\lra
\Ext^q(\ke,\ke)\times\Ext^p(\ke,\ke)\]
les applications canoniques. Rappelons qu'on a
\[tr_{p+q}\circ i_{qp} \ = \ (-1)^{pq}tr_{p+q}\circ i_{pq} .\]
Plus g\'en\'eralement, si $\kf$ est un autre faisceau coh\'erent sur $X$, le
diagramme canonique suivant
\xmat{
\Ext^p(\ke,\kf)\times\Ext^q(\kf,\ke)\ar[r]\ar[d] & 
\Ext^{p+q}(\ke,\ke)\ar[d]^{tr_{p+q}(\ke)} \\
\Ext^{p+q}(\kf,\kf)\ar[r]^{tr_{p+q}(\kf)} & H^{p+q}(\ko_X)
}
est commutatif si $pq$ est pair et anticommutatif si $pq$ est impair.

\sepprop

\pagebreak[2]\begin{subsub}{\bf D\'efinition : }\label{F_liss}
Soit $E$ un faisceau coh\'erent sur $X$. On dit que $E$ est {\em 2-lisse} 
si le morphisme trace
\[\Ext^2(E,E)\lra H^2(\ko_X)\]
est injectif. C'est alors un isomorphisme si $\ke$ est sans torsion.
\end{subsub}

\sepprop

\pagebreak[2]\begin{subsub}{\bf Proposition : }\label{prop_tra}
Soient $x\in X$ et $\C_x$ le faisceau gratte-ciel sur $X$ concentr\'e en $x$
et de fibre $\C$ en $x$. Alors compte tenu des isomorphismes canoniques
donn\'es par la dualit\'e de Serre
\[\Ext^n(\C_x,\C_x)\simeq\omega_{X,x}^*, \ \ \
H^n(\ko_X)\simeq H^0(\omega_X)^*\]
la transpos\'ee de \ \m{tr^n:\Ext^n(\C_x,\C_x)\to H^n(\ko_X)} \
est l'\'evaluation \ \m{H^0(\omega_X)\to\omega_{X,x}}.
\end{subsub}

\begin{proof} Cela d\'ecoule des d\'efinitions de la trace et de la 
dualit\'e de Serre (cf. \cite{dr_lp}, 1.4).
\end{proof}
\end{sub}

\sepsub

\Ssect{Construction des {\rm Ext} de faisceaux coh\'erents}{Ext_cst}

Soient $E$, $F$ des faisceaux coh\'erents sur une vari\'et\'e projective
$X$. On suppose donn\'es des faisceaux coh\'erents $F_i$, $i\in\N$ sur
$X$ et une r\'esolution de $F$
\xmat{
\cdots F_2\ar[r]^{f_2} & F_1\ar[r]^{f_1} & F_0\ar[r]^{f_0} & F\ar[r] & 0 
}
Alors on en d\'eduit pour tout $n\in\N$ une application
\[\Hom(\imm(f_n),E)\lra\Ext^n(F,E)\]
qui est surjective si on a \ \m{\Ext^{n-i}(F_i,E)=\nsp} \ pour \
\m{0\leq i < n}. Si \ \m{\Ext^{n-1-i}(F_i,E)=\nsp} \ pour \
\m{0\leq i<n-1}, le noyau de cette application est constitu\'e de l'image
de $\Hom(F_{n-1},E)$ dans $\Hom(\imm(f_n),E)$. Dans ce cas $\Ext^n(F,E)$
s'identifie \`a l'espace vectoriel quotient de l'espace des morphismes \ 
\m{F_n\to E} \ s'annulant sur \m{\imm(f_{n+1})} par l'espace de ceux
qui se factorisent par $F_{n-1}$.
On peut choisir une r\'esolution de $F$ qui est finie, localement libre, et
telle que les propri\'et\'es pr\'ec\'edentes soient v\'erifi\'ees
(pour tout $n$). 
Plus g\'en\'eralement, en choisissant une r\'esolution localement libre
ad\'equate de $E$
\xmat{
\cdots E_2\ar[r]^{e_2} & E_1\ar[r]^{e_1} & E_0\ar[r]^{e_0} & E\ar[r] & 0 
}
on peut repr\'esenter les \'el\'ements de $\Ext^n(F,E)$ par des morphismes
de complexes de degr\'e $-n$ de la r\'esolution de $F$ vers celle de $E$,
c'est-\`a-dire par des suites $(\phi_i)_{i\geq 0}$ de morphismes, avec
$\phi_i : F_{i+n}\to E_i$, telles que pour tout $i\geq 0$ on ait
\[\phi_if_{i+n+1}+(-1)^ne_{i+1}\phi_{i+1} \ = \ 0 ,\]
c'est-\`a-dire que le diagramme suivant est commutatif ou anticommutatif 
suivant la parit\'e de $n$
\xmat{
E_{i+1}\ar[rr]^-{e_{i+1}} & & E_i \\
F_{i+n+1}\ar[rr]^-{f_{i+n+1}}\ar[u]^{\phi_{i+1}} & & F_{i+n}\ar[u]^{\phi_i}
}

\end{sub}

\sepsub

\Ssect{Polygone de Harder-Narasimhan}{polHN}

Soient $X$ une surface alg\'ebrique projective irr\'eductible et lisse, et
\m{\ko_X(1)} un fibr\'e en droites tr\`es ample sur $X$. 
Soit $E$ un faisceau coh\'erent sans torsion sur $X$. Alors $E$ poss\`ede une
{\em filtration de Harder-Narasimhan}
$$0 = E_0\subset E_1\subset\cdots\subset E_n=E$$
telle que pour $1\leq i\leq n$, $E_1/E_{i-1}$ soit semi-stable, que
$$\mu(E_i/E_{i-1}) \ \geq \ \mu(E_{i+1}/E_i)$$
(si $i\geq 1$), et qu'en cas d'\'egalit\'e
$$\frac{\chi(E_i/E_{i-1})}{rg(E_i/E_{i-1})} \ > \ 
\frac{\chi(E_{i+1}/E_i)}{rg(E_{i+1}/E_i)} .$$
Cette filtration est unique et $E$ est semi-stable si et seulement si $n=1$. On
dit que $n$ est la {\em longueur} de $E$.

On consid\`ere le polygone convexe dans $\R^2$ de sommets \m{O, 
Q_1=(rg(E_1),c_1(E_1).\ko_X(1)),\ldots,}\hfil\break
\m{Q_n=(rg(E_n),c_1(E_n).\ko_X(1))}.

\bigskip

\hskip 2cm
\includegraphics{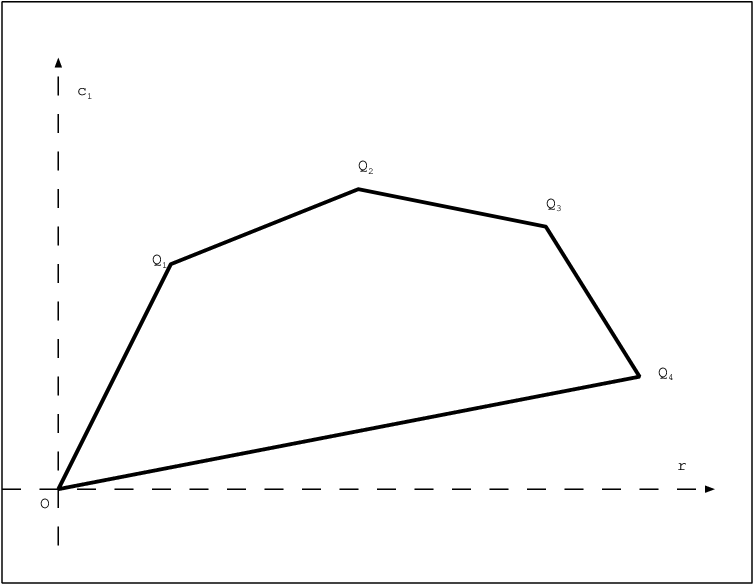}

\medskip

\centerline{{\ttx Figure 1} - Polygone de Harder-Narasimhan d'un faisceau de 
longueur 4}

\bigskip

On l'appelle le {\em polygone de Harder-Narasimhan} de $E$ et 
on le note $P(E)$.
Cette d\'efinition est l\'eg\`erement diff\'erente de celle de \cite{dr_lp}, et
s'inspire de celle de polygones de Harder-Narasimhan des fibr\'es vectoriels sur
les courbes (cf. \cite{sh}, \cite{jlp0}).

Soient $r=rg(E)$, $c_1=c_1(E).\ko_X(1)$. Soit $\kp(r,c_1)$ l'ensemble des 
polygones
convexes de sommets $O$, \m{(x_1,y_1),\ldots,(x_p,y_p)=(r,c_1)} \`a
coordonn\'ees enti\`eres et tels que \ \m{0<x_1<\cdots<x_p}. Si
\m{P\in\kp(r,c_1)} on notera aussi $P$ la fonction \m{[0,r]\to\R} associ\'ee.
Si \m{P,P'\in\kp(r,c_1)} on \'ecrira \m{P<P'} lorsque l'in\'egalit\'e est
v\'erifi\'ee par les fonctions associ\'ees, et de m\^eme pour les autres types
d'in\'egalit\'es.

\sepprop

\pagebreak[2]\begin{subsub}{\bf Proposition : }\label{polHN1}
Soient $P\in\kp(r,c_1)$, $S$  une vari\'et\'e alg\'ebrique irr\'eductible,
\m{s_0\in S}, $\kf$ un faisceau coh\'erent sur \m{S\times X}, plat sur $S$,
tel que \ \m{\kf_{s_0}\simeq E} \ et que pour \m{s\in S\backslash\lbrace
s_0\rbrace}, on ait \ \m{P(\kf_s)=P}. Alors on a \ \m{P\leq P(E)}. 
\end{subsub}

Autrement dit si $E$ est une sp\'ecialisation de faiceaux de polygone de 
Harder-Narasimhan $P$, alors on a \ \m{P\leq P(E)}.

\begin{proof}
Analogue au cas de $\P_2$ trait\'e dans \cite{dr_lp}. Voir aussi \cite{sh},
\cite{jlp0}, expos\'e 4, pour le cas des courbes.
\end{proof}

\end{sub}

\sepsub

\Ssect{Vari\'et\'es de modules fins de faisceaux coh\'erents}{modfin}

Les d\'efinitions et r\'esultats de \ref{modfin} proviennent de \cite{dr}. Soit
$X$ une vari\'et\'e projective lisse et irr\'eductible.

\sepprop

\pagebreak[2]\begin{subsub}{\bf D\'efinition : }
On appelle {\em famille de faisceaux sur $X$ param\'etr\'ee 
par $S$} un
faisceau coh\'e-\break rent $\kf$ sur \ \m{S\times X}, plat sur $S$.
\end{subsub}

\sepprop

\pagebreak[2]\begin{subsub}{\bf D\'efinition : }
On appelle {\em polyfamille de faisceaux sur $X$ 
param\'etr\'ee par $S$} la
donn\'ee d'un recouvrement ouvert \m{(U_i)_{i\in I}} de $S$, et pour tout 
\m{i\in I} d'une famille \m{\ke_i} de faisceaux sur $X$ param\'etr\'ee par
\m{U_i}, tels que pour tous \m{i,j\in I} et tout \m{x\in U_i\cap U_j} on ait
\ \m{\ke_{i,x}\simeq\ke_{j,x}}. 
\end{subsub}

\sepprop

\noindent{\em Notations : } Pour toute sous-vari\'et\'e localement 
ferm\'ee
\m{S'} de $S$, on note
$$\kf_{S'} \ = \ \kf_{\mid S'\times X}.$$
Si \m{S'} est r\'eduite \`a un point $s$ on notera plus simplement \
\m{\kf_{S'} = \kf_s}.
On note \m{p_S} (resp. \m{p_X}) la projection \ \m{S\times X\lra S} \
(resp.\ \m{S\times X\lra X}). Si $\ke$, $\kf$, $\kg$ sont des
faisceaux coh\'erents sur \ \m{S\times X}, $S$ et $X$
respectivement, on notera plus simplement
$$\ke\ot p_S^*(\kf)\ot p_X^*(\kg) \ = \ \ke\ot\kf\ot\kg.$$
Si \ \m{f : T\lra S} \ est un morphisme de vari\'et\'es alg\'ebriques, et
$\kf$ une famille de faisceaux sur $X$ param\'etr\'ee par $S$, on note
$$f^\sharp(\kf) \ = \ (f\times I_{X})^*(\kf).$$
C'est une famille de faisceaux sur $X$ param\'etr\'ee par $T$.

Si $\kx$ est un ensemble non vide de classes d'isomorphisme de faisceaux
coh\'erents sur $X$, on appelle {\em famille de faisceaux de 
$\kx$ param\'etr\'ee par $S$} une famille $\kf$ de faisceaux coh\'erents sur $X$
param\'etr\'ee par $S$ telle que pour tout point $s$ de $S$ la classe
d'isomorphisme de \m{\kf_s} soit dans $\kx$.

\medskip

Soient $r$, \m{c_i\in H^i(X,\Z)}, \m{1\leq i\leq d}, avec \ \m{r\geq 0}. 
Soit $\kx$ un ensemble non vide de classes d'isomorphisme de faisceaux
coh\'erents sur $X$, de rang $r$ et de classes de Chern \m{c_i}. 

\sepprop

\pagebreak[2]\begin{subsub}{\bf D\'efinition : }
On dit que $\kx$ est un {\em ensemble ouvert} si pour toute vari\'et\'e
alg\'ebrique $S$ et toute famille $\kf$ de faisceaux coh\'erents sur $X$
de rang $r$ et de classes de Chern \m{c_i} param\'etr\'ee par $S$,
l'ensemble des points $s$ de $S$ tels que la classe d'isomorphisme 
de \m{\kf_s} soit dans $\kx$ est un ouvert de Zariski de $S$.
\end{subsub}

\sepprop

\pagebreak[2]\begin{subsub}{\bf D\'efinition : }
Supposons que $\kx$ soit un ensemble ouvert. On dit que $\kx$ est
{\em irr\'eductible} si pour toutes polyfamilles de faisceaux de $\kx$
param\'etr\'ees par des vari\'et\'es alg\'ebriques \m{X_1}, \m{X_2}
respectivement, il existe un polyfamille de faisceaux de $\kx$ param\'etr\'ee
par une vari\'et\'e alg\'ebrique irr\'eductible $Y$ contenant tous les
faisceaux des polyfamilles param\'etr\'ees par \m{X_1} et \m{X_2}.
\end{subsub}

\sepprop

Soient $r$, \m{c_i\in H^i(X,\Z)}, \m{1\leq i\leq d}, avec \ \m{r\geq 0}. 
Soit  $\kx$ un ensemble ouvert de faisceaux coh\'erents sur $X$ de rang $r$ 
et de classes de Chern \m{c_i}. 

\sepprop
 
\pagebreak[2]\begin{subsub}{\bf D\'efinition : }
On appelle {\em vari\'et\'e de modules fins globale}, ou plus 
simplement {\em vari\'et\'e de modules fins}
pour $\kx$ la donn\'ee d'une vari\'et\'e alg\'ebrique int\`egre $M$ et 
d'une famille $\ke$ de faisceaux de $\kx$ param\'etr\'ee par $M$ telles que :

\noindent (i) Pour tout \'el\'ement $x$ de $\kx$, il existe un unique point
$m$ de $M$ tel que la classe d'isomorphisme de \m{\ke_m} soit $x$.

\noindent (ii) Pour toute famille $\kf$ de faisceaux de $\kx$ param\'etr\'ee 
par une vari\'et\'e alg\'ebrique $S$, il existe un morphisme \
\m{f:S\lra M} \ tel que pour tout point $s$ de $S$ il existe un ouvert $U$ 
de $S$ contenant $s$ et un isomorphisme \ \m{\kf_U\simeq f^\sharp(\ke)_U}.

\medskip

{\em
On dit aussi dans ce cas que \m{(M,\ke)} est une vari\'et\'e de 
modules fins
pour $\kx$. On remarquera que d'apr\`es (i), il existe une bijection canonique 
entre $\kx$ et l'ensemble des points de $M$. D'autre part le morphisme
$f$ de (ii) est unique : l'image d'un point $s$ de $S$ est le point de $M$
correspondant \`a \m{\kf_s}.}

\medskip

On appelle {\em vari\'et\'e de modules fins d\'efinie 
localement} pour $\kx$
la donn\'ee d'une vari\'et\'e alg\'ebrique int\`egre $M$, d'un recouvrement 
ouvert \m{(U_i)_{i\in I}} de $M$, et pour tout \m{i\in I} d'une famille 
\m{\ke_i} de faisceaux de $\kx$ param\'etr\'ee par \m{U_i} tels que

\noindent (i) Pour tout \'el\'ement $x$ de $\kx$, il existe un unique point
$m$ de $M$ tel que pour tout \m{i\in I} tel que \m{x\in U_i}, la classe
d'isomorphisme de \m{\ke_{im}} soit $x$.

\noindent (ii) Pour toute famille $\kf$ de faisceaux de $\kx$ param\'etr\'ee 
par une vari\'et\'e alg\'ebrique $S$, il existe un morphisme \
\m{f:S\lra M} \ tel que pour tout point $s$ de $S$ et tout \m{i\in I} tel
qu'il existe un \m{x\in U_i} tel que \m{\ke_{ix}\simeq\kf_s} il existe un 
ouvert $U$ de $S$ contenant $s$ tel que \ \m{f(U)\subset U_i} \ et un 
isomorphisme \ \m{\kf_U\simeq f^\sharp(\ke_i)_U}.
\end{subsub}
\end{sub}

\sepsub

\Ssect{Faisceaux purs, faisceaux r\'efl\'exifs et faisceaux parfaits}{grat}

(cf. \cite{lp}, 8.1 et  \cite{hu_le}, 1.1)

Soit $X$ une vari\'et\'e alg\'ebrique projective et irr\'eductible de dimension
\m{n>0}.

Soit $T$ un faisceau coh\'erent sur $X$. Rappelons qu'on appelle 
{\em dimension} de $T$ celle de son support \m{{\rm Supp}(T)}, et on la 
note \m{\dim(T)}. On appelle {\em codimension} de $E$ la codimension de 
\m{{\rm Supp}(T)}, et on la note \m{\codim(T)}.

On dit qu'un faisceau $T$ de dimension $d$ est {\em pur} si pour tout 
sous-faisceau coh\'erent non trivial $S$ de $E$ on a \ \m{\dim(S)=d}. 

Supposons $T$ pur de dimension $d$ et codimension $c=n-d$. Alors pour tout
faisceau localement libre $E$ sur $X$, les
faisceaux $\EExt^q(T,E)$ ont leur support contenu dans \m{{\rm Supp}(T)}
et sont nuls pour $q<c$. De plus on a \
\m{\codim(\EExt^q(T,E)) \ \geq \ q} \
pour $q\geq c$.

Soit $T$ un faisceau pur de dimension $d$ et de codimension $c=n-d$. On pose 
\[T^\vee \ = \ \EExt^c(T,\omega_X) , \ \ \ \ \ \
\wT \ = \ \EExt^c(T,\ko_X) = T^\vee\ot\omega_X^{-1}\]
et on appelle $T^\vee$ le {\em dual} de $T$. Il co\"\i ncide avec le dual
habituel $T^*$ de $T$ dans le cas o\`u $c=0$. 

Il existe en g\'en\'eral un morphisme canonique \ \m{T\to T^{\vee\vee}}. 

On appelle {\em faisceau r\'efl\'exif} sur $X$ un faisceau pur $T$ tel 
que le morphisme canonique \ \m{T\to T^{\vee\vee}} \ soit un isomorphisme. 
Un faisceau pur de dimension 0 est toujours r\'efl\'exif (cf. \cite{hu_le},
prop. 1.1.10).

\sepprop

\pagebreak[2]\begin{subsub}{\bf D\'efinition : }
On appelle {\em faisceau parfait} un faisceau pur $T$ sur $X$ tel que
\hfil\break \m{\EExt^q(T,\ko_X) \ = \ 0} \ pour \ \m{q>\codim(T)}. 
\end{subsub}

\sepprop

Un tel faisceau est alors r\'efl\'exif. 
Les faisceaux de dimension 0 sont parfaits, ainsi que les
faisceaux localement libres.

\sepprop

\pagebreak[2]\begin{subsub}{\bf Lemme : }\label{lemm_parf}
Soient $T$ un faisceau pur de codimension $c$ et $E$ un faisceau localement 
libre sur $X$. Alors on a \ \m{\Ext^i(T,E)=\nsp} \
si $i<c$ et des isomorphismes canoniques
\[\Ext^c(T,E) \ \simeq \ H^0(E\ot\wT), \ \ \ \ 
\Ext^{c+1}(T,E) \ \simeq \ H^1(E\ot\wT). \]
Si $T$ est parfait alors on a
\[\Ext^{c+k}(T,E) \ \simeq \ H^k(E\ot\wT)\]
si \m{k\geq 0}. En particulier, pour tout faisceau ample $\ko_X(1)$ sur $X$ on a
\[\Ext^i(T,E\ot\ko_X(p)) \ = \ \nsp \]
si \ $i>c$ \ et \ $p\gg 0$.
\end{subsub}

\begin{proof}
Cela d\'ecoule de la suite spectrale 
\[E_2^{pq} \ = \ H^p(X,\EExt^q(T,E)) \ \Longrightarrow \Ext^{p+q}(T,E)\]
(cf. \cite{go}, th\'eor\`eme 7.3.3, \cite{gr_ha}, p. 706). La derni\`ere
assertion est une cons\'equence du th\'eor\`eme B de Serre.
\end{proof}

\sepprop

\pagebreak[2]\begin{subsub}{\bf Lemme : }\label{lemm_parf2}
Soient $Y$ est une sous-vari\'et\'e ferm\'ee localement intersection compl\`ete
de $X$, \m{j:Y\to X} l'inclusion et $\kf$ un faisceau localement libre sur 
$Y$. Alors le  faisceau $T=j_*(\kf)$ est parfait, et on a
\[T^\vee \ \simeq \ j_*(\kf^*\ot\omega_Y). \]
\end{subsub}

\begin{proof}
Cela se voit facilement en utilisant le complexe de Koszul
(cf. par exemple \cite{la}, XIV, 10, ou \cite{ha}, p. 245).
\end{proof}

\sepprop

Les faisceaux localement libres sont parfaits. 
Les faisceaux purs de codimension maximale $n$ sont aussi parfaits.

Soit $T$ un faisceau parfait de codimension $c>0$, et
\xmat{
\cdots T_2\ar[r]^{t_2} & T_1\ar[r]^{t_1} & T_0\ar[r]^{t_0} & T\ar[r] & 0
}
une r\'esolution localement libre de $T$. Alors le complexe dual
\xmat{
T_0^*\ar[r]^{{}^tt_0} & T_1^*\ar[r]^{{}^tt_1} & T_2^*\ar[r]^{{}^tt_2} & \cdots
}
est exact sauf en degr\'e $c$, d'o\`u il d\'ecoule que $\ker({}^tt_{c+1})$ 
est localement libre, et on a une r\'esolution de \
\m{\wT=T^\vee\ot\omega_X^{-1}=\EExt^c(T,\ko_X)}
\xmat{
T_0^*\ar[r]^-{{}^tt_1} & T_1 
\cdots\cdots T_{c-1}^*\ar[r]^-{{}^tt_c} & \ker({}^tt_{c+1})\ar[r] & \wT\ar[r] 
& 0
}
En utilisant la construction des Ext par des r\'esolutions localement libres
(cf. \ref{Ext_cst}) on en d\'eduit ais\'ement le r\'esultat suivant : si
$T$, $T'$ sont des faisceaux parfaits de codimension $c$,
on a des isomorphismes canoniques
\[\Ext^k(T,T') \ \simeq \ \Ext^{k}(\widetilde{T'},\wT), \ \ \ 
0\leq k\leq c .\]
Si $T=T'$, ces isomorphismes sont compatibles avec la trace, c'est-\`a-dire
qu'on a un diagramme commutatif
\xmat{
\Ext^k(T,T)\ar[dd]^\simeq\ar[rrd]^{tr_k(T)} \\
& & H^k(\ko_X) \\ \Ext^k(\wT,\wT)\ar[rru]^{tr_k(\wT)}
}

On a, pour tout faisceau localement libre $E$ et tout faisceau parfait
de codimension $c$ sur $X$ un diagramme commutatif canonique
\xmat{
\Ext^c(T,E)\ot\Hom(E,T)\ar[d]^\simeq\ar[r] & \Ext^c(T,T)\ar[d]^\simeq\\
\Hom(E^*,\wT)\ot\Ext^c(\wT,E^*)\ar[r] & \Ext^c(\wT,\wT)
}

\end{sub}

\sepsub

\Ssect{Fibr\'es exceptionnels sur $\P_2$}{except}

Les conditions d'existence des faisceaux semi-stables sur $\P_2$ 
(cf \cite{dr_lp}) s'expriment en fonction des seules variables
$\mu$ et $\Delta$. On montre qu'il existe une unique fonction $\delta(\mu)$
telle qu'on ait \ \m{\dim(M(r,c_1,c_2)) > 0} \ si et seulement si \
\m{\Delta\geq\delta(\mu)}. La fonction \m{\delta(\mu)} est d\'ecrite \`a 
l'aide des {\it fibr\'es exceptionnels}.

On dit qu'un faisceau coh\'erent $\ke$ sur $\P_2$ est {\it exceptionnel} 
si $\ke$ est {\it simple} (c'est-\`a-dire si les seuls endomorphismes de 
$\ke$ sont les homoth\'eties), et si
$$\Ext^1(\ke,\ke) \ = \ \Ext^2(\ke,\ke) \ = \ \lbrace 0\rbrace.$$
Un tel faisceau est alors localement libre et stable, et la vari\'et\'e de
modules de faisceaux semi-stables correspondante contient l'unique point $\ke$.
Il existe une infinit\'e d\'enombrable de fibr\'es exceptionnels, et un
proc\'ed\'e simple permet de les obtenir tous \`a partir des fibr\'es en
droites (cf. \cite{dr1}). Notons qu'un fibr\'e exceptionnel est 
uniquement d\'etermin\'e par sa pente.
Soit $F$ un fibr\'e exceptionnel. On note \m{x_F} la
plus petite solution de l'\'equation
$$X^2-3X+\frac{1}{rg(F)^2} \ = \ 0.$$
Alors on montre que les intervalles \
\m{\rbrack\mu(F)-x_F,\mu(F)+x_F\lbrack} \  
constituent une partition de l'ensemble des nombres rationnels. On va d\'ecrire
la fonction \m{\delta(\mu)} sur cet intervalle. Posons
$$P(X) = \frac{X^2}{2}+\frac{3}{2}X+1.$$
Sur l'intervalle \ \m{\rbrack\mu(F)-x_F,\mu(F)\rbrack}, on a
$$\delta(\mu) \ = \ P(\mu-\mu(F))-\frac{1}{2}(1-\frac{1}{rg(F)^2}),$$
et sur \  \m{\lbrack\mu(F),\mu(F)+x_F\lbrack}, on a
$$\delta(\mu) \ = \ P(\mu(F)-\mu)-\frac{1}{2}(1-\frac{1}{rg(F)^2}).$$
On obtient les courbes repr\'esent\'ees sur la figure qui suit.
Ce sont des segments de coniques. La fonction $\delta(\mu)$ est la partie de ces
courbes situ\'ee au dessus de la droite d'\'equation \m{\Delta=1/2}.

\bigskip 

\hskip 2cm 
\includegraphics{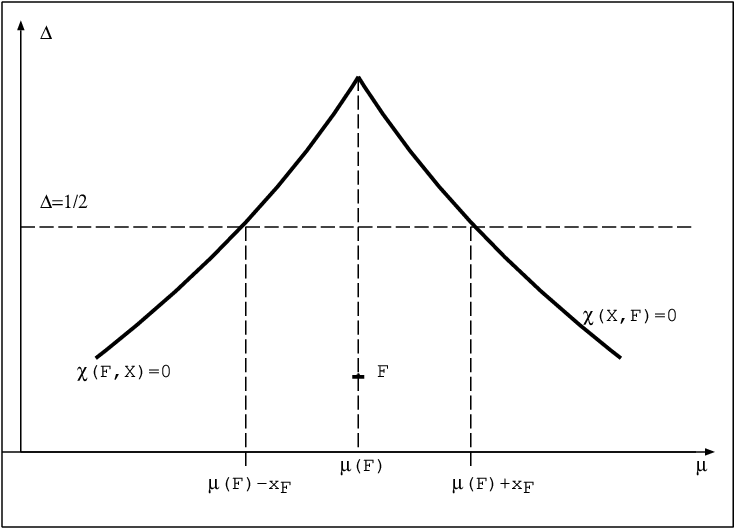}

\medskip

\centerline{{\ttx Figure 2} - La fonction $\delta(\mu)$}

\end{sub}

\sepsub

\Ssect{Quotients alg\'ebriques}{quotalg}

Soient $k$ un entier, $k\geq 1$, \m{V_1,\ldots,V_k,W_1,\ldots,W_k} des
$\C$-espaces vectoriels de dimension finie non nuls, 
$$E=V_1\times\cdots V_k\times W_1\times\cdots\times W_k, \ \ \ \
G_0=(\C^*)^{k+1} .$$
On consid\`ere l'action suivante de $G_0$ sur $E$ :
$$(t,t_1,\ldots,t_k)(v_1,\ldots,v_k,w_1,\ldots,w_k) \ = \
(t_1v_1,\ldots,t_kv_k,tt_1^{-1}w_1,\ldots,tt_k^{-1}w_k) .$$
L'action induite sur $\P(E)$ se prolonge en une action sur $\ko_{\P(E)}(1)$.
Pour que le stabilisateur d'un point g\'en\'erique soit fini, il vaut mieux
consid\'erer l'action du sous-groupe $G$ constitu\'e des \m{(t,t_1,\ldots,t_k)}
tels que \ \m{tt_1\cdots t_k=1}, qui donne les m\^emes orbites dans $\P(E)$ que 
l'action de \m{G_0}.
On peut donc d\'efinir les notions de {\em points semi-stables} et 
{\em points stables} de \m{\P(E)} (cf. \cite{mumf}, \cite{news}). Soit 
\m{\P(E)^{ss}} (resp.
\m{\P(E)^{s}}) l'ouvert des points semi-stables (resp. stables) \m{\P(E)}.
On sait qu'il existe un bon quotient \ \m{\P(E)^{ss}//G} \ et un quotient
g\'eom\'etrique \ \m{\P(E)^{s}/G} \ qui est un ouvert du pr\'ec\'edent.

\sepprop

\pagebreak[2]\begin{subsub}{\bf Proposition : }\label{quotalg1}
Un point \ \m{x=\C(v_1,\ldots,v_k,w_1,\ldots,w_k)} \ de \m{\P(E)} est stable si
et seulement si pour \m{1\leq i\leq k} on a \m{v_i\not=0}, \m{w_i\not=0}, et $x$
est semi-stable si et seulement si pour \m{1\leq i\leq k}, \m{v_i=0} si et
seulement si \m{w_i=0}.
\end{subsub}

\begin{proof}
Cela se voit ais\'ement en utilisant les crit\`eres num\'eriques de
(semi-)stabilit\'e (cf. \cite{mumf}, chapt. 2, \cite{news}, 4.2).
\end{proof}

\sepprop

On notera
$$M(V_1,\ldots,V_k,W_1,\ldots,W_k) \ = \ \P(E)^{ss}//G , \ \ \ \ \ \ 
M^s(V_1,\ldots,V_k,W_1,\ldots,W_k) \ = \ \P(E)^{s}//G .$$
La premi\`ere vari\'et\'e est projective et normale, et la seconde est un ouvert
lisse de la premi\`ere.

Soit $Y$ une vari\'et\'e alg\'ebrique int\`egre, \m{\V_1,\ldots,\V_k},
\m{\W_1,\ldots,\W_k} des fibr\'es vectoriels non nuls sur $Y$. Soit \
\m{\E=\V_1\oplus\cdots\oplus\V_k\oplus\W_1\oplus\cdots\oplus\W_k}. On
consid\`ere l'action de $G$ sur \m{\P(\E)} qui sur chaque fibre est l'action
d\'ecrite pr\'ec\'edemment. Soient \m{\P(\E)^{ss}}, \m{\P(\E)^{s}} les ouverts
de \m{\P(\E)} d\'efinis par : pour tout \m{y\in Y}, \ \m{\P(\E)^{ss}_y=
\P(\E_y)^{ss}} \ et \ \m{\P(\E)^{s}_y=\P(\E_y)^{s}}. On montre ais\'ement qu'il
existe un bon quotient
$$M(\V_1,\ldots,\V_k,\W_1,\ldots,\W_k) \ = \ \P(\E)^{ss}//G$$
et un quotient g\'eom\'etrique
$$M^s(\V_1,\ldots,\V_k,\W_1,\ldots,\W_k) \ = \ \P(\E)^{s}/G .$$
Les morphismes \m{\P(\E)^{ss}\to Y} et \m{\P(\E)^{s}\to Y} sont $G$-invariants
et passent donc au quotient, et on a, pour tout \m{y\in Y}
$$M(\V_1,\ldots,\V_k,\W_1,\ldots,\W_k)_y=
M(\V_{1y},\ldots,\V_{ky},\W_{1y},\ldots,\W_{ky}) ,$$
$$M^s(\V_1,\ldots,\V_k,\W_1,\ldots,\W_k)_y=
M^s(\V_{1y},\ldots,\V_{ky},\W_{1y},\ldots,\W_{ky}).$$
Le morphisme \ \m{M(\V_1,\ldots,\V_k,\W_1,\ldots,\W_k)\to Y} \ est projectif et
localement trivial, et  \ \m{M^s(\V_1,\ldots,\V_k,\W_1,\ldots,\W_k)\to Y} \ est
lisse et localement trivial.

On a
$$\dim(M(\V_1,\ldots,\V_k,\W_1,\ldots,\W_k)) \ = \
\dim(Y)-k-1+\sigg_{i=1}^k rg(\V_i)+\sigg_{i=1}^k rg(\W_i) .$$

\sepprop

\pagebreak[2]\begin{subsub}{\bf Remarque : }\label{quotalg2}\rm
On peut aussi consid\'erer 
\[E' \ = \ \P(V_1)\times\cdots\times\P(V_k)\times(W_1\backslash\lbrace 0\rbrace)
\times\cdots\times(W_k\backslash\lbrace 0\rbrace) ,\]
muni de l'action de \m{\C^*} :
\[ t.(l_1,\ldots,l_k,w_1,\ldots,w_k) \ = \ (l_1,\ldots,l_k,tw_1,\ldots,tw_k) .\]
Alors il existe un quotient g\'eom\'etrique \m{E'/\C^*} et on a \ 
\m{E'/\C^*\simeq\P(E)^s/G}. On peut construire de mani\`ere analogue
\m{\P(\E)^s/G}. Cette m\'ethode cache cependant les r\^oles sym\'etriques
jou\'es par les espaces vectoriels \m{V_i} et \m{W_i}.
\end{subsub}
\end{sub}

\sepsub

\newpage
\Ssect{Fonctions concaves}{concav}

\pagebreak[2]\begin{subsub}{\bf D\'efinition : }
Soient \m{a,b\in\R} tels que \m{a<b}. On dit qu'une fonction \ \m{f:[a,b]\to\R}
\ est {\em concave} si pour tous \m{x,y\in[a,b]} tels que \m{x<y} et tout 
\m{t\in[0,1]} on a
$$f(tx+(1-t)y) \ \geq \ t.f(x)+(1-t).f(y) .$$
\end{subsub}

\sepprop

Si $f$ est deux fois d\'erivable par morceaux, alors $f$ est concave si et
seulement si on a \m{f"(x)\leq 0} pour tout \m{x\in[a,b]} o\`u \m{f"} est
d\'efinie.

\sepprop

\pagebreak[2]\begin{subsub}{\bf Proposition : }\label{concav1}
Soient \ \m{f,g:[a,b]\to\R} \ des fonctions de classes \m{C^2} par morceaux et
telles que \ \m{f(a)=g(a)}, \m{f(b)=g(b)}, \m{f\leq g}. On suppose que $f$ est
concave. Alors on a
$$\int_a^bf'(x)^2dx \ \leq \ \int_a^bg'(x)^2dx .$$ 
\end{subsub}

\begin{proof}
Posons, pour \m{1\leq t\leq 1}, \m{f_t=(1-t).f+t.g}, et
$$\alpha(t) \ = \ \int^b_af'_t(x)^2dx .$$
On a 
$$\alpha(0) \ = \ \int^b_af'(x)^2dx , \ \ \ \ \ \
\alpha(1) \ = \ \int^b_ag'(x)^2dx .$$
Il suffit donc de montrer que $\alpha$ est croissante, c'est-\`a-dire que
\m{\alpha'\geq 0}. On a
$$\alpha'(t) \ = \ 2t\int_a^b(f'(x)-g'(x))^2dx +
\int_a^bf'(x)(g'(x)-f'(x))dx .$$
En int\'egrant par parties on trouve
$$\int^b_af'(x)(g'(x)-f'(x))dx \ = \ f'(g-f)\biggr]^b_a -
\int(g(x)-f(x))f"(x)dx .$$
On a \m{g-f\geq 0}, \m{f"\leq 0}, donc \ \m{\int^b_af'(x)(g'(x)-f'(x))dx\geq 0} 
\ et \m{\alpha'(t)\geq 0}.
\end{proof}
\end{sub}

\sepsec


\pagebreak[4]\section{D\'eformations infinit\'esimales de faisceaux coh\'erents}

\Ssect{D\'eformation semi-universelle d'un faisceau coh\'erent}{S_T}

Les r\'esultats de \ref{S_T} proviennent de \cite{si_tr}. Soient 
$X$ une vari\'et\'e alg\'ebrique projective et $E$ un faisceau coh\'erent
sur $X$. 
Une {\em d\'eformation} de $E$ est la donn\'ee d'un germe $(S,s_0)$ de
vari\'et\'e analytique, d'un faisceau coh\'erent $\ke$ sur \ $S\times X$, 
plat sur $S$ et d'un isomorphisme \ \m{\alpha:\ke_{s_0}\simeq E}. C'est donc 
en fait un quadruplet \ \m{\kd=(S,s_0,\ke,\alpha)}.

Consid\'erons deux d\'eformations \ \m{\kd=(S,s_0,\ke,\alpha)}, 
\m{\kd'=(S,s_0,\ke',\alpha')}
de $E$ (param\'etr\'ees par le m\^eme germe). Un {\em isomorphisme} \
\m{\kd\simeq\kd'} \ est un isomorphisme \ \m{\sigma:\ke\simeq\ke'} \ tel 
que \ \m{\alpha'\circ\sigma_{s_0}=\alpha}.
Si \ \m{\kd=f:(S',s'_0)\to(S,s_0)} \ est un morphisme de germes, on en 
d\'eduit ais\'ement la d\'eformation \ 
\m{f^*(\kd)=(S',s'_0,f^\sharp(\ke),f^\sharp(\alpha))}. 

On dit qu'une d\'eformation \ \m{\kd=(S,s_0,\ke,\alpha)} \ de $E$ est 
{\em semi-universelle} si elle est {\em compl\`ete} 
c'est-\`a-dire que si \m{\kd'} est une
d\'eformation de $E$ (param\'etr\'ee par un germe \m{(S',s'_0)}), 
il existe un morphisme de germes \ 
\m{f:(S',s'_0)\to(S,s_0)} \ et un isomorphisme \ 
\m{f^\sharp(\kd)\simeq\kd'},
et si l'application lin\'eaire tangente \ \m{T_{s'_0}S'\to T_{s_0}S} \ est
uniquement d\'etermin\'ee.
Il existe toujours une d\'eformation semi-universelle \m{(S,s_0,\ke,\alpha)}
de $E$ (\cite{si_tr}, theorem I). 

\end{sub}

\sepsub

\Ssect{D\'eformations infinit\'esimales}{def_elem}

On pose \ \m{A_2=\C[t]/(t^2)} \ et \ 
$Z_2={\rm Spec(A_2)}$.

\sepsubsub

\SSsect{Faisceaux coh\'erents sur la vari\'et\'e
double d'une vari\'et\'e alg\'ebrique}{defo0}

Si $X$ une vari\'et\'e alg\'ebrique sur $\C$, on note
$$X^{(2)} \ = \ X\times Z_2,$$
et on l'appelle la vari\'et\'e {\em double} de $X$.
Le morphisme canonique \ $A_2\to\C$ \ permet de voir $X$ comme 
une sous-vari\'et\'e ferm\'ee de
$X^{(2)}$. On note $i_X$ l'inclusion $X\subset X^{(2)}$. On a  \
\m{p_X\circ i_X=I_X}, $p_X$ d\'esignant la projection \ \m{X^{(2)}\lra X}.
Si $\ke$ est un faisceau coh\'erent sur $X^{(2)}$, on en d\'eduit deux faisceaux
coh\'erents sur $X$ : \m{\ke_{\mid X}} et
\m{\ke/t\ke} (qui est un faisceau sur \m{X^{(2)}} dont le support est contenu
dans \m{X}).

\sepprop

\pagebreak[2]\begin{subsub}{\bf Lemme : }\label{a2plat}
Un $A_2$-module $M$ est plat si et seulement si une des conditions 
\'equivalentes suivantes est r\'ealis\'ee :

(i) La
multiplication par \ \m{t:M/tM\to tM} \ est un
isomorphisme.

(ii) $M$ est libre.
\end{subsub}

\begin{proof} On utilise pour cela le fait que $M$ est plat si et seulement 
si le morphisme canonique \ 
\m{(t)\ot_{A_2}M\to tM} \ est injectif.
\end{proof}

\sepprop

\pagebreak[2]\begin{subsub}{\bf D\'efinition : }
Soit $E$ un faisceau coh\'erent sur $X$. On appelle {\em extension 
double} de $E$ une extension de $E$ par lui-m\^eme : \
\m{0\lra E\lra F\lra E\lra 0} .
\end{subsub}

\sepprop

On d\'efinit de mani\`ere \'evidente la notion d'{\em isomorphisme} 
d'extensions doubles, et on d\'emontre sans peine la

\sepprop

\pagebreak[2]\begin{subsub}{\bf Proposition : }\label{propext}
Soient $E$ un faisceau coh\'erent sur X et $\ke$ un faisceau coh\'erent sur
$X^{(2)}$ plat sur $Z_2$ et tel que \ $\ke_{\mid X}\simeq E$. Alors 
\m{p_{X*}(\ke)} poss\`ede une structure naturelle d'extension double 
de $E$ et l'association
\[ \ke\to p_{X*}(\ke) \]
d\'efinit une bijection entre l'ensemble des classes d'isomorphisme de 
faisceaux coh\'erents sur $X^{(2)}$ plats sur $Z_2$ dont la restriction \`a
$X$ est isomorphe \`a $E$, et l'ensemble des classes d'isomorphismes 
d'extensions doubles de $E$.
\end{subsub}

\sepprop

\pagebreak[2]\begin{subsub}{\bf D\'efinition : }
Soit $E$ un faisceau coh\'erent sur $X$. On appelle {\em d\'eformation
infinit\'esimale double} (ou plus simplement {\em d\'eformation
infinit\'esimale}) de $E$ un faisceau coh\'erent
$\ke$ sur $X^{(2)}$ plat sur $Z_2$ et tel que \ $\ke_{\mid X}\simeq E$. 
\end{subsub}

\sepprop

D'apr\`es la proposition \ref{propext}, les
d\'eformations infinit\'esimales de $E$ sont param\'etr\'ees
naturellement par $\Ext^1(E,E)$.

\sepsubsub

\SSsect{Morphisme de d\'eformation infinit\'esimale de Koda\"\i ra-Spencer}
{K-S}

Soient $S$ une vari\'et\'e alg\'ebrique et $\ke$ un faisceau coh\'erent sur
$S\times X$ plat sur S. Soient $s\in S$, $m_s$ l'id\'eal maximal de $s$, 
$X_s$ la fibre de $S\times X$ au dessus de $s$, $X_{s,2}$ le voisinage
infinit\'esimal d'ordre 2 de $X_s$ et $\ke_s$ la restriction de $\ke$ \`a 
$X_s$, qu'on peut voir aussi comme un faisceau coh\'erent sur $X_{s,2}$. 
Sur $X_{s,2}$ on a une suite exacte de faisceaux
$$0\lra T^*_sS\ot\ke_s\lra\ke/m_s^2\ke\lra\ke_s\lra 0$$
($T^*_sS$ d\'esignant l'espace tangent de $S$ en $s$).
En prenant l'image directe de cette suite exacte sur $X_s$ on obtient une 
suite exacte
$$0\lra T^*_sS\ot\ke_s\lra F\lra\ke_s\lra 0,$$
d'o\`u une application lin\'eaire
$$\omega_s : T_sS\lra \Ext^1_X(\ke_s,\ke_s)$$
appel\'ee {\em morphisme de d\'eformation infinit\'esimale de
Koda\"\i ra-Spencer de $\ke$ en $s$}. La relation 
avec la proposition \ref{propext} 
est la suivante : soit \ \m{\alpha\in T_sS}, \m{\alpha\not=0}, qu'on peut 
voir comme une forme $\C$-lin\'eaire sur \m{m_s/m_s^2}. On en d\'eduit un 
morphisme
$$\ov{\alpha} : Z_2\lra {\rm Spec}(\ko_s/m_s^2)$$
provenant du morphisme d'anneaux \ \m{f:\ko_s/m_s^2\lra\C[t]/(t^2)} \ 
d\'efini par
$$f(\phi) \ = \ \phi(s)+\alpha(\phi-\phi(s)).t .$$
Alors $\omega_s(\alpha)$ n'est autre que l'extension de $\ke_s$ par 
lui-m\^eme provenant de la d\'eformation 
$\ov{\alpha}^\sharp(\ke_{\mid X_{s,2}})$ d'ordre 2 de $\ke_s$.

On dit que $\ke$ est une {\em d\'eformation compl\`ete} de $\ke_s$ si 
$\omega_s$ est surjective. On dit que $\ke$ est une {\em famille 
compl\`ete} si pour tout point ferm\'e $s$ de $S$, $\ke$ est une d\'eformation 
compl\`ete de $\ke_s$.

La notion de morphisme de d\'eformation infinit\'esimale de 
Koda\"\i ra-Spencer s'\'etend aux familles plates de faisceaux analytiques
coh\'erents sur $X$ param\'etr\'ees par des vari\'et\'es analytiques ou 
des germes de vari\'et\'es analytiques. Si $E$ est un faisceau coh\'erent 
sur $X$ et \m{(S,s_0,\ke,\alpha)} une d\'eformation semi-universelle de $E$ 
(cf. \ref{S_T}), le morphisme de d\'eformation infinit\'esimale de 
Koda\"\i ra-Spencer en $s_0$ est un isomorphisme \ 
\m{T_{s_0}S\simeq\Ext^1(E,E)} \ (cf. \cite{si_tr}).

\end{sub}

\sepsub

\Ssect{D\'eformations triples}{triple}

On pose \ \m{A_3=\C[t]/(t^3)} \ et \ 
$Z_3={\rm Spec(A_3)}$. 

\sepprop

\pagebreak[2]\begin{subsub}{\bf Lemme : }
Un $A_3$-module $M$ est plat si et seulement si une des conditions 
\'equivalentes suivantes est r\'ealis\'ee :

(i) les multiplication par $t$  \m{tM/t^2M\to t^2M} \ et par
$t^2$ \m{M/tM\to t^2M} sont des isomorphismes.

(ii) $M$ est libre.
\end{subsub}

(Analogue au lemme \ref{a2plat}).

\sepprop

Soient $X$ une vari\'et\'e alg\'ebrique projective et \
\m{X^{(3)} \ = \ X\times Z_3} .
On a des immersions naturelles \ \m{X\subset X^{(2)}\subset X^{(3)}}.

\sepprop

\pagebreak[2]\begin{subsub}{\bf D\'efinition : }
Soit $E$ un faisceau coh\'erent sur $X$. On appelle {\em d\'eformation
infinit\'esimale triple} de $E$ un faisceau coh\'erent
$\ke$ sur $X^{(3)}$ plat sur $Z_3$ et tel que \ $\ke_{\mid X}\simeq E$. 
\end{subsub}

\sepprop

 Les restrictions \`a $X^{(2)}$ des
d\'eformations triples de $E$ donnent des d\'eformations doubles de $E$.
Les d\'eformations triples de $E$ correspondent aux {\em extensions 
triples} de $E$, c'est-\`a-dire aux triplets $(F,(F_i),\tau)$ constitu\'e de
\begin{enumerate}
\item[-] un faisceau coh\'erent $F$ sur $X$,
\item[-] une filtration
\[0=F_0\subset F_1\subset F_{2}\subset F_3=F\]
dont tous les gradu\'es $F_i/F_{i-1}$ sont isomorphes \`a $E$, 
\item[-] un isomorphisme $\tau : F_3/F_1\to F_2$ dont la restriction \`a 
$F_2/F_1$ induit un isomorphisme $F_2/F_1\simeq F_1$.
\end{enumerate}
On d\'emontre sans peine la

\sepprop

\pagebreak[2]\begin{subsub}{\bf Proposition : }\label{def_tripl}
Une d\'eformation double de $E$, d\'efinie par \ \m{\sigma\in\Ext^1(E,E)} \
s'\'etend en une d\'eformation triple de $E$ si et seulement si on a \
\m{\sigma^2=0} \ dans \m{\Ext^2(E,E)}.
\end{subsub}
\end{sub}

\sepsub

\Ssect{Module formel}{Mod_Extens}

Soient $E$ un faisceau coh\'erent sur $X$ et \m{(S,s_0,\ke,\alpha)} une
d\'eformation semi-universelle de $E$. L'anneau
\[R_E \ = \ \widehat{\ko}_{S,s_0}\]
s'appelle le {\em module formel} de $E$. Soient 
\[A \ = \ \widehat{\ko}_{\Ext^1(E,E),0} \]
et \m{m_A} l'id\'eal maximal de $A$. D'apr\`es \cite{lau}, il existe une
application lin\'eaire \hfil\break
 \m{\omega : \Ext^2(E,E)^*\to m_A^2}  \ telle que
 \ \m{R_E\simeq A/A\imm(\omega)}.
On a un isomorphisme canonique
\[m_A^2/m_A^3 \ \simeq \ S^2(\Ext^1(E,E)^*) .\]
Soit \ \m{\omega_2(E) : \Ext^2(E,E)^*\to S^2(\Ext^1(E,E)^*)} \ la compos\'ee de
l'application quotient \hfil\break
 \m{m_A^2\to m_A^2/m_A^3} \ et de $\omega$. On a alors

\sepprop

\pagebreak[2]\begin{subsub}{\bf Proposition : }\label{Laud1} Supposons que 
\m{X=\P_n}. Alors on a, pour tous \m{\alpha,\beta\in\Ext^1(E,E)},
\[{}^t\omega_2(E)(\alpha.\beta) \ = \
\frac{1}{2}(\mu_0(\alpha,\beta)+\mu_0(\beta,\alpha)) , \]
$\mu_0$ d\'esignant le produit canonique \
\m{\Ext^1(E,E)\times\Ext^1(E,E)\to\Ext^2(E,E)} .
\end{subsub}

(cf. \cite{lau}, \cite{st2}).

\end{sub}

\sepsub

\Ssect{D\'eformations d'extensions}{Def_Extens}

(Voir le chapitre \ref{Extens} pour la d\'efinition et des propri\'et\'es des
extensions)

Soit $X$ une vari\'et\'e alg\'ebrique projective lisse et connexe.

\sepprop

\pagebreak[2]\begin{subsub}{\bf Proposition : }\label{Def_Extens1}
Soient $E$, $F$ des faisceaux coh\'erents simples sur $X$. On 
suppose que \ \m{\Hom(E,F)=\Ext^2(F,E)=\nsp}. Soit
\xmat{ 0\ar[r] & E\ar[r] & \ke\ar[r] & F\ar[r] & 0 }
une extension, associ\'ee \`a \ \m{\sigma\in\Ext^1(F,E)}. Soient
\[ f_\sigma : \Ext^1(E,\ke)\lra\Ext^2(F,\ke), \ \ 
g_\sigma : \Ext^1(\ke,F)\lra\Ext^2(\ke,F) \]
les multiplications par $\sigma$. Alors il existe un diagramme commutatif 
canonique avec lignes et colonnes exactes
\xmat{ & 0\ar[d] & 0\ar[d] & 0\ar[d] \\
0\ar[r] & \Ext^1(F,E)/\C\sigma\ar[r]\ar[d] & \Ext^1(\ke,E)\ar[r]\ar[d] & 
\Ext^1(E,E)\ar[r]\ar[d]& 0 \\
0\ar[r] & \Ext^1(F,\ke)\ar[r]\ar[d] & \Ext^1(\ke,\ke)\ar[r]\ar[d] & 
\ker(f_\sigma)\ar[r] & 0 \\
0\ar[r] & \Ext^1(F,F)\ar[r]\ar[d] & \ker(g_\sigma)\ar[d] \\
 & 0 & 0}
\end{subsub}

\begin{proof} Imm\'ediat.\end{proof}

\sepprop

On note $D(\sigma)$ le sous-espace vectoriel \ 
\m{\Ext^1(\ke,E)+\Ext^1(F,\ke)} \ de \m{\Ext^1(\ke,\ke)}.

Soit $\E$ (resp. $\F$) une famille de faisceaux coh\'erents lisses et simples 
sur $X$ param\'etr\'ee par une vari\'et\'e alg\'ebrique affine $S$ (resp. $T$). 
On suppose que pour tous points ferm\'es $s\in S$, $t\in T$,  on a
\[ \Hom(\E_s,\F_t) \ = \ \Ext^2(\F_t,\E_s) \ = \ \nsp . \]
On suppose que la dimension de
\m{\Ext^1(\F_t,\E_s)} est ind\'ependante de $s\in S$ et $t\in T$. On 
d\'emontre alors comme dans le lemme \ref{ext_reg0} que 
\[ \kv \ = \ \EExt^1_{p_{S\times T}}(p^\sharp_S(\F),p_T^\sharp(\E))\]
est un faisceau localement libre sur \ $S\times T$ (en chaque \
\m{(s,t)\in S\times T} \ la fibre $\kv_{(s,t)}$ du fibr\'e vectoriel 
correspondant s'identifie \`a $\Ext^1(\F_t,\E_s)$). Soit
\ \m{\pi : \kv\to S\times T} \ la projection. Il existe une 
{\em extension universelle} sur \ \m{\kv\times X}
\xmat{ 0\ar[r] & \pi^\sharp(p_S^\sharp(\E))\ar[r] & \ke\ar[r] 
& \pi^\sharp(p_T^\sharp(\F))\ar[r] & 0 }
telle que pour tout $v\in\kv$ au dessus de \ $(s,t)\in S\times T$ \ la 
restriction de la suite exacte pr\'ec\'edente \`a \ 
\m{\lbrace v\rbrace\times X} \
\xmat{ 0\ar[r] & \E_s\ar[r] & \ke_v\ar[r] & \F_t\ar[r] & 0 }
soit associ\'ee \`a \ \m{v\in\Ext^1(\F_t,\E_s)}.

\sepprop

\pagebreak[2]\begin{subsub}{\bf Proposition : }\label{Def_Extens2}
Soit $v$ un point ferm\'e de $\kv$, au dessus de \ \m{(s,t)\in S\times T}. 
Alors 

1 - Le morphisme de d\'eformation infinit\'esimale de Koda\"\i ra-Spencer
\[ \omega_v : T\kv_v\lra\Ext^1(\kv_v,\kv_v) \]
est \`a valeurs dans $D(v)$.

2 - Si $\E$, $\F$ sont des d\'eformations compl\`etes de $\E_s$, $\F_t$
respectivement, l'image de $\omega_v$ est \'egale \`a $D(v)$.
\end{subsub}

\medskip

La proposition \ref{Def_Extens2} est une cons\'equence imm\'ediate des deux
lemmes qui suivent.

Soient \ \m{\kv^s=\kv_{\mid \lbrace s\rbrace\times T}},
\m{\ke^s=\ke_{\mid\kv^s}}, \m{\kv^t=\kv_{\mid S\times\lbrace t\rbrace}},
\m{\ke^t=\ke_{\mid\kv^t}}. Soient
\[ \omega_v^s : T\kv_v^s\lra\Ext^1(\kv_v,\kv_v), \ \ \ \
\omega_v^t : T\kv_v^t\lra\Ext^1(\kv_v,\kv_v) \]
les morphismes de d\'eformation infinit\'esimale de Koda\"\i ra-Spencer de
$\ke^s$, $\ke^t$ respectivement, au point $h$, qui sont des restrictions
de $\omega_v$. Soient enfin
\[ \omega_s : T\kv_v^s\lra\Ext^1(\F_t,\F_t), \ \ \ \
\omega_t : T\kv_v^t\lra\Ext^1(\E_s,\E_s) \] 
les morphismes de d\'eformation infinit\'esimale de Koda\"\i ra-Spencer de
$p_S^\sharp(\F)$, $p_T^\sharp(\E)$ respectivement, au point $v$.

\sepprop

\pagebreak[2]\begin{subsub}{\bf Lemme : }\label{Def_Extens3}
(i) L'application $\omega_v^s$ est \`a valeurs dans le sous-espace vectoriel 
\m{\Ext^1(\F_t,\ke_v)} \ de \m{\Ext^1(\ke_v,\ke_v)}.

(ii) La compos\'ee 
\xmat{
T\kv_v^s\ar[r]^-{\omega_v^s} & \Ext^1(\F_t,\ke_v)\flon[r] & \Ext^1(\F_t,\F_t) }
est \'egale \`a $\omega_s$.

(iii) L'image de $\omega_v^s$ est exactement \
\m{\Ext^1(\F_t,\ke_v)} \ si $\F$ est une d\'eformation
compl\`ete de $\F_t$.
\end{subsub}

\sepprop

\pagebreak[2]\begin{subsub}{\bf Lemme : }\label{Def_extens4}
(i) L'application $\omega_v^t$ est \`a valeurs dans le sous-espace vectoriel 
\ \m{\Ext^1(\ke_v,\E_s)} \ de \m{\Ext^1(\ke_v,\ke_v)}.

(ii)  La compos\'ee 
\xmat{
T\kv_v^t\ar[r]^-{\omega_v^t} & \Ext^1(\ke_v,\E_s)\flon[r] & \Ext^1(\E_s,\E_s) }
est \'egale \`a $\omega_t$.

(iii) L'image de $\omega_v^t$ est exactement \
\m{\Ext^1(\ke_v,\E_s)} \ si $\E$ est une d\'eformation
compl\`ete de $\E_s$.
\end{subsub}

\medskip

Les d\'emonstrations de ces deux lemmes sont analogues \`a celles des lemmes
\ref{ext_reg2} et \ref{ext_reg3}.
\end{sub}

\sepsec


\pagebreak[4]\section{Extensions de faisceaux coh\'erents}\label{Extens}\rm

\Ssect{D\'efinition}{Extensdef}

Soient $E$, $F$ des faisceaux coh\'erents sur une vari\'et\'e 
projective lisse $X$.
Rappelons qu'une {\em extension} de $F$ par $E$ est une suite exacte 
$$0\lra E\lra \ke \lra F\lra 0 ,$$
o\`u $\ke$ est un faisceau coh\'erent sur $X$. Deux extensions
\[ 0\lra E\lra \ke \lra F\lra 0 , \ \ \ \ 0\lra E\lra \ke' 
\lra F\lra 0 \]
sont dites {\em isomorphes} s'il existe un diagramme commutatif
\xmat{
0\ar[r] & E\ar@{=}[d]\ar[r] & \ke\ar[r]\ar[d] & F\ar[r]\ar@{=}[d] & 0 \\
0\ar[r] & E\ar[r] & \ke'\ar[r] & F\ar[r] & 0 \\
}
la fl\`eche verticale du milieu \'etant alors un isomorphisme. Il est bien
connu qu'il existe une bijection canonique entre l'ensemble des classes
d'isomorphisme d'extensions de $F$ par $E$ et l'espace vectoriel
\m{\Ext^1(F,E)}. On va construire explicitement cette bijection (la 
construction s'inspire de celle donn\'ee dans \cite{gr_ha}).
\end{sub}

\sepsub

\Ssect{Construction des extensions}{Extensconstr}

Si \ \m{0\lra E\lra \ke \lra F\lra 0} \ est une extension de $F$ par $E$, 
on en
d\'eduit le morphisme de liaison \ \m{\delta:\End(E)\to\Ext^1(F,E)} \ 
provenant
de l'application du foncteur $\Hom(-,E)$ \`a la suite exacte pr\'ec\'edente, 
et l'\'el\'ement
de \m{\Ext^1(F,E)} associ\'e \`a l'extension pr\'ec\'edente est 
\m{\delta(I_E)}. 

R\'eciproquement soit \ \m{\sigma\in\Ext^1(F,E)}.
Soit
\xmat{\cdots F_2\ar[r]^{f_2} & F_1\ar[r]^{f_1} & F_0\ar[r]^{f_0} &
F\ar[r] & 0 }
une r\'esolution localement libre de $F$. On en d\'eduit une suite exacte
$$\Hom(F_0,E)\ \hfl{\alpha}{} \ \Hom(F_1,E)\
\hfl{\beta}{} \ \Hom(F_2,E).$$
On peut d'apr\`es \ref{Ext_cst} choisir cette r\'esolution de
telle sorte que $\Ext^1(F,E)$ soit canoniquement isomorphe 
\`a \
\m{\ker(\beta)/\imm(\alpha)}. Donc $\sigma$ est repr\'esent\'e par un 
morphisme
\ \m{\lambda : F_1/\imm(f_2)\lra E}. Soit
\[\mu = \lambda\oplus\ov{f_1}: F_1/\imm(f_2)\lra E\oplus F_0 ,\]
($\ov{f_1}$ \'etant induit par $f_1$) qui est injectif. Soit \
\m{\ke\ = \ (E\oplus F_0)/\imm(\mu)}.
On a alors une extension
$$0\lra E\ \hfl{i}{} \ \ke \ \hfl{p}{} \ F\lra 0 ,$$
avec \ \m{i(e)=[(e,0)]} \ et \ \m{p([(u,v)])=f_0(v)}. C'est l'extension 
associ\'ee \`a $\sigma$.

\medskip

Plus g\'en\'eralement, si
$$F'_2\ \hfl{f'_2}{} \ F'_1\ \hfl{f'_1}{} \ F'_0 \ 
\hfl{f'_0}{} \ F'\lra 0$$
est une suite exacte de faisceaux coh\'erents, on en d\'eduit le complexe
$$\Hom(F'_0,E)\ \hfl{\alpha'}{} \ \Hom(F'_1,E)\
\hfl{\beta'}{} \ \Hom(F'_2,E).$$
Il existe toujours une application canonique \
\m{\ker(\beta')/\imm(\alpha')\lra\Ext^1(F,E)} \
qui n'est pas en g\'en\'eral un isomorphisme. 

L'extension induite par un morphisme \ \m{F'_1/\imm(f'_2)\to E} \ se
construit de la m\^eme mani\`ere que pr\'ec\'edemment.

\end{sub}

\sepsub

\Ssect{Morphismes d'extensions}{Extensmorph}

\pagebreak[2]\begin{subsub}{\bf Proposition : }\label{Extensprop}
Soient $E$, $F$, $F'$ des faisceaux coh\'erents sur $X$, \ 
\m{f:F\to F'} \ un
morphisme, et
\[{
\xymatrix{0\ar[r] & E\ar[r]^i & \ke\ar[r]^p & F\ar[r] & 0\\ }\ \ \ \ , \ \ \ \
\xymatrix{0\ar[r] & E\ar[r]^{i'} & \ke'\ar[r]^{p'} & F'\ar[r] & 0\\ }
}\]
des extensions, associ\'ees respectivement \`a \m{\sigma\in\Ext^1(F,E)} et
\m{\sigma'\in\Ext^1(F',E)}. Alors il existe un diagramme commutatif
\xmat{
0\ar[r] & E\ar[r]^i\ar@{=}[d] & \ke\ar[r]^p\ar[d]^\phi & F\ar[r]\ar[d]^f & 0\\
0\ar[r] & E\ar[r]^{i'} & \ke'\ar[r]^{p'} & F'\ar[r] & 0\\
}
si et seulement si on a \ \m{l_{f}(\sigma')=\sigma}. Dans ce cas $\phi$ est
injectif (resp. surjectif) si et seulement si $f$ l'est.
\end{subsub}

\begin{proof} Supposons qu'un tel diagramme commutatif existe. 
L'\'egalit\'e \ 
\m{l_{f}(\sigma')=\sigma}  \ d\'ecoule du diagramme commutatif
\xmat{
\End(E)\ar[r]^{\delta'}\ar@{=}[d] & \Ext^1(E,F')\ar[d]^{l_{f}} \\
\End(E)\ar[r]^\delta & \Ext^1(E,F)\\
}
($\delta'$ et $\delta$ \'etant les morphismes de liaison provenant des suites
exactes obtenues par l'application du foncteur $\Hom(E,-)$)
et des d\'efinitions de $\sigma$ et $\sigma'$. R\'eciproquement, 
supposons que
\ \m{l_{f}(\sigma')=\sigma}. On consid\`ere des r\'esolutions localement 
libres de $F$ et $F'$ (cf. \ref{Extensdef}) pour calculer 
$\Ext^1(F,E)$ 
et $\Ext^1(F',E)$. D'apr\`es la d\'efinition de ces r\'esolutions on peut les
choisir de telle sorte qu'on ait un diagramme commutatif
\xmat{
F_2\ar[r]^{f_2}\ar[d]^{\phi_2} & F_1\ar[r]^{f_1}\ar[d]^{\phi_1} &
F_0\ar[r]^{f_0}\ar[d]^{\phi_0} & F\ar[d]^f \\
F'_2\ar[r]^{f'_2} & F'_1\ar[r]^{f'_1} & F'_0\ar[r]^{f'_0} & F'\\
}
Supposons que $\sigma'$ soit repr\'esent\'e par un morphisme \
\m{\lambda': F'_1/\imm(f'_2)\to E}.
 Alors \hfil\break \m{\ov{f}(\sigma')} est repr\'esent\'e par \
\m{\lambda=\lambda'\circ\ov{\phi_1}: F_1/\imm(f_2)\to E} , \m{\ov{\phi_1}}
d\'esignant le morphisme \hfil\break
\m{F_1/\imm(f_2)\to F'_1/\imm(f'_2)} \ induit par $\phi_1$. Soient
\[\mu=\lambda\oplus\ov{f_1}:F_1/\imm(f_2)\lra E\oplus F_0, \ \ \ \
\mu'=\lambda'\oplus\ov{f'_1}:F'_1/\imm(f'_2)\lra E\oplus F'_0 .\]
Alors on a \
\m{\ke = (E\oplus F_0)/\imm(\mu)}, \m{\ke' = (E\oplus F'_0)/\imm(\mu')}.
Soit \ \m{\phi : \ke\to\ke'} \ le morphisme induit par \m{I_E\oplus\phi_0}. 
Il
est ais\'e de voir que le diagramme de la proposition \ref{Extensprop} est
commutatif. La derni\`ere assertion de la proposition \ref{Extensprop} est
\'evidente. 
\end{proof}

\sepprop

\pagebreak[2]\begin{subsub}{\bf Proposition : }\label{Extensprop2}
Soient $G$, $\ke$, $F$ des faisceaux coh\'erents sur $X$, \ 
\m{h:\ke\to F} \ un
morphisme, et
\[{
\xymatrix{0\ar[r] & F\ar[r]^i & F'\ar[r]^p & G\ar[r] & 0\\ }\ \ \ \ , \ \ \ \
\xymatrix{0\ar[r] & \ke\ar[r]^{i'} & \ke'\ar[r]^{p'} & G\ar[r] & 0\\ }
}\]
des extensions, associ\'ees respectivement \`a \m{\sigma\in\Ext^1(G,F)} et
\m{\sigma'\in\Ext^1(G,\ke)}. Alors il existe un diagramme commutatif
\xmat{
0\ar[d] & 0\ar[d] \\
\ke\ar[r]^h\ar[d]^{i'} & F\ar[d]^i \\
\ke'\ar[r]^{\psi}\ar[d]^{p'} & F'\ar[d]^p \\
G\ar@{=}[r]\ar[d] & G\ar[d] \\
0 & 0 \\
}
si et seulement si on a \ \m{r_{h}(\sigma')=\sigma}. Dans ce cas $\psi$ est
injectif (resp. surjectif) si et seulement si $h$ l'est.
\end{subsub}

\begin{proof} Analogue \`a la proposition \ref{Extensprop}. \end{proof} 

\sepprop

\pagebreak[2]\begin{subsub}{\bf Corollaire : }\label{Extenscoro}
Soient $A_0$, $B_0$, $C_0$, $A_1$, $B_1$, $C_1$ des faisceaux coh\'erents sur 
$X$ et
\xmat{0\ar[r] & A_0\ar[r]^{i_0} & B_0\ar[r]^{p_0} & C_0\ar[r] & 0}
\xmat{0\ar[r] & A_1\ar[r]^{i_1} & B_1\ar[r]^{p_1} & C_1\ar[r] & 0}
des suites exactes. Soient \ \m{\alpha : A_0\to A_1}, \m{\gamma : C_0\to C_1} \
des morphismes. Soit \m{\eta_0\in\Ext^1(C_0,A_0)}
(resp. \m{\eta_1\in\Ext^1(C_1,A_1)}) l'\'el\'ement associ\'e \`a la premi\`ere
(resp. seconde) suite exacte. Alors il existe un morphisme \ \m{\beta : B_0\to
B_1} \ tel que le diagramme suivant 
\xmat{
A_0\ar[r]^{i_0}\ar[d]^\alpha & B_0\ar[r]^{p_0}\ar[d]^\beta & C_0\ar[d]^\gamma\\
A_1\ar[r]_{i_1} & B_1\ar[r]_{p_1} & C_1
}
soit commutatif si et seulement si on a \ \m{r_\alpha(\eta_0) = 
l_\gamma(\eta_1)}.
\end{subsub}

\end{sub}

\sepsub

\Ssect{Extensions et diagrammes 3x3}{ext_const}

\pagebreak[2]\begin{subsub}{\bf Proposition : }\label{ext_const1}
Soient $A$, $A'$, $E$, $F$ des faisceaux coh\'erents sur $X$,
\xmat{
(S) \ \ \ \ \ \ 0\ar[r] & A\ar[r] & E\ar[r] & A'\ar[r] & 0 }
une suite exacte. 

1 - Soit \ \m{\delta : \Hom(A,F)\to\Ext^1(A',F)} \ l'application
induite par $(S)$. Soit \ \m{\phi\in\Hom(A,F)}. Alors il existe un diagramme
commutatif avec lignes et colonnes exactes
\xmat{
 & & 0\ar[d] & 0\ar[d] & \\
 & & A\fleq[r]\ar[d]^i & A\ar[d] & \\
0\ar[r] & F\fleq[d]\ar[r] & F\oplus E\ar[r]\ar[d]^\psi & E\ar[d]\ar[r] & 0 \\
0\ar[r] & F\ar[r] & U\ar[r]\ar[d] & A'\ar[r]\ar[d] & 0 \\
 & & 0 & 0 & \\
}
tel que la la suite exacte verticale de droite soit l'extension $(S)$,
que la suite exacte horizontale du bas soit associ\'ee \`a 
\m{\delta(\phi)}, que celle du milieu soit triviale, que le morphisme \
\m{A\to F} \ induit par \ \m{\psi_{\mid E} : E\to U} \ soit 
\'egal \`a $\phi$, et que celui qui provient de $i$ soit \'egal \`a $-\phi$.

R\'eciproquement, supposons donn\'e un tel diagramme, tel que la suite exacte
verticale de droite soit l'extension $(S)$, et que la suite exacte horizontale
du milieu soit triviale. Soit \ \m{\phi : A\to F} \ le morphisme induit par \
\m{\psi_{\mid E} : E\to U}. Alors le morphisme $A\to F$ d\'efini par $i$ est
\'egal \`a $-\phi$ et la suite exacte horizontale du bas est associ\'ee \`a
$\delta(\phi)$.

\medskip

2 - Soit \ \m{\delta' : \Hom(F,A')\to\Ext^1(F,A)} \ l'application
induite par $(S)$. Soit \ \m{\phi'\in\Hom(F,A')}. Alors il existe un 
diagramme commutatif avec lignes et colonnes exactes
\xmat{
 & 0\ar[d] & 0\ar[d] & & \\
0\ar[r] & A\ar[r]\ar[d] & V\ar[r]\ar[d]^{\psi'} & F\ar[r]\fleq[d] & 0 \\
0\ar[r] & E\ar[d]\ar[r] & E\oplus F\ar[r]\ar[d]^p & F\ar[r] & 0 \\
 & A'\ar[d]\fleq[r] & A'\ar[d] & & \\
 & 0 & 0 & & }
tel que la la suite exacte verticale de gauche soit l'extension $(S)$,
que la suite exacte horizontale du haut soit associ\'ee \`a 
\m{\delta'(\phi')}, que celle du milieu soit triviale, que le morphisme \
\m{\alpha' : F\to A'} \ induit par la composante \m{V\to E} de \m{\psi'} \
soit \'egal \`a $\phi'$, et que celui qui provient de $p$ soit \'egal \`a
$-\phi'$.

R\'eciproquement, supposons donn\'e un tel diagramme, tel que la suite exacte
verticale de gauche soit l'extension $(S)$, et que la suite exacte horizontale
du milieu soit triviale. Soit \ \m{\phi' : F\to A'} \ le morphisme induit par \
la composante \m{V\to E} de \m{\psi'}. Alors le morphisme $F\to A'$ d\'efini
par $p$ est \'egal \`a $-\phi'$ et la suite exacte horizontale du haut est 
associ\'ee \`a $\delta'(\phi')$.
\end{subsub}

\begin{proof}
On ne d\'emontrera que 1-, 2- \'etant analogue.
On utilise la construction explicite des extensions donn\'ee dans
 \ref{Extensdef}. On utilise une r\'esolution localement libre ad\'equate
de $A'$ :
\xmat{
A'_2\ar[r]^{f'_2} & A'_1\ar[r]^{f'_1} & A'_0\ar[r]^{f'_0} & A' }
L'extension de de $A'$ par $A$ de la proposition \ref{ext_const1} provient
d'un \'el\'ement de \m{\Ext^1(A',A)} d\'efini par un morphisme
\[ f : A'_1/\imm(f'_2)\lra A . \]
Soit \ \m{\mu=f\oplus\ov{f'_1} : A'_1/\imm(f'_2)\to A\oplus A'_0}. Alors on
a \ \m{E=(A\oplus A'_0)/\imm(\mu)}. 

On construit maintenant le diagramme commutatif de la proposition
\ref{ext_const1}. Soit
\[\eta \ = \ (\phi\circ f)\oplus\ov{f'_1}\  : \
A'_1/\imm(f'_2)\lra F\oplus A'_0 . \]
Alors on a \ \m{U = (F\oplus A'_0)/\imm(\eta)}. On d\'efinit maintenant \
\m{\psi : F\oplus E\to U}, c'est-\`a-dire
\[ \psi : (A\oplus A'_0)/\imm(\mu)\lra 
((A\oplus A'_0)/\imm(\mu)\oplus A'_0) . \]
On consid\`ere le morphisme
\[{\xymatrix@R=6pt{
\psi_0 : F\oplus A\oplus A'_0\ar[r] & F\oplus A'_0 \\
\ \ \ \ (v,a,a'_0)\fmaps[r] & (v+\phi(a),a'_0) }}\]
Il est clair qu'on a \ \m{\psi_0(F\oplus\imm(\mu))\subset\imm(\eta)}. On en
d\'eduit $\psi$. La v\'erification du fait que le diagramme de la proposition
\ref{ext_const1} est commutatif ainsi que les d\'emonstrations des autres
assertions sont laiss\'ees au lecteur.
\end{proof}

\sepprop

\pagebreak[2]\begin{subsub}{\bf Corollaire : }\label{ext_const2}
Soient $A$, $E$, $A'$ des faisceaux coh\'erents sur $X$,
\xmat{
0\ar[r] & A\ar[r] & E\ar[r] & A'\ar[r] & 0 }
une suite exacte, et \ \m{\delta' : \Hom(A,A')\to\Ext^1(A',A')}, \
\m{\delta : \Hom(A,A')\to\Ext^1(A,A)} \ les applications induites. 
Soit \ \m{\phi\in\Hom(A,A')}. Alors il existe un diagramme 3x3
\Bdiag{A}{B}{A}{E}{E\oplus E}{E}{A'}{B'}{A'}
tel que les suites exactes verticales de gauche et de droite soient
identiques \`a la suite exacte pr\'ec\'edente, que la suite exacte horizontale
du milieu soit triviale, que celle du haut soit associ\'ee \`a $\delta(\phi)$,
que celle du bas soit associ\'ee \`a $\delta'(\phi)$,
et que le morphisme \ \m{\alpha : A\to A'} \ induit par la seconde composante
du morphisme \ \m{E\oplus E\to B'} \ soit \'egal \`a $\phi$. 

R\'eciproquement, \'etant donn\'e un tel diagramme, la suite exacte
horizontale du haut est associ\'ee \`a $\delta(\alpha)$ et celle du bas est 
associ\'ee \`a $\delta'(\alpha)$. 
\end{subsub}

\begin{proof}
On applique la proposition \ref{ext_const1} (avec $F=E$) et on d\'eduit de 
$\phi$ un premier diagramme commutatif avec lignes et colonnes exactes :
\xmat{
 & & 0\ar[d] & 0\ar[d] & \\
 & & A\fleq[r]\ar[d] & A\ar[d] & \\
0\ar[r] & A'\fleq[d]\ar[r] & A'\oplus E\ar[r]\ar[d] & E\ar[d]\ar[r] & 0 \\
0\ar[r] & A'\ar[r] & B'\ar[r]\ar[d] & A'\ar[r]\ar[d] & 0 \\
 & & 0 & 0 & }
On en d\'eduit le diagramme 3x3 du corollaire \ref{ext_const2} en 
rempla\c cant le morphisme \m{A'\to B'} par la compos\'ee \
\m{E\to A'\to B'}. Il reste \'a v\'erifier que la suite exacte horizontale du
haut est associ\'ee \`a $\delta(\phi)$. Pour cela on d\'eduit du diagramme
3x3 le diagramme commutatif avec lignes et colonnes exactes
\xmat{
 & 0\ar[d] & 0\ar[d] & & \\
0\ar[r] & A\ar[r]\ar[d] & B\ar[r]\ar[d] & A\ar[r]\fleq[d] & 0 \\
0\ar[r] & E\ar[d]\ar[r] & E\oplus A\ar[r]\ar[d] & A\ar[r] & 0 \\
 & A'\ar[d]\fleq[r] & A'\ar[d] & & \\
 & 0 & 0 & & }
o\`u le morphisme induit $A\to A'$ est $\phi$ et la suite exacte horizontale 
du haut est celle du diagramme 3x3. Toujours d'apr\'es la proposition
\ref{ext_const1}, cette suite exacte est associ\'ee \`a $\delta(\phi)$.

Le reste du corollaire \ref{ext_const2} est laiss\'e au lecteur.
\end{proof}

\end{sub}

\sepsub

\Ssect{Extensions duales}{Ext_dual}

\pagebreak[2]\begin{subsub}{\bf Proposition : }\label{dual_ext}
Soit 
\xmat{0\ar[r] & E\ar[r] & \ke\ar[r] & F\ar[r] & 0 }
une extension de faisceaux localement libres
sur $X$, associ\'ee \`a \ \m{\sigma\in\Ext^1(F,E)}. Soit 
\xmat{ F_2\ar[r]^{f_2} & F_1\ar[r]^{f_1} & F_0\ar[r]^{f_0} & 
F\ar[r] & 0 }
une r\'esolution localement libre de $F$.
On suppose que $\sigma$ provient d'un morphisme \ \m{\eta : F_1\to E} \
tel que \ \m{\eta\circ f_2=0}. Supposons qu'on ait une r\'esolution localement 
libre de $E^*$
\xmat{D_2\ar[r] & D_1\ar[r]^{d_1} & D_0\ar[r]^{d_0} & E^*\ar[r] & 0 }
telle qu'il existe un diagramme commutatif
\xmat{
0\ar[r] & F^*\ \flinc[r]^-{^t{f_0}} & F_0^*\ar[r]^-{^t{f_1}} & F_1^* \\
D_2\ar[u]\ar[r]^{d_2} & D_1\ar[u]^\lambda\ar[r]^{d_1} 
& D_0\ar[u]^\alpha\ar[r]^{d_0} & E^*\ar[u]^{^t\eta} }
Alors le morphisme \ \m{D_1/\imm(d_2)\to F^*} \ induit par $\lambda$ est 
associ\'e \`a l'extension
\xmat{0\ar[r] & F^*\ar[r] & \ke^*\ar[r] & E^*\ar[r] & 0 }
duale de $\sigma$.

Il existe toujours une telle r\'esolution localement libre de $E^*$.
\end{subsub}

\begin{proof}
Pour d\'emontrer l'existence de la r\'esolution donnant un tel diagramme 
commutatif on prend le m\^eme type de r\'esolutions que dans 
\ref{Extensconstr}.
Les constructions de la r\'esolution et du diagramme se font pas \`a pas en 
appliquant le th\'eor\`eme B de Serre. 

\medskip

Soit \ \m{0\lra F^*\lra V\lra E^*\lra 0} \
l'extension d\'eduite de $\lambda$. On va voir que $V$ est isomorphe \`a
$\ke^*$. Rappelons les d\'efinitions 
de $\ke$ et $V$. Soient
\[\mu = \eta\oplus f_1: F_1\lra E\oplus F_0 , \ \ \ \
\nu = \lambda\oplus d_1 : D_1\lra F^*\oplus D_0 ,\]
Alors on a \ \m{\ke=(E\oplus F_0)/\imm(\mu)} \ et \
\m{V=(F^*\oplus D_0)/\imm(\nu)}. Soit
\[{\xymatrix@R=6pt{
I : F^*\oplus D_0\ar[r] & E^*\oplus F_0^* \\
\ \ \ \ (\phi,\delta_0)\fmaps[r] & (d_0(\delta_0),{}^tf_0(\phi) -
\alpha(\delta_0)) \\
}}\]
Il est ais\'e de montrer que $I$ induit un isomorphisme \ 
$\ov{I}:V\simeq\ke^*$. Pour achever la d\'emonstration on v\'erifie 
ais\'ement qu'on a un diagramme commutatif
\xmat{
0\ar[r] & F^*\ar[r]\fleq[d] & V\ar[d]^{\ov{I}}\ar[r] & E^*\fleq[d]\ar[r] 
& 0 \\
0\ar[r] & F^*\ar[r] & \ke^*\ar[r] & E^*\ar[r] & 0 }
\end{proof}

\sepprop

\pagebreak[2]\begin{subsub}{\bf Remarque : }\rm
On peut obtenir un diagramme comme celui de la proposition
\ref{dual_ext} pour tout morphisme \ $D_1\to F^*$ \ repr\'esentant
l'extension duale. En effet, si \ \m{\psi : D_0\to F^*}, on a un
diagramme commutatif
\xmat{
0\ar[rr] & & F^*\ \flinc[rr]^-{^t{f_0}} &  & F_0^*\ar[r]^-{^t{f_1}} & F_1^* \\
D_2\ar[u]\ar[rr]^{d_2} & & D_1\ar[u]^{\lambda+\psi d_1}\ar[rr]^{d_1} &
& D_0\ar[u]^{\alpha+{}^tf_0\psi}\ar[r]^{d_0} & E^*\ar[u]^{^t\eta} }
\end{subsub}
\end{sub}

\sepsub

\Ssect{R\'esolutions d'extensions}{Ext_res}

Soit
\[  
(L) \ \ \ \ \
\xymatrix{0\ar[r] & E\ar[r] & \ke\ar[r] & F\ar[r] & 0 \\ }
\]
une extension de faisceaux coh\'erents sur $X$. Soit
\xmat{\cdots E_2\ar[r]^{e_2} & E_1\ar[r]^{e_1} & E_0\ar[r]^{e_0} &
E\ar[r] & 0 }
une r\'esolution localement libre de $E$. Soit
\xmat{\cdots F_2\ar[r]^{f_2} & F_1\ar[r]^{f_1} & F_0\ar[r]^{f_0} &
F\ar[r] & 0 }
une r\'esolution localement libre de $F$  telle que $(L)$ provienne d'un
morphisme \ \m{\sigma_1 : F_1\lra E} \ s'annulant sur $\imm(f_2)$.
(cf. \ref{Extensconstr}). On en d\'eduit une r\'esolution localement libre
de $\ke$ :
\xmat{\cdots F_2\ar[r]^{f_2} & F_1\ar[r]^{(f_1,\sigma_1)} & 
F_0\oplus E\ar[r] & \ke\ar[r] & 0 }
Cette r\'esolution n'est malheureusement pas arbitrairement n\'egative
(\`a cause du terme $F_0\oplus E$), et ne peut donc pas \^etre utilis\'ee
pour calculer des espaces du type $\Ext^p(\ke,\ku)$, o\`u $\ku$ est un
faisceau coh\'erent sur $X$. En fait seuls les \'el\'ements qui 
proviennent de $\Ext^1(F,\ku)$ peuvent \^etre construits (si $p>0$).
Cependant, si on prend les $F_i$ suffisamment
n\'egatifs, on d\'eduit ais\'ement de la r\'esolution pr\'ec\'edente de 
$\ke$ une nouvelle r\'esolution
\xmat{\cdots F_2\oplus E_2\ar[r]^{\delta_2} & 
F_1\oplus E_1\ar[r]^{\delta_1} & 
F_0\oplus E_0\ar[r]^{\delta_0} & \ke\ar[r] & 0 }
o\`u $\delta_0$ est la compos\'ee
\xmat{F_0\oplus E_0\ar[r]^{(I,e_0)} & F_0\oplus E\ar[r] & \ke}
et $\delta_i$ pour $i>1$ est repr\'esent\'e par une matrice
\[\left(\begin{array}{cc}f_i & 0 \\ \theta_i & e_i
\end{array}\right)\]
avec \ \m{\theta_i : F_i\to E_{i-1}}. Les $\theta_i$ v\'erifient 
l'\'equation \ 
\[\theta_if_{i+1}+e_i\theta_{i+1} \ = \ 0\]
pour $i\geq 1$.
Si les r\'esolutions de $E$ et $F$
sont ad\'equates, tout \'el\'ement $u$ de $\Ext^1(\ke,\ku)$ peut \^etre
repr\'esent\'e par un morphisme
\[(\phi_p,\epsilon_p) :F_p\oplus E_p\to\ku\]
s'annulant sur $\imm(\delta_{p+1})$, et $u$ provient de $\Ext^p(F,\ku)$ si
et seulement si on peut le repr\'esenter par un morphisme tel que \
$\epsilon_p=0$. En g\'en\'eral l'image de $u$ dans $\Ext^p(E,\ku)$ est
repr\'esent\'ee par \ \m{\epsilon_p:E_p\to\ku}.

On a une suite exacte de complexes
\xmat{\cdots E_2\ar[r]^{e_2}\flinc[d] & E_1\ar[r]^{e_1}\flinc[d] 
& E_0\ar[r]^{e_0}\flinc[d] &E\ar[r]\flinc[d] & 0 \\
\cdots F_2\oplus E_2\ar[r]^{\delta_2}\flon[d] & 
F_1\oplus E_1\ar[r]^{\delta_1}\flon[d] & 
F_0\oplus E_0\ar[r]^{\delta_0}\flon[d] & \ke\ar[r]\flon[d] & 0 \\
\cdots F_2\ar[r]^{f_2} & F_1\ar[r]^{f_1} & 
F_0\ar[r]^{f_0} & F\ar[r] & 0}

R\'eciproquement, si on a une suite exacte du type pr\'ec\'edent, les 
morphismes $\delta_i$ sont d\'efinis par des matrices \ 
\m{\left(\begin{array}{cc}f_i & 0 \\ \theta_i & e_i
\end{array}\right)}, avec \ \m{\theta_i : F_i\to E_{i-1}}, et 
\m{(\theta_i)_{i\geq 1}} d\'efinit un \'el\'ement de $\Ext^1(F,E)$ qui
est celui qui est associ\'e \`a l'extension de $F$ par $E$ induite par la
suite exacte de complexes. 

\end{sub}

\sepsec


\pagebreak[4]\section{Faisceaux r\'eguliers}\label{f_reg}

Soit $X$ une vari\'et\'e alg\'ebrique projective lisse et irr\'eductible

\sepsub

\Ssect{D\'efinition des faisceaux r\'eguliers}{reg_def}

\medskip

\pagebreak[2]\begin{subsub}{\bf D\'efinition : }\label{f_reg_def}
Soit $F$ un faisceau coh\'erent $X$. On dit que $F$ est {\em r\'egulier} si
\begin{itemize}
\item[-](i)  $F$ est sans torsion. 
\item[-](ii) $F^{**}$ est simple et 2-lisse (cf. \ref{F_liss}).
\item[-](iii) Si $F$ n'est pas localement libre, $F^{**}/F$ est parfait de
codimension 2 (cf \ref{grat}).
\end{itemize}
\end{subsub}

\sepprop

\pagebreak[2]\begin{subsub}{\bf Remarques : }\rm
1 - La simplicit\'e de $F^{**}$ entraine celle de $F$.

2 - Si $X$ est une surface, on voit ais\'ement en utilisant la dualit\'e de
Serre que la 2-lissit\'e de $F^{**}$ entraine celle de $F$, et que la condition
(iii) est une cons\'equence de (i).
La d\'efinition d'un faisceau r\'egulier est donc plus simple dans ce cas :
c'est un faisceau isomorphe au noyau d'un morphisme \m{E\to T}, o\`u $T$ est un
faisceau de dimension 0 et $E$ un fibr\'e vectoriel simple et 2-lisse.
\end{subsub}

\sepprop

\pagebreak[2]\begin{subsub}{\bf Lemme : }\label{reg_com}
Soient $F$ un faisceau r\'egulier non localement libre sur $X$ et 
\m{T=F^{**}/F}. Alors on a un
diagramme commutatif canonique avec lignes et colonnes exactes
\[
\xymatrix{
 & 0\ar[d] & 0\ar[d] \\ 0\ar[r] & \Hom(F^{**},T)/\End(T)\ar[r]\ar[d] &
 \Ext^1(F^{**},F)/(\End(T)/\C)\ar[d] \\
 0\ar[r] & \Hom(F,T)\ar[r] & \Ext^1(F,F)
}\]
\end{subsub}

\begin{proof}
On d\'eduit de la suite exacte \ \m{0\lra F\lra F^{**}\lra T\lra 0} \
un diagramme commutatif
\[
\xymatrix{
& & 0\ar[d] \\ & 0\ar[d] & \C\ar[d] \\ & \End(T)\ar[r]^\simeq\ar[d] &
\Ext^1(T,F)\ar[d] \\ 0\ar[r] & \Hom(F^{**},T)\ar[r]\ar[d] &
\Ext^1(F^{**},F)\ar[d] \\ 0\ar[r] & \Hom(F,T)\ar[r] & \Ext^1(F,F)
}\]

dont d\'ecoule ais\'ement le lemme.
\end{proof}

\sepprop

Soit $F$ un faisceau r\'egulier sur $X$. On notera
\[T(F) \ = \ F^{**}/F ,\]
\[ D(F) \ \ = \ \ \Ext^1(F^{**},F)/(\End(T)/\C) \ + \ \Hom(F,T) \ \ 
\subset \ \Ext^1(F,F) . \]
On va voir dans \ref{Def_reg} que \m{D(F)} correspond aux d\'eformations de $f$
obtenues en d\'eformant \m{F^{**}}, $T$ et le morphisme surjectif 
\m{F^{**}\to T}.
\end{sub}

\sepsub

\Ssect{D\'eformations de faisceaux r\'eguliers}{Def_reg}

Soit $X$ une vari\'et\'e alg\'ebrique projective lisse et connexe.
Soient $\kt$ une famille de faisceaux parfaits de codimension 2 sur $X$, 
param\'etr\'ee par une vari\'et\'e alg\'ebrique irr\'eductible $Y$ et $\kf$ une
famille de faisceaux localement libres sur $X$ param\'etr\'ee par une 
vari\'et\'e alg\'ebrique irr\'eductible $Z$ (cf. \ref{modfin}). On note \ 
\m{p_Y:Y\times Z\to Y}, 
\m{p_Z:Y\times Z\to Z}, \m{p_{Y\times Z}:Y\times Z\times X\to Y\times Z}
 \ les projections. On suppose que les propri\'et\'es suivantes sont 
v\'erifi\'ees :

(i) \m{\dim(\Hom(\kf_z,\kt_y))} \ est ind\'ependant des points ferm\'es $y$ de 
$Y$ et $z$ de $Z$. 

(ii) si $Y$ n'est pas r\'eduite, pour tout \m{y\in Y}, $\kt_y$ est 2-lisse. 

(iii) pour tout \m{z\in Z}, $\kf_z$ est simple et 2-lisse.

\sepprop

\pagebreak[2]\begin{subsub}{\bf Lemme : }\label{ext_reg0}
Le faisceau coh\'erent \ 
\m{\kh=p_{Y\times Z*}(\HHom(p_Y^\sharp(\kf),p_Z^\sharp(\kt)))}
 \ sur \ \m{Y\times Z} \ est localement libre.
\end{subsub}

\begin{proof}
La dimension de $\Hom(\kf_z,\kt_y)$ ne d\'epend pas des points $y$ de $Y$ et
$z$ de $Z$ d'apr\`es (i). Donc si $Y$ et $Z$ sont r\'eduites $\kh$ est bien
localement libre (cf. par exemple \cite{ha}, cor. II.12.9). Supposons que $Y$ 
soit non r\'eduite. Pour voir que $\kh$ est
localement libre on montre que c'est le cas au voisinage
de toute paire de points ferm\'es $(y,z)$. On utilise des d\'eformations 
semi-universelles de $\kt_y$ et $\kf_z$ (cf. \ref{S_T}). Elles sont
param\'etr\'ees par des germes lisses $U_y$ et $V_z$ respectivement (car $\kf_z$
est lisse et $\kt_y$ aussi d'apr\`es (ii)).
Le faisceau $\kh'$ analogue \`a $\kh$ sur \ $U_y\times V_z$ \
est alors localement libre. Il en est donc de m\^eme de $\kh$ au voisinage
de $(x,y)$ puisque que c'est l'image r\'eciproque de $\kh'$ par le
morphisme universel \ \m{Y\times Z\to U_y\times V_z} \ d\'efini au voisinage 
de $(y,z)$. Dans le cas o\`u est $Y$ est r\'eduite, il suffit de consid\`erer
une d\'eformation semi-universelle de \m{\kf_z}.
\end{proof}

\sepprop

On note $\kh^0$ l'ouvert de $\kh$ (vu comme fibr\'e
vectoriel sur \ $Y\times Z$) constitu\'e des morphismes surjectifs \
\m{\kf_z\to\kt_y}. Si $p$ d\'esigne la
projection \ \m{\kh^0\to Y\times Z} \ on a donc un morphisme universel surjectif
de faisceaux sur \m{\kh^0\times X}
\[ p^\sharp(p_Z^\sharp(\kf))\lra p^\sharp(p_Y^\sharp(\kt)) \]
dont le noyau $\kg$ est une famille de faisceaux r\'eguliers sur $X$.

Soient \m{\P(\kh^0)} l'ouvert du fibr\'e en espaces projectifs \m{\P(\kh)}
constitu\'e des points au dessus de \m{\kh^0} et \ \m{q:\P(\kh^0)\to Y\times Z},
\m{p_\P : \P(\kh^0)\times X\to\P(\kh^0)} \
les projections. On a comme pr\'ec\'edemment un morphisme universel surjectif
de faisceaux sur \m{\P(\kh^0)\times X}
\[ q^\sharp(p_Z^\sharp(\kf))\ot p_\P^*(\ko_{\P(\kh)}(-1))\lra 
q^\sharp(p_Y^\sharp(\kt)) \]
dont le noyau $\kg'$ est une famille de faisceaux r\'eguliers sur $X$.

\sepprop

\pagebreak[2]\begin{subsub}{\bf Proposition : }\label{ext_reg}
Soit $h$ un point ferm\'e de $\kh^0$, au dessus de \ \m{(y,z)\in Y\times Z}. 
Alors 

1 - Le morphisme de d\'eformation infinit\'esimale de Koda\"\i ra-Spencer
\[ \omega_h : T\kh_h\lra\Ext^1(\kg_h,\kg_h) \]
est \`a valeurs dans $D(\kg_h)$.

2 - Si $\kt$, $\kf$ sont des d\'eformations compl\`etes de $\kt_y$, $\kf_z$
respectivement, l'image de $\omega_h$ est \'egale \`a $D(\kg_h)$.
\end{subsub}

\medskip 

La proposition \ref{ext_reg} est une cons\'equence imm\'ediate des deux
lemmes qui suivent.

Soient \ \m{\kh^y=\kh_{\mid \lbrace y\rbrace\times Z}},
\m{\kg^y=\kg_{\mid\kh^y}}, \m{\kh^z=\kh_{\mid Y\times\lbrace z\rbrace}},
\m{\kg^z=\kg_{\mid\kh^z}}. Soient
\[ \omega_h^y : T\kh_h^y\lra\Ext^1(\kg_h,\kg_h), \ \ \ \
\omega_h^z : T\kh_h^z\lra\Ext^1(\kg_h,\kg_h) \]
les morphismes de d\'eformation infinit\'esimale de Koda\"\i ra-Spencer de
$\kg^y$, $\kg^z$ respectivement, au point $h$, qui sont des restrictions
de $\omega_h$. Soient enfin
\[ \omega_z : T\kh_h^y\lra\Ext^1(\kf_z,\kf_z), \ \ \ \
\omega_y : T\kh_h^z\lra\Ext^1(\kt_y,\kt_y) \] 
les morphismes de d\'eformation infinit\'esimale de Koda\"\i ra-Spencer de
$p_Z^\sharp(\kf)$, $p_Y^\sharp(\kt)$ respectivement, au point $h$.

\sepprop

\pagebreak[2]\begin{subsub}{\bf Lemme : }\label{ext_reg2}
(i) L'application $\omega_h^y$ est \`a valeurs dans le sous-espace vectoriel 
\hfil\break
\m{\Ext^1(\kg_h^{**},\kg_h)/(\End(\kt_y)/\C)} \ de \m{\Ext^1(\kg_h,\kg_h)}.

(ii) La compos\'ee 
\xmat{
T\kh_h^y\ar[r]^-{\omega_h^y} &
\Ext^1(\kg_h^{**},\kg_h)/(\End(\kt_y)/\C)\flon[r] & 
\Ext^1(\kg_h^{**},\kg_h^{**}) }
est \'egale \`a $\omega_z$.

(iii) L'image de $\omega_h^y$ est exactement \
\m{\Ext^1(\kg_h^{**},\kg_h)/(\End(\kt_y)/\C)} \ si $\kf$ est une 
d\'eformation compl\`ete de $\kf_z$.
\end{subsub}

\begin{proof}
D\'emontrons (i). Soit \ \m{\phi : Z_2\to\kh^y} \ le 
morphisme correspondant \`a un \'el\'ement $u$ de \m{T\kh_h}. On a une 
suite exacte sur \ $\kh^y\times X$
\xmat{
0\ar[r] & \kg^y\ar[r] & p_Z^\sharp(\kf)\ar[r] & p_X^*(\kt_y)\ar[r] & 0 \\  }
(o\`u \ \m{p_Z : \kh_y\to Z} \ et \ \m{p_X : \kh_y\to X} \ sont les 
projections). On en d\'eduit la suite exacte
\xmat{
0\ar[r] & \phi^\sharp(\kg^y)\ar[r] & \phi^\sharp(p_Z^\sharp(\kf))\ar[r] & 
p_X^*(\kt_y)\ar[r] & 0 \\  }
sur \ \m{Z_2\times X}. En projetant sur $X$ on obtient le diagramme 3x3
suivant :
\Bdiag{\kg_h}{\ku}{\kg_h}{\kf_z}{\kv}{\kf_z}{\kt_y}
{\kt_y\oplus\kt_y}{\kt_y}
o\`u \
\m{\ku = p_{X*}(\phi^\sharp(p_Z^\sharp(\kg^y)))} , \ 
\m{\kv = p_{X*}( \phi^\sharp(p_Z^\sharp(\kf)))} .

Dans le diagramme 3x3 pr\'ec\'edent on a \
\m{\kt_y\oplus\kt_y \ = \ p_{X*}(p_X^*(\kt_y))} ,
et la suite exacte horizontale du bas est la suite exacte triviale. 
L'\'el\'ement de \m{\Ext^1(\kg_h,\kg_h)} associ\'e \`a la suite exacte
horizontale du haut n'est autre que \m{\omega^y_h(u)}, d'apr\`es \ref{K-S}.
Soit \ \m{\kv'\subset\kv} \ le sous-faisceau image inverse du second facteur 
\m{\kt_y} par le morphisme \ \m{\kv\to\kt_y\oplus\kt_y}. On a un diagramme 
commutatif canonique avec lignes exactes
\xmat{
0\ar[r] & \kg_h\fleq[d]\ar[r] & \ku\ar[r]\ar[d] & \kg_h\ar[r]\flinc[d] & 0 \\
0\ar[r] & \kg_h\ar[r] & \kv'\ar[r] & \kf_z\ar[r] & 0 }

Soit $\eta$ l'\'el\'ement de \ \m{\Ext^1(\kf_z,\kg_h)= 
\Ext^1(\kg_h^{**},\kg_h)} \ correspondant \`a la suite exacte du bas. Alors
le diagramme commutatif pr\'ec\'edent montre que $\omega^y_h(u)$ est l'image 
de $\eta$ par l'application
\[\Ext^1(\kg_h^{**},\kg_h)\lra\Ext^1(\kg_h,\kg_h)\]
provenant de l'inclusion \ \m{\kg_h\subset\kg_h^{**}}. Ceci d\'emontre (i).

Pour d\'emontrer l'assertion (ii) du lemme \ref{ext_reg2}, on remarque que
$\omega_z(u)$ n'est autre que l'extension horizontale du milieu du
diagramme 3x3 pr\'ec\'edent, et (ii) d\'ecoule du diagramme commutatif avec 
lignes exactes
\xmat{
0\ar[r] & \kg_h\flinc[d]\ar[r] & \kv'\ar[r]\flinc[d] & \kf_z\ar[r]\fleq[d] 
& 0 \\
0\ar[r] & \kf_z\ar[r] & \kv\ar[r] & \kf_z\ar[r] & 0 \\ }
Pour d\'emontrer (iii) on se restreint \`a l'ouvert
$U$ de $\Hom(\kf_z,\kt_y)$ correspondant aux morphismes surjectifs. On montre
que la restriction de $\omega_h^y$ \`a $TU_h$ est \`a valeurs dans
\[\Hom(\kf_z,\kt_y) \ \subset \ 
\Ext^1(\kg_h^{**},\kg_h)/(\End(\kt_y)/\C)\]
et que c'est m\^eme l'identit\'e de $\Hom(\kf_z,\kt_y)$. On en d\'eduit 
l'assertion (iii) \`a l'aide de (ii).
\end{proof}

\sepprop

\pagebreak[2]\begin{subsub}{\bf Lemme : }\label{ext_reg3}
(i) L'application $\omega_h^z$ est \`a valeurs dans le sous-espace vectoriel
 \ \m{\Hom(\kg_h,\kt_y)} \ de \m{\Ext^1(\kg_h,\kg_h)}.

(ii)  La compos\'ee 
\xmat{
T\kh_h^z\ar[r]^-{\omega_h^z} & \Hom(\kg_h,\kt_y)\flon[r] & \Ext^1(\kt_y,\kt_y)}
est \'egale \`a $\omega_y$.

(iii) L'image de $\omega_h^z$ est exactement \
\m{\Hom(\kg_h,\kt_y)} \ si $\kt$ est une d\'eformation
compl\`ete de $\kt_y$.
\end{subsub}

\begin{proof}
La d\'emonstration de (i) est analogue \`a celle du (i) du lemme 
\ref{ext_reg2}. On aboutit \`a un diagramme 3x3
\Bdiag{\kg_h}{\ku}{\kg_h}{\kf_z}{\kf_z\oplus\kf_z}{\kf_z}{\kt_y}{\kw}
{\kt_y}
o\`u la suite exacte horizontale du haut est associ\'ee \`a $\omega^z_h(u)$.
On en d\'eduit un diagramme commutatif avec lignes exactes
\xmat{
0\ar[r] & \kg_h\flinc[d]\ar[r] & \ku\ar[r]\ar[d] & \kg_h\ar[r]\fleq[d] & 0 \\
0\ar[r] & \kf_z\ar[r] & \kf_z\oplus\kg_h\ar[r] & \kg_h\ar[r] & 0 }
qui montre que l'image de $\omega_h^z$ dans \ \m{\Ext^1(\kg_h,\kf_z)=
\Ext^1(\kg_h,\kg_h^{**})} \ est nulle, et donc que $\omega^z_h(u)$ provient
de $\Hom(\kg_h,\kt_y)$.

L'assertion (ii) est une cons\'equence du corollaire \ref{ext_const2} et (iii)
se d\'emontre comme le (iii) du lemme \ref{ext_reg2}.
\end{proof}

\sepprop

On a des r\'esultats analogues \`a la proposition \ref{ext_reg} et aux lemmes
\ref{ext_reg2}, \ref{ext_reg3}, en rempla\c cant \m{\kh^0} et $\kg$ par
\m{\P(\kh^0)}, $\kg'$ respectivement.

\end{sub}

\sepsub

\Ssect{Faisceaux r\'eguliers sur les surfaces}{reg_surf}

\pagebreak[2]\begin{subsub}{\bf Proposition : }\label{reg_3x3}
On suppose que $X$ est une surface.
Soient $F$ un faisceau r\'egulier non localement libre sur $X$, et \
\m{T=F^{**}/F}. 

Soient
\[\eta : \Ext^1(F,F^{**})\to\Ext^1(F,T)\simeq\Ext^2(T,T) , \ \ \ \
\lambda : \Ext^2(T,F^{**})\to\Ext^2(T,T) \]
les applications induites par $\pi$. Alors l'image de $\eta$ est 
\m{\Ad^2(T)} et $\lambda$ est surjective. On pose \ \m{B(F)=\ker(\eta)}.
Alors on a un diagramme 3x3 canonique

\Bdiag{\Hom(F^{**},T)/\End(T)}{\Ext^1(F^{**},F)/(\End(T)/\C)}
{\Ext^1(F^{**},F^{**})}{\Hom(F,T)}{\Ext^1(F,F)}
{B(F)}{\Ext^1(T,T)}{\Ext^2(T,F)}{\ker(\lambda)}
\end{subsub}

\begin{proof} 
On a un diagramme 
\xmat{
\Ext^1(F,F^{**})\ar[r]^\eta & \Ext^1(F,T)\ar[r]\ar[d]^\simeq &
\Ext^2(F,F)\ar[d]^{tr_2(F)} \\
& \Ext^2(T,T)\ar[r]^{tr_2(T)} & H^2(\ko_X)
}
o\`u la ligne du haut est exacte et le carr\'e est anticommutatif (cf.
\ref{trace}). Puisque \m{tr_2(F)} est un isomorphisme on a
\ \m{\imm(\eta)=\ker(tr_2(T))}. La seconde assertion d\'ecoule du fait que
\ \m{\coker(\lambda)\subset\Ext^3(T,F)}. 

Montrons maintenant que l'application \m{\Ext^1(F,F)\to\Ext^2(T,F)} est
surjective. On a une suite exacte
\xmat{\Ext^1(F,F)\ar[r] & \Ext^2(T,F)\ar[r] & \Ext^2(F^{**},F)\ar[r]^-\phi
& \Ext^2(F,F)}
donc il suffit de montrer que $\phi$ est injective, et c'est imm\'ediat en
utilisant la dualit\'e de Serre et la lissit\'e de $F$ et \m{F^{**}}.

Le reste du lemme \ref{reg_3x3} d\'ecoule ais\'ement du lemme \ref{reg_com} et
du diagramme commutatif avec lignes et colonnes exactes suivant :
\xmat{
\Ext^1(T,F)\ar[r]\ar[d] & \nsp=\Ext^1(T,F^{**})\ar[d] \\
\Ext^1(F^{**},F)\ar[r]\ar[d] & \Ext^1(F^{**},F^{**})\ar[d]\ar[r] & 
\nsp=\Ext^1(F^{**},T)\ar[d] \\
\Ext^1(F,F)\ar[r]\ar[d] & \Ext^1(F,F^{**})\ar[r]^\eta\ar[d] & 
\Ext^1(F,T)\ar[r]\ar[d]^\simeq & \Ext^2(F,F) \\
\Ext^2(T,F)\ar[r]\ar[d] & \Ext^2(T,F^{**})\ar[r]^\lambda\ar[d] & 
\Ext^2(T,T)\ar[d]\ar[dl]^{tr_2(T)} \\
0 & \Ext^2(F^{**},F^{**})=H^2(\ko_X) & \nsp=\Ext^2(F^{**},T)
}
\end{proof}

\sepprop

\pagebreak[2]\begin{subsub}{\bf Remarque : }\rm On voit que $\lambda$ ne peut 
\^etre bijective que si $F$ est de rang 1. Donc si \m{rg(F)>1}, on a 
\m{D(F)\not=\Ext^1(F,F)}. On retrouve le fait bien connu que sur une surface 
les faisceaux sans torsion de rang sup\'erieur \`a 1 peuvent \^etre 
d\'eform\'es en faisceaux localement libres.
\end{subsub} 

\sepprop

\pagebreak[2]\begin{subsub} D\'eformation des faisceaux r\'eguliers sur les 
surfaces. \rm
La condition (i) de \ref{Def_reg} est toujours v\'erifi\'ee. Soient \m{k>0} un
entier, \m{X_k} l'ouvert de \m{\Hilb^k(X)} correspondant aux sous-sch\'emas
constitu\'es de $k$ points distincts et \m{U_k\subset X^k} l'ouvert
au dessus de \m{X_k}. On suppose que \m{Y=U_k}, \m{\kt} \'etant
le faisceau universel \'evident. Pour \m{1\leq i\leq k}, soient \m{p_i:U_k\to X}
la $i$-\`eme projection, et \ \m{\kf_i=(I_Z\times p_i)^*(\kf)}. Soit enfin
\[W \ = \ \P(\kf_1^*)\times_{U_k\times Z}\cdots\times_{U_k\times Z}\P(\kf_k^*) 
.\]
Pour tous \m{(x_i)\in U_k} et \m{z\in Z}, la fibre de $W$ au dessus de
\m{((x_i),z)} est \ \m{\P(\kf_{z,x_1}^*)\times\cdots\times\P(\kf_{z,x_k}^*)}.
Pour \m{1\leq i\leq k}, on note \m{\pi_i:W\times X\to \P(\kf_i)} la $i$-\`eme
projection. Soient \m{p_Z} la projection \m{W\to Z} et \m{p_U} la projection
\m{W\to U_k}. Alors on a un morphisme 
surjectif \'evident
\[p_Z^\sharp(\kf)\ot\pi_1^*(\ko_{\P(\kf_1^*)}(-1))\ot\cdots\ot
\pi_k^*(\ko_{\P(\kf_k^*)}(-1))\lra p_U^\sharp(\kt) .\]
Son noyau $\kv$ est une famille de faisceaux r\'eguliers. 

Soit \m{\Sigma_k} le groupe des permutations de \m{\lbrace 1,\ldots,k\rbrace}.
On consid\`ere l'action \'evidente de \m{\Sigma_k} sur \m{U_k}, $W$ et $\kv$.
Soient \m{\overline{W}=W/\Sigma_k}, \m{\overline{\kv}=\kv/\Sigma_k}, qui est une
famille de faisceaux r\'eguliers param\'etr\'ee par \m{\overline{W}}. Le
r\'esultat suivant d\'ecoule ais\'ement de la construction de
\m{\overline{\kv}} :
\end{subsub} 

\sepprop

\pagebreak[2]\begin{subsub}{\bf Proposition : }\label{reg_3x3a}
Pour tout point $w$ de \m{\overline{W}}, le morphisme de d\'eformation
infinit\'esimale de Koda\"\i ra-Spencer
\[\omega_w : T\overline{W}_k\lra\Ext^1(\overline{\kv}_w,\overline{\kv}_w)\]
est injectif et a pour image \m{D(\overline{\kv}_w)}.
\end{subsub} 

\sepprop

Soient $F$ un faisceau r\'egulier tel qu'il existe \m{w\in\overline{W}} tel que
\m{\overline{\kv}_w\simeq F}, $\ku$ une d\'eformation semi-universelle de $F$
param\'etr\'ee par un germe \m{(S,s_0)} et \m{f:(V,w)\to(S,s_0)} un morphisme
d\'efini par $\ku$ ($V$ \'etant un voisinage ad\'equat de $w$ dans
\m{\overline{W}}). Alors l'image de $f$ dans $S$ ne d\'epend que de $F$ et de
$\ku$, on la note \m{D(\ku,F)}. Plus pr\'ecis\'ement, si \m{\ku'} est une autre
d\'eformation semi-universelle de $F$, param\'etr\'ee par un germe
\m{(S',s'_0)}, et si \m{\phi:(S',s'_0)\to(S,s_0)} est un isomorphisme tel que
\m{\phi^\sharp(\ku)\simeq\ku'}, on a \m{\phi(D(\ku',F))=D(\ku,F)}. L'espace 
tangent de \m{D(\ku,F)} en \m{s_0} est \m{D(F)}.

Soient $T$ une vari\'et\'e alg\'ebrique, $t$ un point ferm\'e de $T$, $\ke$ un
faisceau coh\'erent sur \m{T\times X}, plat qur $T$ et tel que \m{\ke_t\simeq
F}. On dit que $\ke$ est une {\em d\'eformation r\'eguli\`ere} de 
$F$ en $t$ si
un morphisme \m{\psi:(T,t)\to(S,s_0)} tel que \m{\psi^\sharp(\ku)\simeq\ke} 
est \`a
valeurs dans \m{D(\ku,F)}. On montre ais\'ement que $\ke$ est une d\'eformation 
r\'eguli\`ere de $F$ en $t$ si et seulement si il existe un ouvert $U$ de $T$
contenant $t$ tel que \m{\ke_{U\times X}^*} soit localement libre. Si c'est le 
cas on peut choisir $U$ de telle sorte qu'il existe un unique morphisme
\m{\Psi:U\to\overline{W}} tel que 
\m{\Psi^\sharp(\overline{\kv})\simeq\ke_{U\times X}}. C'est pourquoi on peut 
dire que \m{\overline{W}}, munie de \m{\overline{\kv}}, est une sorte de 
{\em vari\'et\'e de modules} de faisceaux r\'eguliers.

\end{sub}

\sepsub

\Ssect{Faisceaux r\'eguliers sur les vari\'et\'es de dimension sup\'erieure
\`a 2}{reg_big}

\pagebreak[2]\begin{subsub}{\bf Proposition : }\label{reg_3x3b}
On suppose que \ \m{\Ext^1(F^{**},T)=\Ext^2(F^{**},T)=\Hom(F^*,\wT)=\nsp}. 
Alors on a un diagramme commutatif avec lignes et colonnes exactes
\xmat{
& 0\ar[d] & 0\ar[d] \\
0\ar[r] & \Hom(F^{**},T)/\End(T)\ar[r]\ar[d] & \Ext^1(F^{**},F)/(\End(T)/\C)
\ar[r]\ar[d] & \Ext^1(F^{**},F^{**})\ar[r]\ar[d]^\simeq & 0\\
0\ar[r] & \Hom(F,T)\ar[r]\ar[d] & \Ext^1(F,F)\ar[r]\ar[d] & 
\Ext^1(F,F^{**})\ar[r] & 0\\
& \Ext^1(T,T)\ar[d]\ar[r]^\simeq & \Ext^2(T,F)\ar[d]\\ & 0 & 0
}
On a donc \ \m{\Ext^1(F,F)=D(F)}, et toute d\'eformation de $F$ s'obtient en
d\'eformant \m{F^{**}}, $T$ et le morphisme surjectif \m{F^{**}\to T}.
\end{subsub}

\begin{proof} 
Le seul point non imm\'ediat sont la surjectivit\'e des morphismes \hfil\break
\m{\alpha:\Ext^1(F,F)\to\Ext^1(F,F^{**})} \ et \ 
\m{\beta:\Ext^1(F,F)\to\Ext^2(T,F)}. En ce qui concerne le premier on a une 
suite exacte
\xmat{
\Ext^1(F,F)\ar[r]^\alpha & \Ext^1(F,F^{**})\ar[r]^\eta & \Ext^1(F,T)
}
et il suffit de montrer que \m{\eta=0}. Cela d\'ecoule du carr\'e commutatif
\xmat{
\Ext^1(F^{**},F^{**})\ar[r]\ar[d]^\simeq & \Ext^1(F^{**},T)=\nsp\ar[d] \\
\Ext^1(F,F^{**})\ar[r]^\eta & \Ext^1(F,T)
}
o\`u la fl\`eche de gauche est un isomorphisme parce que \ 
\m{\Ext^1(T,F^{**})=\Ext^2(T,F^{**})=\nsp}. On a une suite exacte
\xmat{
\Ext^1(F,F)\ar[r]^\beta & \Ext^2(T,F)\ar[r] & \Ext^1(F^{**},F)\ar[r]^\gamma &
\Ext^2(F,F)
}
et il suffit de montrer que $\gamma$ est injective. Cela d\'ecoule du diagramme
commutatif
\xmat{
\Ext^2(F^{**},F)\ar[d]^\gamma\ar[r]^\simeq & \Ext^2(F^{**},F^{**})\flinc[d] \\
\Ext^2(F,F)\ar[r] & \Ext^2(F,F^{**})
}
o\`u la fl\`eche du haut est un isomorphisme parce que \ \m{\Ext^1(F^{**},T)=
\Ext^2(F^{**},T)=\nsp}, et celle de droite une inclusion parce que \
\m{\Ext^2(T,F^{**})=\nsp}.
\end{proof}

\sepprop

On se place dans la situation de
\ref{Def_reg}, les conditions (i), (ii), (iii) \'etant v\'erifi\'ees. On
d\'eduit ais\'ement de la proposition \ref{reg_3x3b} la

\sepprop

\pagebreak[2]\begin{subsub}{\bf Proposition : }\label{reg_3x3b2} 
Soient $h$ un point ferm\'e de $\kh^0$, au dessus de \ \m{(y,z)\in Y\times Z}
et $q$ le point correspondant de \m{\P(\kh^0)}.
On suppose que \ \m{H^0(\kf_z\ot\widetilde{\kt_y})=\nsp}, et que $\kt$, $\kf$ 
sont des d\'eformations compl\`etes de $\kt_y$, $\kf_z$, respectivement.
Alors on a

(i) $\kg$ (resp. $\kg'$) est une d\'eformation compl\`ete de $\kg_h$ (resp.
$\kg'_q$). 

(ii) Si de plus
$Y$ est r\'eduite en $y$, la base d'une d\'eformation semi-universelle de
$\kg_h$ est aussi r\'eduite.

(iii) Si \m{\kt_y} est simple et $Y$ r\'eduite en $y$, et si $\kf$, $\kt$ sont 
des d\'eformations semi-universelles de $\kf_z$, $\kt_y$ respectivement, 
$\kg'$ est une d\'eformation semi-universelle de $\kg'_q$.
\end{subsub}

\sepprop

et le

\sepprop

\pagebreak[2]\begin{subsub}{\bf Corollaire : }\label{reg_3x3b3}
On suppose que $Y$ munie de $\kt$ et $Z$ munie de $\kf$ sont des vari\'et\'es 
de modules fins (cf. \ref{modfin}), que $Y$ est r\'eduit et que pour tout
\m{y\in Y} le faisceau \m{\kt_y} est simple. Alors \m{\P(\kh^0)} muni de
\m{\kg'} est une vari\'et\'e de modules fins. En particulier les faisceaux
\m{\kg'_q}, \m{q\in\P(\kh^0)}, constituent un ensemble ouvert.
\end{subsub}

\sepprop

\pagebreak[2]\begin{subsub}\label{reg_3x3b4}Exemples sur \m{\P_3}. \rm 
Soient $V$ un $\C$-espace vectoriel de dimension 4 et \m{\P_3=\P(V)}.
Soient $Y_0$ un ouvert
lisse d'une composante irr\'eductible d'un sch\'ema de Hilbert de courbes sur
\m{\P_3}, constitu\'e de courbes lisses et irr\'eductibles et \m{{\bf C}\to Y_0}
la courbe universelle. Soient $d$ le degr\'e des courbes de $Y_0$ et $g$ 
leur genre. Soit $m$ un entier tel que \m{1-g+m} et $d$ soient premiers entre
eux. Soit \m{{\bf J}^m\to Y_0} la jacobienne relative de degr\'e $m$.
D'apr\`es \cite{me_ra} il existe un fibr\'e de Poincar\'e $\kj_m$ sur \m{{\bf
J}^m\times_{Y_0}{\bf C}}. On peut voir $\kj_m$ comme un famille de faisceaux
parfaits de codimension 2 sur $\P_3$ param\'etr\'ee par ${\bf J}^m$. C'est 
m\^eme une vari\'et\'e de modules fins.

Soient d'autre part $\bf M$ une vari\'et\'e de modules fins de fibr\'es 
vectoriels simples et 2-lisses sur $\P_3$, $\kf$ le fibr\'e universel associ\'e.
Les conditions (i), (ii), (iii) de \ref{Def_reg} et celles des propositions
\ref{reg_3x3b}, \ref{reg_3x3b2} sont v\'erifi\'ees si \m{m\gg 0}. Dans ce cas,
d'apr\`es le corollaire \ref{reg_3x3b3}, \m{\P(\kh^0)} muni de \m{\kg'} est une
vari\'et\'e  de modules fins, constitu\'ee de faisceaux r\'eguliers non 
localement libres.

Par exemple on prend pour $Y_0$ la grassmannienne des droites de \m{\P_3}. Le 
fibr\'e de Poincar\'e $\kj_m$ existe quelque soit $m$. Soit \m{{\bf M}(2;0,1,0)}
la vari\'et\'e de modules des faisceaux semi-stables de rang 2 et de classes de
Chern 0,1,0 sur \m{\P_3}. Soit $\bf M$ l'ouvert de \m{{\bf M}(2;0,1,0)}
correspondant aux faisceaux localement libres (ce sont les {fibr\'es de
corr\'elation nulle}) (cf. \cite{O_S_S}, 4.3, ex. 3). Alors $\bf M$ est
canoniquement isomorphe au compl\'ementaire de la grassmannienne des droites de
\m{\P_3} dans \m{\P(\wedge^2V^*)} et il existe un fibr\'e universel $\kf$ sur
\m{{\bf M}\times\P_3}. Alors les conditions (i), (ii), (iii) de \ref{Def_reg} 
et celles des propositions \ref{reg_3x3b}, \ref{reg_3x3b2} sont v\'erifi\'ees 
d\'es que \m{m\geq 0}.

Si on prend pour $\bf M$ la vari\'et\'e de modules fins constitu\'ee du seul
fibr\'e $\ko$, toutes les hypoth\`eses de la proposition \ref{reg_3x3b} ne sont
pas v\'erifi\'ees (on a \m{h^0(\widetilde{\ko_\ell})=h^0(\ko_\ell(2))\not=0})
mais on a quand m\^eme \ \m{D(\ki_\ell)=\Ext^1(\ki_\ell,\ki_\ell)}.
\end{subsub}

\end{sub}

\sepsec


\pagebreak[4]\section{Extensions larges - D\'efinition}\label{Larg0}

Soient $X$ une vari\'et\'e alg\'ebrique projective lisse connexe, $\ko_X(1)$ un
fibr\'e en droites tr\`es ample sur $X$.

\sepsub

\Ssect{Motivation}{extens0}

Soient $E_0$ un faisceau localement libre, $F$ un faisceau r\'egulier sur 
$X$, non localement libre, et $n$ un entier. On consid\`ere des extensions
\xmat{0\ar[r] & E\ar[r] & \ke\ar[r] & F\ar[r] & 0 }
avec \ $E=E_0(n)$ . Soit \ \m{T=F^{**}/F}.

\sepprop

\pagebreak[2]\begin{subsub}{\bf Lemme : }\label{larg_lemm1}
Il existe un entier $n_0$ tel que si $n\geq n_0$, les propri\'et\'es
suivantes soient \hbox{v\'erifi\'ees :} 

1 - On a \ \m{\Ext^i(F^{**},E)=\nsp} \ pour $i\geq 1$ et \
\m{\Ext^i(F,E)=\nsp} \ pour $i\geq 2$.
 
2 - On a \ \m{\Ext^i(E,F^{**})=\nsp} \ pour $i<\dim(X)$.

3 - On a \ \m{\Ext^1(F',E)=\nsp} \ pour toute d\'eformation localement libre
$F'$ de $F$.
\end{subsub}

\begin{proof} Cela d\'ecoule ais\'ement du th\'eor\`eme B de Serre. Le
fait que $\Ext^i(F,E)$ est nul si $i\geq 2$
d\'ecoule du lemme \ref{lemm_parf} : on
choisit d'abord $n_0$ suffisamment grand pour qu'on ait \
\m{\Ext^i(F^{**},E)=0} \ pour $i\geq 1$, et on obtient alors des 
isomorphismes
\[\Ext^i(F,E)\simeq\Ext^{i+1}(T,E)\]
pour $i\geq 2$, et \m{\Ext^{i+1}(T,E)} est nul si $n_0$ est assez grand.
\end{proof}

\sepprop

Il d\'ecoule de 3- que les extensions non triviales de $F$ par $E$ proviennent
essentiellement des singularit\'es de $F$. On a d'ailleurs d'apr\`es le lemme 
pr\'ec\'edent un isomorphisme canonique
\[\Ext^1(F,E) \ \simeq \ \Ext^2(T,E) \ \simeq
H^0(E\ot\omega_X^{-1}\ot T^\vee).\]

\sepprop

\pagebreak[2]\begin{subsub}{\bf Lemme : }\label{lib}
Soit \ $\sigma\in\Ext^1(F,E)$ \ et \ $0\to E\to\ke\to F\to 0$
 \ l'extension associ\'ee. Alors $\ke$ est localement libre si et seulement si
le morphisme \ \m{E^*\to T(F)^\vee\ot\omega_X^{-1}} \ associ\'e \`a
$\sigma$ est surjectif.
\end{subsub}

\begin{proof}
Ce r\'esultat est analogue au lemme 5.1.2 de \cite{O_S_S} et se d\'emontre 
de la m\^eme fa\c con.
\end{proof}
\end{sub}

\sepsub

\Ssect{Dualit\'e et d\'efinition des extensions larges}{dual}

On consid\`ere une extension comme dans le lemme \ref{larg_lemm2} :
\[{\xymatrix{ (L) \ \ \ \ \ \ \ \
0\ar[r] & E\ar[r] & \ke\ar[r] & F\ar[r] & 0 \\ }} ,\]
o\`u $\ke$ est localement libre.
On d\'eduit de la suite exacte pr\'ec\'edente la suite exacte
\xmat{
0\ar[r] & F^*\ar[r] & \ke^*\ar[r] & E^*\ar[r] & \EExt^1(F, \ko_X)\simeq 
 T^\vee\ot\omega_X^{-1}\ar[r] & 0 .} 
Le morphisme surjectif pr\'ec\'edent
\ \m{\phi : E^*\to T^\vee\ot\omega_X^{-1}} est exactement 
l'\'el\'ement de 
\hfil\break
$\Ext^1(F,E)\simeq\Hom(E^*,\wT\ot\omega_X^{-1})$ associ\'e
\`a $\lambda$ (cf. lemme \ref{lib}). On pose \ \m{G=\ker(\phi)}, donc
\m{G^{**}/G= T^\vee\ot\omega_X^{-1}} . On a alors \ 
\m{E\simeq G^*}. L'extension $(L)$ s'\'ecrit donc
\xmat{
0\ar[r] & G^*\ar[r] & \ke\ar[r] & F\ar[r] & 0 }
et on obtient aussi l'extension suppl\'ementaire
\[{\xymatrix{(L^*) \ \ \ \ \ \ \ \
0\ar[r] & F^*\ar[r] & \ke^*\ar[r] & G\ar[r] & 0 \\ }.}\]
L'extension $(L)$ est enti\`erement d\'etermin\'ee par les fibr\'es 
vectoriels $F^*$, $G^*$, le faisceau parfait $T$ et deux morphismes 
surjectifs
\begin{itemize}
\item Celui qui d\'efinit $F$, \ \  $\tau : F^{**}\to T$, 
\item Celui provenant de l'\'el\'ement $\sigma$ de $\Ext^1(F,E)$ 
associ\'e \`a
$(L)$, 
\[ \phi : G^{**}\lra T^\vee\ot\omega_X^{-1} .\]
\end{itemize}
De la m\^eme mani\`ere l'extension $(L^*)$ est enti\`erement 
d\'etermin\'ee par les fibr\'es vectoriels $F^*$, $G^*$, le faisceau 
parfait
\ $\wT= T^\vee\ot\omega_X^{-1}$ \ et les morphismes surjectifs
\begin{itemize}
\item Celui qui d\'efinit $G$, \ \ $\tau^* : G^{**}\to \wT$, 
\item Celui provenant de l'\'el\'ement $\sigma^*$ de $\Ext^1(G,F^*)$ 
associ\'e
\`a $(L^*)$, 
\[\phi^* : F^{**}\lra \wT^\vee\ot\omega_X^{-1} =T.\]
\end{itemize}

\sepprop

\pagebreak[2]\begin{subsub}{\bf Proposition : }
Les morphismes correspondant \`a $(L)$ et $(L^*)$ sont les m\^emes : 
on a
\[\tau^* \ = \ \phi, \ \ \ \ \phi^* \ = \ \tau .\]
\end{subsub}

\begin{proof}Cela d\'ecoule ais\'ement de l'associativit\'e des Ext.
\end{proof}

\sepprop

\pagebreak[2]\begin{subsub}{\bf D\'efinition : }\label{larg_def}
On appelle {\em extension large} une extension
\xmat{
0\ar[r] & G^*\ar[r] & \ke\ar[r] & F\ar[r] & 0 }
comme pr\'ec\'edemment, l'extension associ\'ee 
\xmat{
0\ar[r] & F^*\ar[r] & \ke^*\ar[r] & G\ar[r] & 0 }
\'etant appel\'ee {\em extension duale},
telle que les propri\'et\'es du lemme
\ref{larg_lemm1} soient v\'erifi\'ees, c'est-\`a-dire
\[\Ext^i(F^{**},G^*) \ = \ \nsp \ \ \ \ {\rm si} \ \ i\geq 1 , \]
\[\Ext^i(F,G^*) \ = \ \Ext^i(G,F^*) \ = \ \nsp \ \ \ \ 
{\rm si} \ \ i\geq 2 , \]
\[\Ext^i(G^*,F^{**}) \ = \ \nsp \ \ \ \ {\rm si} \ \ i<\dim(X) , \]
et que :
\begin{enumerate}
\item $\ke$ est localement libre. 
\item $F$ et $G$ sont r\'eguliers. 
\end{enumerate}

\medskip

On appelle aussi extensions larges les fibr\'es vectoriels $\ke$.
\end{subsub}

\sepprop

Examinons maintenant comment on peut construire des extensions larges. Soient
$G_0$, $F_0$ des faisceaux localement libres et $T_0$ un faisceau parfait sur
$X$. On suppose que pour $p\gg 0$ il existe des morphismes surjectifs \
\m{F_0\lra T_0(p)}, \m{G_0\lra \wT_0(p)}.
Alors, si $n\gg m\gg 0$, il existe des morphismes surjectifs
\[\pi : F_0\lra T_0(m), \ \ \ \ \rho : G_0(-n)\lra \wT_0(-m)\]
d\'efinissant une extension large, avec \ \m{F=\ker(\pi)}, \m{G=\ker(\rho)}
et \m{T=T_0(m)}.
\end{sub}

\sepsub

\Ssect{Exemples}{Larg_exemp}

\medskip

\SSsect{Exemple 1 - Fibr\'es instables de rang 2 sur une surface}{ex_larg_1}

On suppose que $X$ est une surface. Soient $Y\subset X$ un sous-sch\'ema de
dimension 0 et
$\ki_Y$ le faisceau d'id\'eaux de $Y$.

Soient $\kl_0$, $\kl_1$   des fibr\'es en droites sur $X$ tels que
\[h^1(\kl_0^*\ot\kl_1^*)=h^2(\kl_0^*\ot\kl_1^*)=h^0(\kl_0\ot\kl_1)=
h^1(\kl_0\ot\kl_1)=h^0(\kl_0\ot\kl_1\ot\ki_Y\ot\omega_X)=0 .\]
Alors il existe des extensions larges
\xmat{
0\ar[r] & \kl_0^*\ar[r] & \ke\ar[r] & \kl_1\ot\ki_Y\ar[r] & 0 }

En particulier, soit $m>0$ un entier tel que
\[h^0(\ki_Y\ot\omega_X(-m)) = h^1(\ko_X(m)) = h^2(\ko_X(m)) = h^1(\ko_X(-m)) = 
\nsp.\]
Alors il existe des extensions larges 
\xmat{
0\ar[r] & \ko_X(m)\ar[r] & \ke\ar[r] & \ki_Y\ar[r] & 0 }
Le cas de $\P_2$ est \'etudi\'e dans \cite{st}. Dans ce cas la seule
condition est $m>0$.

\sepsubsub

\SSsect{Exemple 2 - Fibr\'es instables de rang 2
sur une vari\'et\'e de dimension $n>2$}{ex_larg_2}

Soit \ \m{n=\dim(X)}.
Soit $Y\subset X$ une sous-vari\'et\'e ferm\'ee localement intersection
compl\`ete, de faisceau d'id\'eaux $\ki_Y$,
$\kl_0$, $\kl_1$ des fibr\'es en droites sur $X$ tels que \
\m{\kl_0\ot\kl_1 \ \simeq \ \det(\kn_{Y/X})}
($\kn_{Y/X}$ d\'esignant le fibr\'e normal de $Y$). On suppose que
\[h^i(\kl_0^*\ot\kl_1^*) = 0 \ \ \ \ {\rm si} \ \ i\geq 1 , \ \ \ \ \ \ 
h^i(\kl_0\ot\kl_1) = 0 \ \ \ \ {\rm si} \ \ i<n , \]
\[h^i(\kl_0\ot\kl_1\ot\ki_Y\ot\omega_X) = 0 \ \ \ \ {\rm si} 
\ \ i\leq n-2 .\]
Alors il existe des extensions larges
\xmat{
0\ar[r] & \kl_0^*\ar[r] & \ke\ar[r] & \kl_1\ot\ki_Y\ar[r] & 0 }
Remarquons que si $Y$ est irr\'eductible, il existe essentiellement une
seule telle extension.

En particulier, soit $m>0$ un entier tel que
\[h^i(\ko_X(m)) = 0 \ \ \ {\rm si} \ \ i\geq 1 ,\ \ \ \ 
 h^i(\ko_X(-m)) = 0 \ \ \ {\rm si} \ \ i<n ,\ \ \ \
h^i(\ki_Y\ot\omega_X(-m)) = 0 \ \ \ \ {\rm si} \ \ 
i\leq n-2 .\]
Alors il existe des extensions larges 
\xmat{
0\ar[r] & \ko_X(m)\ar[r] & \ke\ar[r] & \ki_Y\ar[r] & 0 }

Dans le cas de $\P_n$ la seule condition est $m>0$. Pr\'ecisons que 
l'existence des sous-vari\'et\'es $Y$ n'est pas assur\'ee si $n>3$, car
l'hypoth\`ese \ \m{\det(\kn_{Y/X})\simeq\ko_X(-m)} \ entraine que $Y$ ne
peut pas \^etre une intersection compl\`ete. On
conjecture en fait qu'il n'existe pas de telles $Y$ si $n\geq 5$ 
(cf. \cite{ha}).

\sepsubsub

\SSsect{Exemple 3 - Fibr\'es de rang \'elev\'e sur $\P_2$}{ex_larg_3}

Soit $E$ un fibr\'e exceptionnel sur $\P_2$ (cf. \ref{except}). Soient
$m\geq 1, r\geq 2, c_1, c_2$ des entiers et tels que
\[\mu(E) \ < \ \frac{c_1}{r} \ < \ \mu(E)+x_E\]
et que les faisceaux coh\'erents $E'$ de rang $r$ et de classes de Chern 
$c_1$, $c_2$ v\'erifient l'\'equation \ \m{\chi(E',E)=0}. 
Soit $T$ un faisceau de dimension 0 tel qu'il existe des morphismes 
surjectifs \ $E\to T$ \ et \ $\ko\ot\C^r\to T$. D'apr\`es \cite{dr}, 
prop. 3.6,
il existe des fibr\'es stables $E'$ de rang $r$ et de classes de Chern
$c_1$, $c_2$ tels que \ \m{\Ext^1(E,E')=\nsp}. Alors il existe des 
extensions larges
\[0\lra G^*\lra\ke\lra F\lra 0\]
o\`u $G$ (resp. $F$) est le noyau d'un morphisme surjectif \ \m{{E'}^*\to \wT}
(resp. \ $E\to T$).

Dans la figure ci-dessous on a repr\'esent\'e la position du point de $\R^2$
correspondant aux fibr\'es $E'$.

\bigskip

\hskip 2cm
\includegraphics{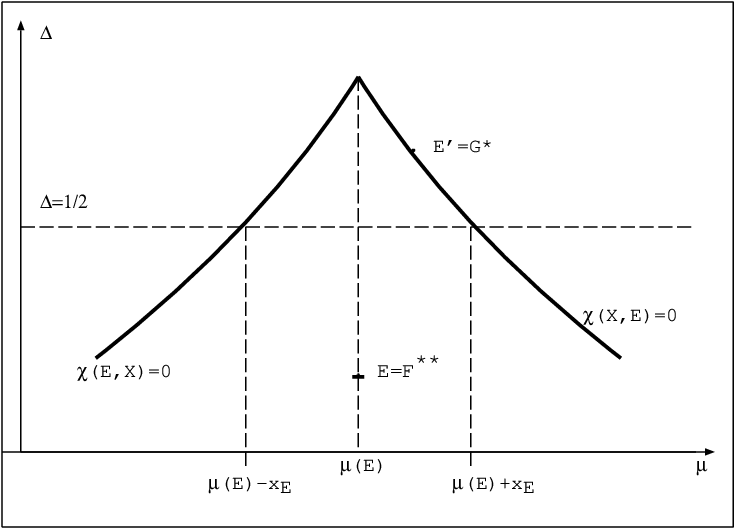}

\medskip

\centerline{{\ttx Figure 3} - Position de $E'$ dans le plan de coordonn\'ees 
$(\mu,\Delta)$}

\bigskip
 
Rappelons qu'on dit qu'un faisceau coh\'erent $\kf$ sur $\P_2$ est {\em
prioritaire} s'il est sans torsion et qu'on a \ 
\m{\Ext^2(\kf,\kf(-1))=\nsp}
(cf. \cite{hi_la}). Il est facile de voir que le faisceau $\ke$ pr\'ec\'edent
est prioritaire. On a de plus \ \m{\mu(E) < \mu(\ke) < \mu(F')}. D'autre part,
on a
\[\Delta(\ke) \ = \ C + \frac{h^0(T)}{rg(E)+r} ,\]
avec
\[C \ = \ \frac{1}{rg(E)+r}(c_2+c_2(E)+c_1c_1(E)-
\frac{rg(E)+r-1}{2(rg(E)+r)}(c_1+c_1(E))^2)). \]
Donc si \m{h^0(T)} est suffisamment grand, le point du plan \m{\R^2}
correspondant \`a $\ke$ se trouve au dessus du  graphe de la fonction
\m{\Delta=\delta(\mu)} (cf. \ref{except}). Puisque le champ des faisceaux
prioritaires est irr\'eductible d'apr\`es \cite{hi_la}, cela entraine que $\ke$
se d\'eforme en fibr\'es stables. Dans ce cas les d\'eformations de $\ke$ ne
sont pas des extensions larges du m\^eme type (cf \ref{def_larg_0}).

\end{sub}

\sepsub

\Ssect{Construction des extensions larges}{Larg_constr}

On refait maintenant les constructions pr\'ec\'edentes dans un ordre
diff\'erent. On consid\`ere des fibr\'es vectoriels $\F$, $\G$, et un 
faisceau parfait $T$ sur $X$. On pose
\[\wT \ = \ \EExt^2(T,\ko_X) \ \simeq \ 
T^\vee\ot\omega_X^{-1} .\]
La transposition d\'efinit un isomorphisme canonique \ 
\m{\End(T)\simeq\End(\wT)}.
Soient
\[\pi : \F\lra T, \ \ \ \ \rho : \G\lra\wT\]
des morphismes surjectifs. On pose \
\m{F = \ker(\pi)}, \m{G = \ker(\rho)},
de telle sorte que \ \m{\F=F^{**}}, \m{\G=G^{**}}. On se place dans le cas 
o\`u comme dans \ref{dual}, $\pi$ et $\rho$ d\'efinissent une extension 
large
\[\xymatrix{ (L) \ \ \ \ \ \ \ \
0\ar[r] & G^*\ar[r]^i & \ke\ar[r]^p & F\ar[r] & 0 }\]
associ\'ee \`a \ \m{\sigma\in\Ext^1(F,G^*)} \ et l'extension duale
\[\xymatrix{ (L^*) \ \ \ \ \ \ \ \
0\ar[r] & F^*\ar[r]^{{}^tp} & \ke^*\ar[r]^{{}^ti} & G\ar[r] & 0 }\]
associ\'ee \`a \ \m{\sigma^*\in\Ext^1(G,F^*)}. 
Rappelons qu'on a des isomorphismes canoniques
\[\delta(\pi):\EExt^1(G,\ko_X)\to T, \ \ \ \
\delta'(\rho):\EExt^1(F,\ko_X)\to\wT\]
provenant des suites exactes
\[\xymatrix{0\ar[r] & F\ar[r] & \F\ar[r]^\pi & T\ar[r] & 0}, \ \ \
\xymatrix{0\ar[r] & G\ar[r] & \G\ar[r]^\rho & \wT\ar[r] & 0}\]
et des isomorphismes \ \m{\wT\simeq\EExt^2(T,\ko_X)}, 
\m{T\simeq\EExt^2(\wT,\ko_X)}.

D'apr\`es \ref{Extensconstr} il existe une r\'esolution localement libre 
de $F$
\xmat{ \cdots F_2\ar[r]^{f_2} & F_1\ar[r]^{f_1} & F_0\ar[r]^{f_0} &
F\ar[r] & 0 }
et une r\'esolution localement libre de $G$
\xmat{ \cdots G_2\ar[r]^{g_2} & G_1\ar[r]^{g_1} & G_0\ar[r]^{g_0} &
G\ar[r] & 0 }
telles que $(L)$ provienne d'un morphisme \
\m{\sigma_1 : F_1\lra G^*} \
s'annulant sur l'image de $f_2$, et que $(L^*)$ provienne d'un morphisme
\m{\sigma_1^* : G_1\lra F^*} \ s'annulant sur l'image de $g_2$.

On retrouve les donn\'ees de d\'epart \`a partir de $\sigma_1$, $\sigma_1^*$
de la fa\c con suivante : on a un isomorphisme canonique
\[\EExt^1(G,\ko_X) \ \simeq \ \ker({}^tg_2)/\imm({}^tg_1) , \]
${}^t\sigma_1^*$ est \`a valeurs dans $\ker({}^tg_2)$ et $\pi$ n'est 
autre que la compos\'ee
\xmat{ F^{**}\ar[r]^-{{}^t\sigma_1^*} & \ker({}^tg_2)\flon[r] &
\ker({}^tg_2)/\imm({}^tg_1) .}
On retrouve de m\^eme $\rho$ \`a partir de $\sigma_1$.

On obtient alors des r\'esolutions localement libres de $\ke$ et $\ke^*$
\xmat{ \cdots F_2\ar[r]^{f_2} & F_1\ar[r]^-{(f_1,\sigma_1)} & 
F_0\oplus G^*\ar[r] & \ke\ar[r] & 0 }
\xmat{ \cdots G_2\ar[r]^{g_2} & G_1\ar[r]^-{(g_1,\sigma_1^*)} & 
G_0\oplus F^*\ar[r] & \ke^*\ar[r] & 0 .}
On a des diagrammes commutatifs
\[
\xymatrix{
F_0\oplus G^*\flon[r]\flon[rd]^{f_0} & \ke\flon[d]^p \\ & F \\
} \ , \ \ \ \
\xymatrix{
F_0\oplus G^*\flon[r] & \ke \\ & G^*\flinc[ul]^-{{}^tg_0}\flinc[u]^i \\
}
\]
et des diagrammes analogues pour $\ke^*$.

On d\'eduit des r\'esolutions pr\'ec\'edentes de $\ke$ et $\ke^*$ la suite
exacte
\xmat{ \cdots F_2\ar[r]^{f_2} & F_1\ar[rr]^-{(f_1,\sigma_1)} & &
F_0\oplus G^*\ar[r]^\theta & 
F^{**}\oplus G_0^*\ar[rr]^-{{}^t\sigma_1^*+{}^tg_1} & & G_1^*\ar[r] & G_2^*
\cdots }
Des diagrammes commutatifs pr\'ec\'edents on d\'eduit que $\theta$ est de la
forme
\[\theta \ = \ \left(\begin{array}{cc}f_0 & 0 \\ -\alpha & {}^tg_0
\end{array}\right) ,\]
avec \ \m{\alpha : F_0\to G_0^*}. De l'exactitude de la suite pr\'ec\'edente
on d\'eduit le diagramme commutatif suivant :
\xmat{
F_2\ar[r]^{f_2}\ar[d] & F_1\ar[r]^{f_1}\ar[d]^{\sigma_1} & 
F_0\ar[r]^{f_0}\ar[d]^{\alpha} & F\ar[r]\ar[d]^{{}^t\sigma_1^*} & 0\ar[d] \\
0\ar[r] & G^*\ar[r]^{{}^tg_0} & G_0^*\ar[r]^{{}^tg_1} & 
G_1^*\ar[r]^{{}^tg_2} & G_2^* }
On obtient de m\^eme le diagramme dual
\xmat{
G_2\ar[r]^{g_2}\ar[d] & G_1\ar[r]^{g_1}\ar[d]^{\sigma_1^*} & 
G_0\ar[r]^{g_0}\ar[d]^{{}^t\alpha} & G\ar[r]\ar[d]^{{}^t\sigma_1} & 0\ar[d] \\
0\ar[r] & F^*\ar[r]^{{}^tf_0} & F_0^*\ar[r]^{{}^tf_1} & 
F_1^*\ar[r]^{{}^tf_2} & F_2^* }
Remarquons que ces diagrammes induisent des morphismes de r\'esolutions de $T$
et $\wT$ :
\xmat{
F_2\ar[r]^{f_2}\ar[d] & F_1\ar[r]^{f_1}\ar[d]^{\sigma_1} & 
F_0\ar[r]^{f_0}\ar[d]^{\alpha} & F^{**}\ar[r]^\pi\ar[d]^{{}^t\sigma_1^*} & 
T\fleq[d] \\
0\ar[r] & G^*\ar[r]^{{}^tg_0} & G_0^*\ar[r]^-{{}^tg_1} & 
\ker({}^tg_2)\ar[r]^-{\pi'} & T }
\xmat{
G_2\ar[r]^{g_2}\ar[d] & G_1\ar[r]^{g_1}\ar[d]^{\sigma_1^*} & 
G_0\ar[r]^{g_0}\ar[d]^{{}^t\alpha} & G^{**}\ar[r]^\rho\ar[d]^{{}^t\sigma_1} 
& \wT\fleq[d] \\
0\ar[r] & F^*\ar[r]^{{}^tf_0} & F_0^*\ar[r]^-{{}^tf_1} & 
\ker({}^tf_2)\ar[r]^-{\rho'} & \wT }

\medskip

Ces r\'esolutions permettent d'expliciter partiellement les isomorphismes
\[\Ext^2(T,W)\simeq\Hom(W^*,\wT), \ \ \ \
\Ext^2(\wT,W)\simeq\Hom(W^*,T) ,\]
$W$ \'etant un fibr\'e vectoriel sur $X$. Soit \ $\mu:F_1\to W$ \ un
morphisme s'annulant sur $\imm(f_2)$. Ce morphisme d\'efinit un \'el\'ement
$u$ de $\Ext^2(T,W)$ d'apr\`es la r\'esolution pr\'ec\'edente de $T$. 
Le morphisme \ \m{{}^t\mu:W^*\to F_1^*} \ est \`a valeurs 
dans \m{\ker({}^tf_2)}, et la composition avec le morphisme \
\m{\ker({}^tf_2)\to\wT} \ du diagramme pr\'ec\'edent donne le morphisme \
\m{W^*\to\wT} \ associ\'e \`a $u$.
\end{sub}

\sepsec


\pagebreak[4]\section{Propri\'et\'es des extensions larges}\label{Larg00}

On utilise ici les notations et r\'esultats de \ref{Larg0}, et
plus particuli\`erement ceux de \ref{Larg_constr}. On consid\`ere donc
une extension large

\[\xymatrix{ (L) \ \ \ \ \ \ \ \
0\ar[r] & G^*\ar[r]^i & \ke\ar[r]^p & F\ar[r] & 0 }\]

et le sujet principal de ce chapitre est l'\'etude de la structure de
$\Ext^1(\ke,\ke)$ induite par la suite exacte pr\'ec\'edente et la suite 
exacte duale.

\sepsub

\Ssect{L'application canonique \ $\End(T)\to\Ext^1(\ke,\ke)$}{end_inc}

Examinons les actions de \ $\Aut(T)=\Aut(\wT)$ \ sur
$\Ext^1(F,G^*)$ et $\Ext^1(G,F^*)$. \hfil\break
Ces actions proviennent des isomorphismes canoniques \
\m{\Ext^1(F,G^*)\simeq\Hom(G^{**},\wT)} \ et \ \m{
\Ext^1(G,F^*)\simeq\Hom(F^{**},T)} \ respectivement.

Soit \ \m{\theta\in\Aut(T)}. Alors on a \
\m{\delta'({}^t\theta\circ\rho) =  \theta^{-1}\circ\delta'(\rho)}, 
\m{\delta(\theta^{-1}\circ\pi) = {}^t\theta\circ\delta(\pi)} .
Il en d\'ecoule que
\m{\delta'({}^t\theta\circ\rho)^{-1}\circ\theta^{-1}\circ\pi 
= \delta'(\rho)^{-1}\circ\pi},
\m{\delta(\theta^{-1}\circ\pi)^{-1}\circ{}^t\theta\circ\rho = 
\delta'(\rho)^{-1}\circ\pi} .
On en d\'eduit que les morphismes $\sigma_1$, $\sigma_1^*$, $\alpha$
correspondant \`a \ \m{(\theta^{-1}\circ\pi,{}^t\theta\circ\rho)} \ sont
les m\^emes que ceux qui correspondent \`a \m{(\pi, \rho)}.

Soient
\[\xi_0 :\End(\wT)\lra\Ext^1(\ke,\ke) , \ \ \ \ \ \
\xi^*_0 :\End(T)\lra\Ext^1(\ke^*,\ke^*)\]
les applications compos\'ees :
\xmat{
\End(\wT)\ar[r]^-{\nu} & 
\Hom(G^{**},\wT)\ar[r]^-{\simeq} & \Ext^1(F,G^*)\ar[r]^-{a} & 
\Ext^1(\ke,G^*)\ar[r]^-{b} & \Ext^1(\ke,\ke)
}
\xmat{
\End(T)\ar[r]^-{\nu'} & 
\Hom(F^{**},T)\ar[r]^-\simeq & \Ext^1(G,F^*)\ar[r]^-{a'} & 
\Ext^1(\ke^*,F^*)\ar[r]^-{b'} & \Ext^1(\ke^*,\ke^*)
}
o\`u $a$, $b$, $a'$, $b'$ sont induits par $(L)$, $(L^*)$, $\nu$ est la
multiplication par $\rho$, et $\nu'$ celle par $\pi$.

\sepprop

\pagebreak[2]\begin{subsub}{\bf Proposition : }\label{id1}
On a \ \m{\xi_0=\xi_0^*}, compte tenu des identifications \ 
\m{\End(T)\simeq\End(\wT)},\hfil\break 
\m{\Ext^1(\ke,\ke)\simeq\Ext^1(\ke^*,\ke^*)}.
\end{subsub}

\begin{proof}
On consid\`ere les morphismes
\[\sigma'_1 : F_1\lra G^*, \ \ \ {\sigma'_1}^* : G_1\lra F^*, \ \ \
\alpha' : F_0\lra G_0^*\]
correspondant \`a
\[\theta\circ\pi : \F\lra T, \ \ \ \ 
\rho : \G\lra \wT .\]
Ce sont les m\^emes que ceux qui correspondent \`a
\[\pi : \F\lra T, \ \ \ \ 
{}^t\theta\circ\rho : \G\lra \wT .\]
L'\'el\'ement $\nu(\theta)$ de $\Ext^1(F,G^*)$ provient de $\sigma'_1$, et
\m{a\circ\nu(\theta)} aussi, compte tenu de la r\'esolution de $\ke$ :
\[\ldots F_2\lra F_1\lra F_0\oplus G^*\lra\ke\lra 0 .\]
Donc \m{b\circ a\circ\nu(\theta)} provient de la compos\'ee
$\xymatrix{ F_1\ar[r]^-{\sigma'_1} & G^*\flinc[r] & \ke  \ .}$
De m\^eme \m{b'\circ a'\circ\nu'(\theta)} provient de la compos\'ee
$\xymatrix{ G_1\ar[r]^-{{\sigma'_1}^*} & F^*\flinc[r] & \ke^*}$.
On obtient alors le diagramme commutatif suivant :
\xmat{
G_1\ar[r]\ar[d]^{{\sigma'_1}^*} & G_0\oplus F^*\ar[r]\ar[dd]^\Theta & 
\ke^*\ar[d]^{{}^ti} \\
F^*\ar[d]^{{}^tp} & & G^{**}\ar[d]^{{}^t\sigma'_1} \\
\ke^*\ar[r] & F_0^*\oplus G^{**}\ar[r] & F_1^*
}
o\`u les ligne du haut et du bas proviennent des r\'esolutions de $\ke$ et
$\ke^*$ vues pr\'ec\'edemment, et $\Theta$ est d\'efini par la matrice
\[\left(\begin{array}{cc} {}^t\alpha'-{}^t\alpha & {}^tf_0 \\
g_0 & 0 \\ \end{array}\right) \ . \]
On en d\'eduit avec la proposition \ref{dual_ext}, que \
\m{b'\circ a'\circ\nu'(\theta)=b\circ a\circ\nu(\theta)}.
\end{proof}

\end{sub}

\sepsub

\Ssect{L'application canonique \ $\Ext^1(\ke,\ke)\to\Ad^2(T)$}{on_ext}

On utilise les notations et r\'esultats de \ref{Larg_constr}. 
L'application \m{\Hom(G^*,T)\to\Ext^1(G^*,F)} d\'eduite de $\pi$ est un
isomorphisme, car \m{\Hom(G^*,F^{**})=\Ext^1(G^*,F^{**})=\nsp}. On note $B$ son
inverse. De m\^eme, l'application \m{\Hom(F^*,\wT)\to\Ext^1(F^*,G)} d\'eduite
de $\rho$ est un isomorphisme. On note $B'$ son inverse.

\sepprop

\pagebreak[2]\begin{subsub}{\bf Lemme : }\label{on_ext1}
Les applications compos\'ees
\xmat{
\Ext^1(\ke,\ke)\ar[r]^-A & \Ext^1(G^*,F)\ar[r]^-B & \Hom(G^*,T)\ar[r]^-\simeq
& \Ext^2(\wT,G^{**})\ar[r]^-C & \Ext^2(\wT,\wT)
}
\xmat{
\Ext^1(\ke^*,\ke^*)\ar[r]^-{A'} & \Ext^1(F^*,G)\ar[r]^-{B'} & 
\Hom(F^*,\wT)\ar[r]^-\simeq & \Ext^2(T,F^{**})\ar[r]^-{C'} & \Ext^2(T,T)
}
o\`u $A$, $A'$ sont induits par $(L)$, $(L')$, $C'$ par $\pi$, et
$C$ par $\rho$ sont \`a valeurs dans \m{\Ad^2(\wT)}, \m{\Ad^2(T)} 
respectivement.
\end{subsub}

\begin{proof}
Dans le premier cas
l'application $A$ est \`a valeurs dans le noyau de la multiplication par
$\sigma$ \m{\Ext^1(G^*,F)\to\Ext^2(F,F)}, et le r\'esultat d\'ecoule de la
proposition \ref{id4}. Le second cas est analogue.
\end{proof}

\sepprop

On note 
\[\xi_2 : \Ext^1(\ke,\ke)\lra\Ad^2(\wT) , \ \ \ \ \ \ 
\xi_2^* : \Ext^1(\ke^*,\ke^*)\lra\Ad^2(T) \]
les applications compos\'ees du lemme \ref{on_ext1}.

\sepprop

\pagebreak[2]\begin{subsub}{\bf Proposition : }\label{id4b}
On a \ \m{\xi_2^*= - \xi_2}, compte tenu des 
identifications \hfil\break \m{\Ext^2(T,T)\simeq\Ext^2(\wT,\wT)},
\m{\Ext^1(\ke,\ke)\simeq\Ext^1(\ke^*,\ke^*)}.
\end{subsub}

\begin{proof}
On utilise les notations de \ref{Larg_constr}.
On choisit comme pr\'ec\'edemment des r\'esolutions localement libres
suffisamment "n\'egatives" :
\xmat{
\cdots F_2\ar[r]^{f_2} & F_1\ar[r]^{f_1} & F_0\ar[r]^{f_0} & F\ar[r] & 0
}
\xmat{
\cdots G_2\ar[r]^{g_2} & G_1\ar[r]^{g_1} & G_0\ar[r]^{g_0} & G\ar[r] & 0
}
\xmat{
\cdots F'_2\ar[r]^{f'_2} & F'_1\ar[r]^{f'_1} & F'_0\ar[r]^{f'_0} & 
F^*\ar[r] & 0
}
\xmat{
\cdots G'_2\ar[r]^{g'_2} & G'_1\ar[r]^{g'_1} & G'_0\ar[r]^{g'_0} & 
G^*\ar[r] & 0
}
de telle sorte que les extensions larges soient construites comme dans
\ref{Larg_constr} et qu'on ait des r\'esolutions localement libres
\xmat{
\cdots F'_2\oplus G_2\ar[r]^{\delta_2^*} & F'_1\oplus G_1\ar[r]^{\delta_1^*} 
& F'_0\oplus G_0\ar[r]^-{\delta_0^*} & \ke^*\ar[r] & 0
}
\xmat{
\cdots G'_2\oplus F_2\ar[r]^{\delta_2} & G'_1\oplus F_1\ar[r]^{\delta_1} 
& G'_0\oplus F_0\ar[r]^-{\delta_0} & \ke\ar[r] & 0
}
\`a partir desquelles les $\Ext^i(\ke,-)$ et $\Ext^i(\ke^*,-)$ puissent
\^etre construits (cf. \ref{Ext_res}). Rappelons que pour $i>0$, $\delta_i^*$
provient d'une matrice
\[\left(\begin{array}{cc}
f'_i & \sigma^{0*}_i \\ 0 & g_i
\end{array}\right) \ , \ \ \ \ 
\sigma^{0*}_i : G_i\lra F'_{i-1},\]
et $\delta_i$ d'une matrice
\[\left(\begin{array}{cc}
g'_i & \sigma^{0}_i \\ 0 & f_i
\end{array}\right) \ , \ \ \ \ 
\sigma^{0*}_i : F_i\lra G'_{i-1}.\]
On a aussi \
\m{\sigma_1 = g'_0\circ\sigma_1^0}, \m{\sigma_1^* = f'_0\circ\sigma_1^{0*}} .

\medskip

{\em \'Etape 1. } On montre d'abord que le diagramme
\xmat{
\Ext^1(\ke,\ke)\fleq[rr]\ar[d]^f & & \Ext^1(\ke^*,\ke^*)\ar[d]^g \\
\Ext^1(G^*,\ke)=\Ext^1(\ke^*,G^{**}) & & \Ext^1(\ke^*,G)\ar[ll]^-{h}
}
est commutatif. Soit \ \m{u\in\Ext^1(\ke,\ke)}, provenant de \
\m{u_1 : G'_1\oplus F_1\lra\ke} .
D'apr\`es \ref{Ext_dual}, on peut supposer qu'il existe un diagramme
commutatif
\[(D_0) \ \ \ \ \
\xymatrix{
G'_1\oplus F_1\ar[r]^{\delta_1}\ar[d]^{u_1} & 
G'_0\oplus F_0\ar[r]^-{\delta_0}\ar[d] & \ke\ar[d]^{u_1^*} \\
\ke\ar[r]^-{{}^t\delta_0^*} & {F'_0}^*\oplus G_0^*\ar[r] & 
{F'_1}^*\oplus G_1^*}\]
et \ \m{{}^tu_1^* : F'_1\oplus G_1\to\ke^*} \ repr\'esente $u$, vu comme
\'el\'ement de $\Ext^1(\ke^*,\ke^*)$. D'autre part, $f(u)$ est 
repr\'esent\'e par \ \m{u_{1\mid G'_1} : G'_1\to F_1}, et $g(u)$ par le
morphisme compos\'e
\xmat{\phi : F'_1\oplus G_1\ar[r]^-{{}^tu_1^*} & \ke^*\ar[r] & G .}
D'apr\`es \ref{Ext_dual}, pour v\'erifier que \ \m{h\circ g(u)=f(u)} \ il
suffit de trouver un diagramme commutatif
\[\xymatrix{
G'_1\ar[r]^{g'_1}\ar[d]_{u_{1\mid G'_1}} & G'_0\ar[d]\ar[r]^{g'_0} &
G^*\ar[d]^{{}^t\phi} \\
\ke\ar[r]^-{{}^t\delta_0^*} & {F'_0}^*\oplus G_0^*\ar[r]^{{}^t\delta_1^*} &
{F'_1}^*\oplus G_1^*
} \ .\]
Mais un tel diagramme d\'ecoule ais\'ement du diagramme pr\'ec\'edent 
$(D_0)$.

Il suffit donc de montrer que le diagramme
\xmat{
\Ext^1(G^*,\ke)\ar[d]^{A_0} & & \Ext^1(\ke^*,G)\ar[ll]^h\ar[d]^{A'_0} \\
\Ext^1(G^*,F)\ar[d]^B & & \Ext^1(F^*,G)\ar[d]^{B'} \\
\Hom(G^*,T)\ar[d]^\simeq & & \Hom(F^*,\wT)\ar[d]^\simeq \\
\Ext^2(\wT,G^{**})\ar[d]^C & & \Ext^2(T,F^{**})\ar[d]^{C'} \\
\Ext^2(\wT,\wT)\fleq[rr] & & \Ext^2(T,T)
}
o\`u les fl\`eches proviennent de $(L)$, $(L^*)$, $\pi$ et $\rho$ est 
commutatif.

\medskip

{\em \'Etape 2 .} On part de \ $v\in\Ext^1(\ke^*,G)$, repr\'esent\'e par
un morphisme
\[(\lambda,\gamma) : F'_1\oplus G_1\lra G\]
s'annulant sur $\imm(\delta_2^*)$. D'apr\`es \ref{Ext_dual}, on peut supposer
qu'il existe un diagramme commutatif
\[(D_1) \ \ \ \ \xymatrix{ 
F'_1\oplus G_1\ar[rr]^{\delta_1^*}\ar[d]^{(\lambda,\gamma)} & &
F'_0\oplus G_0\ar[rr]^-{\delta_0^*}\ar[d]^{(\epsilon_0,\chi_0)} & &
\ke^*\ar[d]^\theta \\
G^{**}\ar[rr]^-{{}^tg'_0} & & {G'_0}^*\ar[rr]^-{{}^tg'_1} & & {G'_1}^*
}\]
et \ \m{{}^t\theta : G'_1\to\ke} \ s'annule sur $\imm(g'_2)$ et repr\'esente
$h(v)$.

D'autre part, \m{\lambda : F'_1\to G^{**}} \ s'annule sur $\imm(f'_2)$ et
repr\'esente $A'_0(v)$, $A_0\circ h(v)$ est repr\'esent\'e par le morphisme
compos\'e \
\m{\xymatrix{\mu : G'_1\ar[r]^{{}^t\theta} & \ke\ar[r]^p & G^{**}}} .

\medskip

{\em \'Etape 3 . } Puisque \ $\Ext^1(F^*,G^*)=\Ext^1(G^*,F^*)=\nsp$, il existe
des diagrammes commutatifs \'etendant $\lambda$ et $\mu$
\xmat{
F'_1\ar[r]^{f'_1}\ar[d]^\lambda & F'_0\ar[r]^{f'_0}\ar[d]^{\lambda_0}
& F^*\ar[d]^{\ov{\lambda}} \\
G\flinc[r] & G^{**}\ar[r]^\rho & \wT
}
\xmat{
G'_1\ar[r]^{g'_1}\ar[d]^\mu & G'_0\ar[r]^{g'_0}\ar[d]^{\mu_0}
& G^*\ar[d]^{\ov{\mu}} \\
F\flinc[r] & F^{**}\ar[r]^\pi & T
}
et on a \
\m{\ov{\lambda} = B'\circ A'_0(v)}, \m{\ov{\mu} = B\circ A_0\circ h(v)} .
D'apr\`es $(D_1)$ et les deux diagrammes\hfil\break pr\'ec\'edents, on a \
\m{\epsilon_0f'_1 = {}^tg'_0\lambda = {}^tg'_0\lambda_0f'_1} .
Donc \ \m{(\epsilon_0-{}^tg'_0\lambda_0)f'_1=0}. Il existe donc un
morphisme \ \m{\psi : F^*\to {G'_0}^*} \ tel que \ 
\m{\epsilon_0-{}^tg'_0\lambda_0=\psi f'_0.} On peut alors prendre \ 
\m{\mu_0={}^t\psi}, car \
\m{{}^tf'_0{}^t\psi g'_1 = ({}^t\epsilon_0-{}^t\lambda_0g'_0)g'_1 =
{}^t\epsilon_0g'_1 = {}^t\mu} .
On peut donc supposer qu'on a \
\m{\epsilon_0 \ = \ {}^tg'_0\lambda_0 + {}^t\mu_0f'_0}.

\medskip

{\em \'Etape 4 . } On utilise maintenant les r\'esolutions canoniques de 
$T$ et $\wT$ (cf. la fin de \ref{Larg_constr}). On va montrer qu'il existe
un morphisme \ \m{\mu' : G^*\lra \ker({}^tg_2)} \
tel que \ \m{\pi'\circ\mu'=\ov{\mu}} . 

D'apr\`es la construction des extensions (cf. \ref{Extensconstr}) il existe
un diagramme commutatif avec lignes exactes
\xmat{
G_1\ar[r]^{g_1}\ar[d]^{\sigma_1^*} & G_0\ar[r]^{g_0}\ar[d] & G\fleq[d] \\
F^*\ar[r] & \ke^*\ar[r] & G
}
En le dualisant on obtient le diagramme commutatif avec lignes exactes suivant :
\xmat{
0\ar[r] & G^*\ar[r]^i\fleq[d] & \ke\ar[r]^p\ar[d] & 
F^{**}\ar[r]^\pi\ar[d]^{{}^t\sigma_1^*} & T\ar[r]\fleq[d] & 0 \\
0\ar[r] & G^*\ar[r]^{{}^tg_0} & G_0^*\ar[r]^-{{}^tg_1} & 
\ker({}^tg_2)\ar[r]^-{\pi'} & T\ar[r] & 0
}
Cela prouve que les deux suites exactes de ce diagramme induisent le m\^eme
\'el\'ement de $\Ext^2(T,G^*)$ (c'est en fait \'evidemment $\rho$). Ceci
prouve que le diagramme suivant, o\`u les fl\`eches sont induites par le
diagramme pr\'ec\'edent, est commutatif :
\xmat{
\Hom(G^*,T)\ar[r]^-\nu\ar[d]^{B^{-1}} & \Ext^1(G^*,G_0^*/G^*)\ar[d]^{\nu'} \\
\Ext^1(G^*,F)\ar[r]^\eta & \Ext^2(G^*,G^*)
}
On a \ \m{\eta(B^{-1}(\ov{\mu}))=0} (car $B^{-1}(\ov{\mu})$ provient de 
$\Ext^1(G^*,\ke)$). Donc on a \ \m{\nu'\circ\nu(\ov{\mu})=0}. Mais $\nu'$
est injective, car \ \m{\Ext^1(G^*,G_0^*)=\nsp} (les r\'esolutions sont
suffisamment "n\'egatives"). Donc \ \m{\nu(\ov{\mu})=0}. Ceci prouve 
l'existence de $\mu'$.

On montre de m\^eme qu'il existe un morphisme \
\m{\lambda' : F^*\lra F_1^*} \
tel que \ \m{\rho'\circ\lambda'=\ov{\lambda}} . Ici la situation est
l\'eg\`erement diff\'erente, car ${B'}^{-1}(\ov{\lambda})$ provient de
$\Ext^1(\ke^*,G)$, et il faudrait qu'il provienne de $\Ext^1(F^*,\ke^*)$.
Mais on a un diagramme commutatif
\xmat{
\Ext^1(\ke^*,\ke^*)\ar[r]^{g'}\ar[d]^g & \Ext^1(F^*,\ke^*)\ar[d] \\
\Ext^1(\ke^*,G)\ar[r] & \Ext^1(F^*,G)
}
o\`u $g$ est surjective d'apr\`es la proposition \ref{larg_lemm2}.
Donc si ${B'}^{-1}(\ov{\lambda})$ provient de $\Ext^1(\ke^*,G)$, il provient
aussi de $\Ext^1(F^*,\ke^*)$.

\medskip

{\em \'Etape 5 .} Il d\'ecoule de ce qui pr\'ec\`ede que \
\m{C\circ B\circ A_0\circ h(v)} \ est repr\'esent\'e par \hfil\break
\m{\mu' : G^*\lra \ker({}^tg_2)} \
compte tenu de la r\'esolution de $T$
\[\xymatrix{
0\ar[r] & G^*\ar[r]^{{}^tg_0} & G_0^*\ar[r]^{{}^tg_1} & 
\ker({}^tg_2)\ar[r]^{\pi'} & T\ar[r] & 0
} \ , \]
et \ \m{C'\circ B'\circ A'_0(v)} \ est repr\'esent\'e par \
\m{\lambda' : F^*\lra F_1^*} \ compte tenu de la r\'esolution de $T$
\[\xymatrix{
\cdots F_2\ar[r]^{f_2} & F_1\ar[r]^{f_1} & F_0\ar[r]^{F_0} & 
F^{**}\ar[r]^\pi & T\ar[r] & 0
} \ . \]
On utilise maintenant le diagramme commutatif
d\'efini \`a la fin de \ref{Larg_constr}
\xmat{
F_2\ar[r]^{f_2}\ar[d] & F_1\ar[r]^{f_1}\ar[d]^{\sigma_1} & 
F_0\ar[r]^{f_0}\ar[d]^{\alpha} & F^{**}\ar[r]^\pi\ar[d]^{{}^t\sigma_1^*} & 
T\fleq[d] \\
0\ar[r] & G^*\ar[r]^{{}^tg_0} & G_0^*\ar[r]^-{{}^tg_1} & 
\ker({}^tg_2)\ar[r]^-{\pi'} & T }
Il suffit donc de trouver des morphismes \ \m{H : F_1\to G_0^*}, 
\m{K : F_0\to G_1^*} \ tels que \hfil\break
\m{{}^t\sigma_1^*{}^t\lambda' + \mu'\sigma_1 = Kf_1+{}^tg_1H} , 
ou ce qui revient au m\^eme \
\m{\lambda'\sigma_1^* + {}^t\sigma_1{}^t\mu' = {}^tf_1{}^tK+{}^tHg_1} .

\medskip

{\em \'Etape 6 . } On a
\begin{eqnarray*}
\rho'({}^t\sigma_1^0({}^tg'_0\lambda_0){\sigma_1^0}^* - 
\lambda'\sigma_1^*) & = & \rho'({}^t\sigma_1\lambda_0 -
\lambda'f'_0){\sigma_1^0}^* \\
& = & (\rho\lambda_0 - \ov{\lambda}f'_0){\sigma_1^0}^*\\
& = & 0 .
\end{eqnarray*}
Donc \ \m{{}^t\sigma_1^0({}^tg'_0\lambda_0){\sigma_1^0}^* - 
\lambda'\sigma_1^*} \ est \`a valeurs dans \ \m{\ker(\rho')=\imm({}^tf_1)}.
Puisque \ \m{\Ext^1(G_1,F^*)=\nsp}, on peut \'ecrire
\[ {}^t\sigma_1^0({}^tg'_0\lambda_0){\sigma_1^0}^* - 
\lambda'\sigma_1^* \ = \ {}^tf_1X ,\]
avec \ \m{X : G_1\to F_0^*}. De m\^eme, on a
\[\pi'({}^t{\sigma_1^0}^*({}^tf'_0\mu_0)\sigma_1^0 -
\mu'\sigma_1) \ = \ 0 , \]
donc puisque \ \m{\Ext^1(F_1,G^*)=\nsp}, on peut \'ecrire
\[ {}^t{\sigma_1^0}^*({}^tf'_0\mu_0)\sigma_1^0 -
\mu'\sigma_1 \ = \ {}^tg_1Y ,\]
avec \ \m{Y : F_1\to G_0^*}, c'est-\`a-dire
\[ {}^t\sigma_1^0({}^t\mu_0f'_0){\sigma_1^0}^* -
{}^t\sigma_1{}^t\mu' \ = \ {}^tYg_1 .\]
On a donc
\[ \lambda'\sigma_1^* + {}^t\sigma_1{}^t\mu' \ = \ 
{}^t\sigma_1^0\epsilon_0{\sigma_1^0}^* - {}^tf_1X - {}^tYg_1 \]
(cf. la fin de l'\'etape 3). 

Il suffit donc de montrer qu'il existe des morphismes \ \m{H' : G_0\to F_1^*}
\ et \ \m{K' : G_1\to F_0^*} \ tels que
\[ {}^t\sigma_1^0\epsilon_0{\sigma_1^0}^* \ = \
{}^tf_1K' + H'g_1 . \]
On utilise maintenant le diagramme commutatif $(D_1)$ et celui de la fin de
l'\'etape 5. Puisque \ \m{\Ext^1(G_1,\ker(g_0))=\nsp}, il existe un morphisme
\ \m{\gamma_0 : G_1\to G_0} \ tel que \ \m{\gamma = g_0\gamma_0}. On a
alors
\begin{eqnarray*}
{}^t\sigma_1^0\epsilon_0{\sigma_1^0}^* - 
{}^tf_1\alpha\gamma_0 & = & {}^t\sigma_1^0\epsilon_0{\sigma_1^0}^* -
{}^t\sigma_1g_0\gamma_0 \\
& = & {}^t\sigma_1^0\epsilon_0{\sigma_1^0}^* -
{}^t\sigma_1\gamma \\
& = & {}^t\sigma_1^0(\epsilon_0{\sigma_1^0}^* - 
{}^tg'_0\gamma) \\
& = & -{}^t\sigma_1^0\chi_0g_1
\end{eqnarray*}
(d'apr\`es le carr\'e commutatif de gauche de $(D_1)$). On obtient
finalement
\[ {}^t\sigma_1^0\epsilon_0{\sigma_1^0}^* \ = \
{}^tf_1\alpha\gamma_0 - {}^t\sigma_1^0\chi_0g_1 ,\]
ce qui d\'emontre la proposition \ref{id4b}.
\end{proof}
\end{sub}

\sepsub

\Ssect{Diagrammes 3x3 canoniques}{s_ext_0}

On utilise les notations et r\'esultats de \ref{Larg_constr}. 
 On s'int\'eresse \`a
la structure induite sur $\Ext^1(\ke,\ke)$ par les suites exactes $(L)$ et
$(L')$.

On note $A_2(\sigma)$ le noyau de la multiplication \`a gauche par $\sigma$
\[\sigma\times \ : \Ext^1(G^*,F)\lra\Ext^2(F,F) .\]
Si $X$ est une surface
$A_2(\sigma)$ est aussi le noyau de la
multiplication \`a droite par $\sigma$
\[\times\sigma \ : \Ext^1(G^*,F)\lra\Ext^2(G^*,G^*) \]
\`a cause du diagramme anticommutatif
\xmat{
\Ext^1(G^*,F)\ar[rr]\ar[d] &  & \Ext^2(F,F)\ar[d]^{tr_2(F)} \\
\Ext^2(G^*,G^*)\ar[rr]^-{tr_2(G^*)} & & H^2(\ko_X)
}
et du fait que $tr_2(F)$ et $tr_2(G^*)$ sont des isomorphismes.

On note $B_2(\sigma)$ le noyau de la multiplication \`a gauche par $\sigma$
\[\sigma\times \ : \Ext^1(G^*,\ke)\lra\Ext^2(F,\ke) .\]

\sepprop

\pagebreak[2]\begin{subsub}{\bf Lemme : }\label{larg_lemm2_0}
Si $X$ est une surface, on a \ \m{B_2(\sigma)=\Ext^1(G^*,\ke)}.
\end{subsub}

\begin{proof}
On a un diagramme commutatif
\xmat{
\Ext^2(F,\ke)\ar[rr]^-g\ar[d]^f & & \Ext^2(\ke,\ke)\ar[d]^{tr_2(\ke)} \\
\Ext^2(F,F)\ar[rr]^{tr_2(F)} & & H^2(\ko_X)
}
L'application $f$ est un isomorphisme car \ 
\m{\Ext^2(F,G^*)=\nsp} et
\m{tr_2(F)} est un isomorphisme. On en d\'eduit
que $g$ est injective, et le lemme d\'ecoule du fait que le noyau de $g$ est
l'image de la multiplication \`a gauche par $\sigma$ \
\m{\Ext^1(G^*,\ke)\to\Ext^2(F,\ke)} .
\end{proof}

\sepprop

\pagebreak[2]\begin{subsub}{\bf Proposition : }\label{larg_lemm2}
1 - L'application canonique 
\xmat{
\Hom(F,G^*)\oplus\C\ar[r] & \End(\ke) }
\[{ \ \ \ \ \ \ \ \ \ \ \ \xymatrix{
(f,t)\ar@{|-_{>}}[r] & i\circ f\circ p+tI_\ke \\ }}\]
est un isomorphisme, et \ \m{\Hom(G^*,\ke)=\Hom(\ke,F)=\C}.

2 - On a un diagramme 3x3 canonique
\xmat{
& 0\ar[d] & 0\ar[d] & 0\ar[d] \\
0\ar[r] & \Ext^1(F,G^*)/\C\sigma\ar[r]\ar[d] & \Ext^1(\ke,G^*)\ar[d]\ar[r] & 
\Ext^1(G^*,G^*)\ar[d]\ar[r] & 0 \\
0\ar[r] & \Ext^1(F,\ke)\ar[r]\ar[d] & \Ext^1(\ke,\ke)\ar[d]\ar[r] & 
B_2(\sigma)\ar[d]\ar[r] & 0 \\
0\ar[r] & \Ext^1(F,F)\ar[r]\ar[d] & \Ext^1(\ke,F)\ar[d]\ar[r] & 
A_2(\sigma)\ar[d]\ar[r] & 0 \\
& 0 & 0 & 0 & }

3 - Si $X$ est une surface, on a une suite exacte 
\xmat{
0\ar[r] & H^2(\ko_X)\ar[r]^{j} & \Ext^2(\ke,\ke)\ar[r]^q
& \Ext^2(G^*,F^{**})\ar[r] & 0 ,
}
$j$ \'etant induit par l'inclusion \ \m{\ko_X\subset\ke^*\ot\ke}, et $q$
par $(L)$. Cette suite exacte se scinde naturellement.
\end{subsub}

\begin{proof} L'assertion 1 est imm\'ediate. 
D\'emontrons 2. 
En consid\'erant le
diagramme commutatif
\xmat{
\Ext^2(\ke,G^*)\ar[rr]^k\ar[d]^h & & \Ext^2(\ke,\ke)\ar[d]^{tr_2(\ke)} \\
\Ext^2(G^*,G^*)\ar[rr]^{tr_2(G^*)} & & H^2(\ko_X)
}
on montre comme dans le lemme \ref{larg_lemm2_0} que le morphisme \
\m{\Ext^1(\ke,\ke)\to\Ext^1(\ke,F)} \ est
surjectif, en utilisant la lissit\'e de \m{G^*}. Le diagramme 3x3 d\'ecoule
alors du diagramme commutatif suivant
\xmat{
& 0\ar[d] \\
\Ext^1(\ke,G^*)\ar[r]\ar[d] & \Ext^1(G^*,G^*)\ar[r]\ar[d] & 0\ar[d]\\
\Ext^1(\ke,\ke)\ar[r]\ar[d] & \Ext^1(G^*,\ke)\ar[r]\ar[d] & 
\Ext^2(F,\ke)\ar[d]\\
\Ext^1(\ke,F)\ar[d]\ar[r] & \Ext^1(G^*,F)\ar[d]\ar[r] & \Ext^2(F,F)\ar[d]\\
0 & \Ext^2(G^*,G^*)\fleq[r] & H^2(\ko_X)
}

L'assertion 3 est imm\'ediate.

\end{proof}

\sepprop

Posons
\[M=\Ext^1(\ke,G^*) \ \subset \ N=\Ext^1(G,\ke^*) \ \subset 
\Ext^1(\ke,\ke) ,\]
\[M_*=\Ext^1(\ke^*,F^*) \ \subset \ N_*=\Ext^1(F,\ke) \ \subset 
\Ext^1(\ke,\ke) .\]

D'apr\`es la proposition \ref{larg_lemm2}, ces sous-espaces vectoriels 
s'ins\`erent dans les diagrammes 3x3 canoniques suivants
\Bdiag{\Ext^1(F,G^*)/\C\sigma}{M}{\Ext^1(G^*,G^*)}
{N_*}{\Ext^1(\ke,\ke)}{B_2(\sigma)}
{\Ext^1(F,F)}{\Ext^1(\ke,F)}{A_2(\sigma)}
\Bdiag{\Ext^1(G,F^*)/\C\sigma^*}{M_*}{\Ext^1(F^*,F^*)}
{N}{\Ext^1(\ke^*,\ke^*)}{B_2(\sigma^*)}
{\Ext^1(G,G)}{\Ext^1(\ke^*,G)}{A_2(\sigma^*)}

\sepprop

\pagebreak[2]\begin{subsub}{\bf Lemme : }\label{struc_ext_1}
Le quotient $N/M$ (resp. $N_*/M_*$) s'identifie naturellement à un sous-espace
vectoriel de \m{\Hom(\ke,T)} (resp. \m{\Hom(\ke^*,\wT)}), qui est 
\m{\Hom(\ke,T)} (resp. \m{\Hom(\ke^*,\wT)}) tout entier si $X$ est une surface.
On a des isomorphismes canoniques
\[ M\cap N_* \ \simeq \ \Ext^1(F,G^*)/\C\sigma \ \simeq \ 
\Hom(G^{**},\wT)/\C\rho, \] 
\[ M_*\cap N \ \simeq \Ext^1(G,F^*)/\C\sigma^* \ \simeq \ 
\Hom(F^{**},T)/\C\pi ,\]
\[ \Ext^1(\ke,\ke)/(M+N_*) \ \simeq \ A_2(\sigma) , \]
\[ \Ext^1(\ke,\ke)/(M_*+N) \ \simeq \ A_2(\sigma^*) . \]
\end{subsub}

\begin{proof} Montrons que \ \m{N/M\subset\Hom(\ke,T)}. De la suite
exacte 
\[0\lra G\lra G^{**}\lra\wT\lra 0\]
on d\'eduit la suite exacte
\xmat{
\Ext^1(G^{**},\ke^*)\flinc[r] & \Ext^1(G,\ke^*)\ar[r] 
& \Ext^2(\wT,\ke^*)\ar[r] &\Ext^2(G^{**},\ke^*)\ar[r]^f & \Ext^2(G,\ke^*)}
d'o\`u l'inclusion d\'ecoule. Si $X$ est une surface, il faut
montrer que $f$ est injective. En utilisant le fait que
\[\Ext^i(G^{**},F^*)=\nsp \ \ \ \ {\rm si} \ i\geq 1\]
et que  $T$ est de dimension 0, on voit que les applications canoniques
\[\Ext^2(G^{**},\ke^*)\lra \Ext^2(G^{**},G), \ \ \ \
\Ext^2(G^{**},G)\lra \Ext^2(G^{**},G^{**})\]
induites par $(L)$, $(L^*)$ sont des isomorphismes. Mais puisque 
$G^{**}$ est 2-lisse, l'application trace $tr_2(G^{**})$ est
est un isomorphisme. Le fait que $f$ est injective d\'ecoule donc du
diagramme commutatif suivant :
\xmat{
\Ext^2(G^{**},\ke^*)\ar[r]\ar[dd]^f & 
\Ext^2(G^{**},G^{**})\ar[drr]^{tr_2(G^{**})}\\
 & & & H^2(\ko_X)\\
 \Ext^2(G,\ke^*)\ar[r] &\Ext^2(G,G)\ar[urr]^{tr_2(G)}
}
Le cas de \m{N_*/M_*} \ est semblable.

Les autres isomorphismes d\'ecoulent des diagrammes pr\'ec\'edents.
\end{proof}

\end{sub}

\sepsub

\Ssect{Quelques diagrammes commutatifs ou anticommutatifs}{s_ext4b}

Les d\'emonstrations de quelques uns des r\'esultats suivants sont omises. Elles
sont analogues \`a celles des propositions \ref{id1} et \ref{id4b}.

\sepprop

\pagebreak[2]\begin{subsub}{\bf Proposition : }\label{id2}
On consid\`ere le diagramme
\xmat{
\Ext^2(T,G^*)\ar[r]^-\simeq & \Hom(G^{**},\wT)\ar[dd]^{\delta_0} \\
\Ext^1(F,G^*)\ar[u]^{\delta_1}\ar[d]^a & \\
\Ext^1(\ke,G^*)=\Ext^1(G^{**},\ke^*)\ar[r]^-{b} & \Ext^1(G^{**},G) }
o\`u $\delta_0$, $\delta_1$ sont les morphismes de liaison provenant 
des suites exactes longues provenant des suites exactes 
\[0\to G\to G^{**}\to \wT\to 0, \ \ \ \ \
0\to F\to F^{**}\to T\to 0\]
respectivement, $a$ \'etant induit par le morphisme \ $\ke\to F$ \ et $b$ 
par le morphisme \ $\ke^*\to G$. 

Alors ce diagramme est anticommutatif.
\end{subsub}

\sepprop

On a bien entendu un diagramme anti-commutatif analogue "dual" de celui de la
proposition \ref{id2} :
\xmat{
\Ext^2(\wT,F^*)\ar[r]^-\simeq & \Hom(F^{**},T)\ar[dd] \\
\Ext^1(G,F^*)\ar[u]\ar[d] & \\
\Ext^1(\ke^*,F^*)=\Ext^1(F^{**},\ke)\ar[r] & \Ext^1(F^{**},F) }

On utilisera en \ref{s_ext4} la

\sepprop

\pagebreak[2]\begin{subsub}{\bf Proposition : }\label{id3}
On consid\`ere le diagramme
\xmat{
\Ext^1(F,\ke)\ar[r]^-c\ar[d]^d
& \Ext^1(\ke,\ke)\simeq\Ext^1(\ke^*,\ke^*)\ar[r]^-{A'}
& \Ext^1(F^*,G)\ar[r]^-{B'} & \Hom(F^*,\wT)\ar[d]^\simeq \\
\Ext^1(F,F)\ar[r]^f & \Ext^2(T,F)\ar[rr]^g & & \Ext^2(T,F^{**})
}
o\`u $c$ et $d$ proviennent de $(L)$ et $f$ de $\pi$. Alors ce diagramme est 
commutatif.
\end{subsub}

\sepprop

On a bien entendu un diagramme commutatif analogue "dual" de celui de la
proposition \hbox{\ref{id3} :}
\xmat{
\Ext^1(G,\ke^*)\ar[r]\ar[d]
& \Ext^1(\ke^*,\ke^*)\simeq\Ext^1(\ke,\ke)\ar[r]
& \Ext^1(G^*,F)\ar[r] & \Hom(G^*,T)\ar[d]^\simeq \\
\Ext^1(G,G)\ar[r] & \Ext^2(\wT,G)\ar[rr] & & \Ext^2(\wT,G^{**})
}

\sepprop

\pagebreak[2]\begin{subsub}{\bf Proposition : }\label{id4c}
Le diagramme canonique
\xmat{
\Hom(G^*,T)\ar[r]^-{C_0}\ar[d]^\simeq & \Ext^1(F,T)\ar[dd]^{D_0} \\
\Ext^2(\wT,G^{**})\ar[d]^C \\
\Ext^2(\wT,\wT)\fleq[r] & \Ext^2(T,T)
}
(o\`u $C_0$ provient de $(L)$ et $D_0$ de $\pi$) est commutatif.
\end{subsub}

\begin{proof}
On peut d\'eduire ce r\'esultat de la d\'emonstration de la proposition 
\ref{id4b} ou le d\'emontrer directement de mani\`ere analogue.
\end{proof}

\sepprop

\pagebreak[2]\begin{subsub}{\bf Proposition : }\label{surjx2}
Si $X$ est une surface, alors les applications $\xi_2$ et $\xi_2^*$ sont 
surjectives.
\end{subsub}

\begin{proof}
Il suffit de prouver que $\xi_2$ est surjective. Compte tenu de la 
d\'efinition
de $\xi_2$ dans \ref{on_ext}, cela d\'ecoule du fait que l'image de $A$ est
\m{A_2(\sigma)}, $B$ est un isomorphisme, et $C$ est surjective, car \
\m{\coker(C)\subset\Ext^3(\wT,G)}.
\end{proof}

\sepprop

\pagebreak[2]\begin{subsub}{\bf Proposition : }\label{id4}
Le diagramme suivant
\xmat{
\Ext^1(G^*,F)\ar[rr]^{\times\sigma}& & \Ext^2(F,F)\ar[rrd]^-{tr_2(F)} \\
\Hom(G^*,T)\ar[u]^{B^{-1}}\ar[d]^\simeq & & & & H^2(\ko_X) \\
\Ext^2(\wT,G^{**})\ar[rr]^{C} & & 
\Ext^2(\wT,\wT)\ar[rru]^{tr_2(\wT)}
}
o\`u $\times\sigma$ est la multiplication par $\sigma$ est anticommutatif.
\end{subsub}

\begin{proof}
Cela d\'ecoule du diagramme commutatif
\[\xymatrix{
\Ext^1(G^*,F)\ar[r] & \Ext^2(F,F) \\
\Hom(G^*,T)\ar[u]^\simeq\ar[r] & \Ext^1(F,T)\ar[u]
} \ \ ,\]
du diagramme anticommutatif
\[\xymatrix{
\Ext^1(F,T)\ar[r]\ar[d]^\simeq & \Ext^2(F,F)\ar[d]^{tr_2(F)}\\
\Ext^2(T,T)\ar[r]^-{tr_2(T)} & H^2(\ko_X)
} \ \ ,\]
de la compatibilit\'e de l'isomorphisme \ 
\m{\Ext^2(T,T)\simeq\Ext^2(\wT,\wT)} \ avec la trace, et du diagramme
commutatif canonique
\xmat{
\Ext^2(T,G^*)\ot\Hom(G^*,T)\fleq[d]\ar[r] & \Ext^2(T,T)\fleq[d] \\
\Hom(G^{**},\wT)\ot\Ext^2(\wT,G^{**})\ar[r] & \Ext^2(\wT,\wT)
}
(cf. la fin de \ref{grat}).
\end{proof}

\sepprop

\pagebreak[2]\begin{subsub}{\bf Corollaire : }\label{coro_id4}
Le diagramme suivant
\xmat{
\Ext^2(\wT,G^{**})\ar[r]^{\ov{\rho}} & \Ext^2(G^{**},G^{**})\fleq[dd] \\
\Hom(G^*,T)\ar[u]^\simeq\ar[d]^\simeq \\
\Ext^1(G^*,F)\ar[r]^{\delta_{G^*}} & \Ext^2(G^*,G^*)
}
o\`u $\ov{\rho}$ provient de $\rho$ et $\delta_{G^*}$ de $(L)$ est
commutatif.
\end{subsub}

\begin{proof}
Cela d\'ecoule de la proposition \ref{id4}, des diagrammes
\[
\xymatrix{
\Ext^1(G^*,F)\ar[r]^{\delta_{G^*}}\ar[d]^{\delta_F} & 
\Ext^2(G^*,G^*)\ar[d]^{tr_2(G^*)} \\
\Ext^2(F,F)\ar[r]^{tr_2(F)} & H^2(\ko_X)
} \ \ \ \
\xymatrix{
\Ext^2(\wT,G^{**})\ar[r]^{\ov{\rho}}\ar[d]^C 
& \Ext^2(G^{**},G^{**})\ar[d]^{tr_2(G^{**})} \\
\Ext^2(\wT,\wT)\ar[r]^{tr_2(\wT)} & H^2(\ko_X)
}
\]
dont le premier est anticommutatif et le second commutatif,
et de la lissit\'e de $G^*$.
\end{proof}

\end{sub}

\sepsub

\Ssect{Sous-espaces vectoriels canoniques de $\Ext^1(\ke,\ke)$}{s_ext4}

On a des morphismes canoniques injectifs \'evidents
\[\End(T)/\C I_{T} \lra \Hom(F^{**},T)/\C\pi , \ \ \ \
\End(\wT)/\C I_\wT \lra \Hom(G^{**},\wT)/\C\rho ,\]
et $\End(T)$ et $\End(\wT)$ sont canoniquement isomorphes. On verra 
donc \m{\End(T)} comme un sous-espace
vectoriel de \m{\Hom(F^{**},T)} et \m{\Hom(G^{**},\wT)}, et
\m{\End(T)/\C} comme un sous-espace vectoriel de $M\cap N_*$ et $M_*\cap N$.
Le fait que ces deux sous-espaces vectoriels co\" \i ncident dans
\m{\Ext^1(\ke,\ke)} d\'ecoule de la proposition \ref{id1}.

\sepprop

\pagebreak[2]\begin{subsub}{\bf Proposition : }\label{struc_ext_2}
1 - Le diagramme canonique 
\xmat{
M=\Ext^1(\ke,G^*) \flinc[rd]\ar[ddd]^\simeq &           & \\
                   & N=\Ext^1(G,\ke^*)\ar[rd] & \\ 
                   &           & \Ext^1(G,G) \\
\Ext^1(G^{**},G)\ar[urr] }
(ou la fl\`eche verticale provient de \ \m{{}^ti:\ke^*\to G}) est commutatif. 
De m\^eme, le diagramme cano-\break nique
\xmat{
M_*=\Ext^1(\ke^*,F^*) \flinc[rd]\ar[ddd]^\simeq &           & \\
                   & N_*=\Ext^1(F,\ke)\ar[rd] & \\ 
                   &           & \Ext^1(F,F) \\
\Ext^1(F^{**},F)\ar[urr] }
(ou la fl\`eche verticale provient de \ \m{p:\ke\to F}) est commutatif.

2 - Soit \ \m{Q\subset N} l'image r\'eciproque (par le morphisme canonique
\ \m{N\to\Ext^1(G,G)}) de l'image de $\Ext^1(G^{**},G)$ dans $\Ext^1(G,G)$.
Alors on a \ \m{M\subset Q}, \m{M_*\cap Q=M_*\cap N} \ et \
\m{Q=M+M_*\cap N}.

De m\^eme, Soit \ \m{Q_*\subset N_*} l'image r\'eciproque (par le morphisme 
canonique
\ \m{N_*\to\Ext^1(F,F)}) de l'image de $\Ext^1(F^{**},F)$ dans $\Ext^1(F,F)$.
Alors on a \ \m{M_*\subset Q_*}, \m{M\cap Q_*=M\cap N_*} \ et \
\m{Q_*=M_*+M\cap N_*}.

3 - On a \ \m{M\cap M_* = \End(T)/\C} \ et \ \m{Q+M_* = Q_*+M = M+M_*}.

4 - On a \ \m{Q\cap Q_* = M_*\cap N + M\cap N_*} .
\end{subsub}

\begin{proof} 1- d\'ecoule ais\'ement de l'associativit\'e des Ext. 

D\'emontrons 2-. Il est clair que $M\subset Q$. La seconde assertion 
d\'ecoule
du fait que le noyau de \ \m{N\to\Ext^1(G,G)} \ est \m{M_*\cap N}, et la
troisi\`eme de 1-.

L'assertion sur $Q_*$ se d\'emontre de la m\^eme fa\c con.

La derni\`ere assertion de 3- d\'ecoule imm\'ediatement de 2-. Il reste \`a
montrer que \hfil\break
 \m{M\cap M_* = \End(T)/\C}. D'apr\`es le lemme \ref{reg_3x3} et ce
qui pr\'ec\`ede on a un diagramme commutatif avec lignes exactes
\xmat{
0\ar[r] & \End(T)/\C\ar[d]^j\ar[r] & M\flinc[d]\flon[r] & 
\Ext^1(G^{**},G)/(\End(T)/\C)\ar[r]\fleq[d] & 0 \\
0\ar[r] & M_*\cap N\ar[r]\fleq[d] & Q\flinc[d]\flon[r] &
\Ext^1(G^{**},G)/(\End(T)/\C)\ar[r]\flinc[d] & 0 \\
0\ar[r] & M_*\cap N\ar[r] & N\flon[r] & \Ext^1(G,G)\ar[r] & 0
}
o\`u $j$ est injectif. D'apr\`es la proposition \ref{id1} et la proposition 
\ref{id2}, pour tout \ \m{\theta\in\End(T)}, si $\ov{\theta}$ est l'image de
$\theta$ dans \m{\End(T)/\C}, \m{j(\ov{\theta})} est l'image de
\m{-\theta\circ\pi} dans \ \m{M_*\cap N=\Hom(F^{**},T)/\C}. La premi\`ere 
assertion de 3- en d\'ecoule imm\'ediatement.

L'assertion 4- d\'ecoule imm\'ediatement de 2-.
\end{proof}

\sepprop

Le sous-espace vectoriel \m{Q+M_*=Q_*+M=M+M_*} \ 
de $\Ext^1(\ke,\ke)$ correspond aux d\'e-\break
formations de $\ke$ qui proviennent
des d\'eformations de l'extension large \ $0\to G^*\to \ke\to F\to 0$ \
lorsque $T$ reste fixe, et $Q$ (resp. $Q_*$) correspond aux d\'eformations
de l'extension large lorsque seuls $G^{**}$ (resp. $F^{**}$), $\rho$ et
$\pi$ bougent.

\sepprop

\pagebreak[2]\begin{subsub}{\bf Proposition : }\label{struc_ext_3}
Soit \ \m{P\subset N} l'image r\'eciproque (par le morphisme canonique
\hfil\break
\ \m{N\to\Ext^1(G,G)}) de l'image de \m{\Hom(G,\wT)} 
dans \m{\Ext^1(G,G)}.
Soit \ \m{P_*\subset N_*} l'image r\'eciproque (par le morphisme canonique
\ \m{N_*\to\Ext^1(F,F)}) de l'image de \m{\Hom(F,T)} dans \m{\Ext^1(F,F)}.
On voit $P$ et $P_*$ comme des sous-espaces vectoriels de \m{\Ext^1(\ke,\ke)}.
Alors on a
\[ P \ = \ P_* \ = N\cap N_* ,\]
\[P\cap(M+M_*) \ = \ P\cap Q \ = \ P\cap Q_* \ = \ Q\cap Q_* \ =
\ M\cap N_* + M_*\cap N ,\]
et des isomorphismes canoniques
\[\Hom(F,T) \ \simeq \ (N\cap N_*)/(M\cap N_*), \ \ \ \
\Hom(G,\wT) \ \simeq \ (N\cap N_*)/(M_*\cap N) .\]
\end{subsub}

\begin{proof}
On montre d'abord que \ \m{P\subset N_*}. Il suffit de montrer que l'image de 
$P$ par le morphisme canonique \ \m{\Ext^1(\ke,\ke)\to\Ext^1(G^*,\ke)} \ est
nulle. 

L'image de \m{\Hom(G,\wT)} dans \m{\Ext^1(G,G)} est exactement le 
noyau
de l'application canonique \ \m{\Ext^1(G,G)\to\Ext^1(G,G^{**})}. L'assertion
d\'ecoule alors du diagramme commutatif
\xmat{
\Ext^1(G,\ke^*)\ar[r]\ar[d] & \Ext^1(\ke^*,\ke^*)\ar[r] & 
\Ext^1(\ke^*,G^{**})=\Ext^1(G^*,\ke) \\
\Ext^1(G,G)\ar[rr] & & \Ext^1(G,G^{**})\ar[u] }
d\'eduit de l'associativit\'e des Ext.

On a donc \ \m{P\subset N\cap N_*}. Montrons maintenant que \ 
\m{N\cap N_*\subset P}. Cela d\'ecoule du diagramme commutatif d\'eduit de
l'associativit\'e des Ext
\xmat{
\Ext^1(G,\ke^*)\ar[r]\ar[d] & \Ext^1(\ke^*,\ke^*)=\Ext^1(\ke,\ke)\ar[dd] \\
\Ext^1(G,G)\ar[d] & \\
\Ext^1(G,G^{**})\ar[r] & \Ext^1(G^*,\ke)=\Ext^1(\ke^*,G^{**}) }
et du fait que la fl\`eche horizontale du bas est injective, car \
\m{\Hom(F^*,G^{**})=\nsp}.

On a donc \ \m{P=N\cap N_*}. De m\^eme on a aussi \ \m{P_*=N\cap N_*}.
Une autre d\'emonstration d\'ecoulera de celle de la proposition
\ref{struc_ext_3b}.

La seconde assertion de la proposition \ref{struc_ext_3} d\'ecoule ais\'ement
de la premi\`ere et de la proposition \ref{struc_ext_2}, et la derni\`ere est
une cons\'equence de la premi\`ere et des d\'efinitions de $P$ et \m{P_*}.
\end{proof}

\sepprop

Le sous-espace vectoriel $P$ de $\Ext^1(\ke,\ke)$ correspond aux 
d\'eformations de l'extension large \ $0\to G^*\to \ke\to F\to 0$ \
lorsque $F^{**}$ et $G^{**}$ restent fixes.

\sepprop

\pagebreak[2]\begin{subsub}{\bf Proposition : }\label{kerx2}
On a \ $N+N_* \ = \ \ker(\xi_2) \ = \ \ker(\xi_2^*)$.
\end{subsub}

($\xi_2$ et $\xi_2^*$ sont d\'efinis en \ref{on_ext}).

\begin{proof}
On reprend les notations de \ref{on_ext}. Le noyau de \
\m{ A':\Ext^1(\ke,\ke)\lra\Ext^1(F^*,G)} \ 
est $N+M_*$. Donc \ \m{\xi_2^*(N+M_*)=\nsp}. D'autre part, on a d'apr\`es la
proposition \ref{id3}
\[\xi_2^*(N) \ = \ C'\circ g\circ f\circ d(N) \ = \ \nsp\]
car \ \m{C'\circ g=0}. On a donc \ \m{N+N_* \ \subset \ \ker(\xi_2^*)} .
L'inclusion inverse provient du diagramme de la proposition \ref{id3} et de la
surjectivit\'e de $d$ et $f$ (d'apr\`es le lemme \ref{reg_3x3}). On a donc
\ \m{\ker(\xi_2^*)=N+N_*}. On a de m\^eme \ \m{\ker(\xi_2)=N+N_*}.
\end{proof}

\sepprop

Rappelons qu'on a d'apr\`es le lemme \ref{struc_ext_1} des inclusions
canoniques
\[N/M \ \subset \ \Hom(\ke,T), \ \ \ \ N_*/M_* \ \subset
\ \Hom(\ke^*,\wT) .\]
Soient $P_0\subset N$ l'image r\'eciproque de \ 
\m{\Hom(F,T)\subset\Hom(\ke,T)}, et $P_0^*\subset N_*$ 
l'image r\'eciproque de \ \m{\Hom(G,\wT)\subset\Hom(\ke^*,\wT)}.

\sepprop

\pagebreak[2]\begin{subsub}{\bf Proposition : }\label{struc_ext_3b}
Soient
\[s_P : P\lra\Hom(G,\wT), \ \ \ \ s_{P_*}: P_*=P\lra\Hom(F,T) \]
les projections (d\'efinies par les surjections \m{N\to\Ext^1(G,G)},
\m{N_*\to\Ext^1(F,F)} respectivement), et
\[s_{P_0} : P_0\lra\Hom(F,T), \ \ \ \
s_{P_0^*} : P_0^*\lra\Hom(G,\wT) \]
les projections canoniques. Alors on a \ \m{P=P_0\cap P_0^*}, et les
diagrammes suivants
\[\xymatrix{
P\ar[rrd]^{s_P}\flinc[dd]^{P_0^*} \\
& & \Hom(G,\wT) \\
P_0^*\ar[rru]^{s_{P_0^*}}
} \ \ \ \
\xymatrix{
P\ar[rrd]^{s_{P_*}}\flinc[dd]^{P_0} \\
& & \Hom(F,T) \\
P_0\ar[rru]^{s_{P_0}}
}\]
sont commutatifs.
\end{subsub}

\begin{proof}
On utilise les r\'esultats et notations de \ref{Larg_constr}. Soit \
\m{\nu\in\Ext^1(F,\ke)=N_*}, repr\'esent\'e par un morphisme \
\m{\nu_1 : F_1\to\ke} \ s'annulant sur $\imm(f_2)$. L'image $\nu'$ de $\nu$
dans \m{\Ext^1(F,F)} est repr\'esent\'ee par \ \m{p\nu_1:F_1\to F}. 
Supposons que \m{\nu'\in\Hom(F,T)}. Ceci \'equivaut \`a dire que l'image de
$\nu'$ dans \m{\Ext^1(F,F^{**})} est nulle. Cette image est repr\'esent\'ee
par le compos\'e
\xmat{ \nu''_1 : F_1\ar[r]^{\nu'_1} & F\flinc[r] & F^{**}}
Il existe donc \ \m{\nu_0 : F_0\to F^{**}} \ tel que \
\m{\nu''_1 = \nu_0f_1}. On a alors \ \m{\nu_0(\imm(f_1))\subset F}, donc
$\nu_0$ induit un morphisme \ \m{\ov{\nu} : F\to T} \ tel qu'on ait un
diagramme commutatif
\xmat{
F_0\ar[r]^{\nu_0}\ar[d]^{f_0} & F^{**}\ar[d]^\pi\\
F\ar[r]^{\ov{\nu}} & T }
et on a \ \m{\ov{\nu}=\nu'}. 
On consid\`ere la r\'esolution de $\ke$
\xmat{
\cdots F_2\ar[r]^-{f_2} & F_1\ar[r]^-{(f_1,\sigma_1)} & 
F_0\oplus G^*\ar[r]^{\epsilon_0} &\ke\ar[r] & 0 }
et la r\'esolution de $\ke^*$ (cf. la d\'emonstration de la proposition
\ref{id4b})
\xmat{
\cdots F'_2\oplus G_2\ar[r]^-{\delta_2^*} & F'_1\oplus G_1\ar[r]^{\delta_1^*} 
& F'_0\oplus G_0\ar[r]^-{\delta_0^*} & \ke^*\ar[r] & 0
}
Soit \ \m{\eta_1=(\gamma_1,\psi_1) : F'_1\oplus G_1\to\ke^*} \ un morphisme
s'annulant sur $\imm(\delta_2^*)$ repr\'esentant $\nu$, vu comme
\'el\'ement de \m{\Ext^1(\ke^*,\ke^*)}. C'est un \'el\'ement de 
\m{\Ext^1(G,\ke^*)} d'apr\`es la proposition \ref{struc_ext_3}. On peut
donc supposer que \ \m{\gamma_1=0}. D'apr\`es \ref{Ext_dual} on peut
supposer qu'il existe un diagramme commutatif
\xmat{
F_1\ar[rr]^-{(f_1,\sigma_1)}\ar[d]^{\nu_1} & &
F_0\oplus G^*\ar[r]^-{\delta_0}\ar[d]^\Theta &\ke\ar[d]^{{}^t\eta_1}\\
\ke\ar[rr]^-{{}^t\delta_0^*} & & {F'_0}^*\oplus G_0^*\ar[r]^-{{}^t\delta_1^*}
& {F'_1}^*\oplus G_1^*
}
Soit \m{\left(\begin{array}{cc}a & b\\ c & d\end{array}\right)} la matrice
de $\Theta$ et \m{\left(\begin{array}{cc}0 & 0\\ v & w\end{array}\right)}
celle de \ \m{{}^t\eta_1\delta_0 : F_0\oplus G^*\to {F'_1}^*\oplus G_1^*}.
Du diagramme commutatif pr\'ec\'edent on d\'eduit \
\m{{}^tf_1'a = 0}, \m{{}^tf_1'b = 0} .
Donc $b$ se factorise
\xmat{b : G^*\ar[r] & F^{**}\ar[r]^{{}^tf'_0} & {F'_0}^*}
et comme \ \m{\Hom(G^*,F^{**})=\nsp}, on a \ \m{b=0}. De m\^eme $a$ se
factorise
\[\xymatrix{a : F_0\ar[r]^{a_0} & F^{**}\ar[r]^{{}^tf'_0} & {F'_0}^*} .\]
On en d\'eduit le diagramme commutatif
\[(E) \ \ \ \ \xymatrix{
F_1\ar[rr]^-{(f_1,\sigma_1)}\ar[d]^{\nu_1} & &
F_0\oplus G^*\ar[rr]^-{\delta_0}\ar[d]^{\Theta'} & & \ke\ar[d]^{{}^t\eta_1}\\
\ke\flinc[rr] & & F^{**}\oplus G_0^*\ar[rr]^-{({}^t\sigma_1^*,{}^tg_1)} & &
G_1^*}\]
o\`u $\Theta'$ a pour matrice 
\m{\left(\begin{array}{cc}a_0 & 0\\ c & d\end{array}\right)} et \ 
\m{{}^t\eta_1\delta_0=(v,w)}. De $(E)$ on d\'eduit que \ 
\m{w={}^tg_1d}. Il en d\'ecoule que l'image de $\nu$ (vu comme
\'el\'ement de \m{\Ext^1(G,\ke^*)}) dans \m{\Ext^1(G,G)} appartient \`a
\m{\Hom(G,\wT)}, et que le morphisme \ \m{\ov{\eta}:G\to\wT} \
correspondant est induit par \m{{}^td}, c'est-\`a-dire qu'on a \
\m{{}^td(\imm(g_1))\subset G}, et un diagramme commutatif
\xmat{ G^0\ar[r]^{{}^td}\flon[d] & G^{**}\flon[d]^\rho \\
G\ar[r]^{\ov{\eta}} & \wT }
On retrouve ainsi le r\'esultat de la proposition \ref{struc_ext_3}.

De $(E)$ on d\'eduit aussi que \ \m{\nu_0f_1=a_0f_1}, donc on peut
\'ecrire \ \m{a_0 = \nu_0 + tjf_0},
o\`u $t\in\C$ et $j$ est l'inclusion $F\subset F^{**}$. 
On consid\`ere maintenant l'image $\nu''$ de $\nu$ par le morphisme
\[\Ext^1(G,\ke^*)\lra\Ext^2(\wT,\ke^*)=\Hom(\ke,T)\]
induit par $\rho$. Elle provient de ${}^t\eta_1$ compte tenu de la
r\'esolution de $T$
\xmat{
G^{*}\ar[r]^{{}^tg_0} & G_0^*\ar[r]^-{{}^tg_1} & \ker({}^tg_2)\ar[r]^-{\pi'}
& T\ar[r] & 0 }
(c'est-\`a-dire que \ \m{\nu'=\pi'{}^t\eta_1}). L'image de $\nu''$ dans
\m{\Hom(G^*,T)} provient donc de \ \m{{}^tg_1d : G^*\to G_1^*}, et est donc
nulle. Donc \m{\nu''\in\Hom(F,T)}. Cela montre d\'ej\`a que $P\subset P_0$.
On a un diagramme commutatif
\xmat{
G_0^*\ar[r]^-{{}^tg_1} & \ker({}^tg_2)\ar[r]^-{\pi'} & T\\
G^*\ar[u]^d\flinc[r] & F_0\oplus G^*\ar[u]^{{}^t\eta_1\delta_0}\flon[r]
& F_0\ar[u]^{\nu''_0}
}
et $\nu''$ est induit par $\nu''_0$, et provient donc de \ 
\m{v : F_0\to G_1^*}. On a d'apr\`es $(E)$
\[v \ = \ {}^t\sigma_1^*a_0 + {}^tg_1c . \]
Donc $\nu''$ provient de \ \m{{}^t\sigma_1^*a_1 : F_0\to G_1^*}, et donc 
aussi de \ \m{a_0 : F_0\to F^{**}}, compte tenu de l'autre r\'esolution de 
$T$
\xmat{ \cdots F_1\ar[r]^-{f_1} & F_1\ar[r]^{f_0} & F^{**}\ar[r]^\pi
& T & 0 }
et de l'isomorphisme canonique entre les deux r\'esolutions de $T$ donn\'e 
\`a la fin de \ref{Larg_constr}. 

Comme \ \m{a_0=\nu_0+tjf_0}, $\nu''$ provient aussi de $\nu_0$. On a donc
\ $\nu''=\ov{\nu}$. Ceci prouve que le second diagramme de la proposition
\ref{struc_ext_3b} est commutatif.
\end{proof}

\sepprop

\pagebreak[2]\begin{subsub}{\bf Corollaire : }\label{coro_struc_ext_3b}
On a \m{\Hom(F,T)\subset N/M} et \m{\Hom(G,\wT)\subset N_*/M_*}, compte tenu
des inclusions \m{N/M\subset\Hom(\ke,T)}, \m{N_*/M_*\subset\Hom(\ke^*,\wT)}.
\end{subsub}

\begin{proof} Cela d\'ecoule du fait que \m{s_P} et \m{s_{P_*}} sont
surjectives. \end{proof}
\end{sub}

\sepsub

\Ssect{Le tangent \`a l'espace des extensions}{s_ext_6}

On pose
\[{\bf T} \ = \ {\bf T}(\pi,\rho) \ = \ P+M+M_* 
\ = \ (M+N_*)\cap(M_*+N) .\]
On note $H$ (resp. $H_*$) le noyau de l'application canonique
\[\Hom(G^*,T)\lra\Ext^1(F,T) \ \ \ \ \
{\rm (resp. \ }\Hom(F^*,\wT)\lra\Ext^1(G,\wT) \ {\rm )} \]
(cf. prop. \ref{id4}).

\sepprop

\pagebreak[2]\begin{subsub}{\bf Proposition : }\label{tan2}
Le quotient \ \m{(N+N_*)/{\bf T}} \ s'identifie naturellement \`a la somme
directe d'un sous-espace vectoriel de $H$ et d'un sous-espace vectoriel de
\m{H_*}.
\end{subsub}

\begin{proof}
Rappelons que \m{N/M} (resp. \m{N_*/M_*}) s'identifie naturellement \`a un
sous-espace vectoriel de \m{\Hom(\ke,T))} (resp. \m{\Hom(\ke^*,\wT)}).
On a un morphisme bien d\'efini
\[\Phi : N/(M+N\cap N_*)\oplus N_*/(M_*+N\cap N_*)\lra(N+N_*)/{\bf T}\]
tel que \ \m{\Phi([n],[n_*]) = [n+n_*]} \
pour tous $n$, $n'$ dans $N$, $N_*$ respectivement. Il est clair que c'est un
isomorphisme, compte tenu du fait que \ \m{P=N\cap N_*} \ d'apr\`es la
proposition \ref{struc_ext_3}. Il suffit donc de trouver des isomorphismes 
canoniques
\[N/(M+N\cap N_*)\simeq H, \ \ \ \
N_*/(M_*+N\cap N_*)\simeq H_* . \]
On ne d\'efinira que le premier, le second \'etant analogue. C'est une
cons\'equence de la suite exacte
\xmat{
0\ar[r] & \Hom(F,T)\ar[r] & \Hom(\ke,T)\ar[r] & \Hom(G^*,T)\ar[r] & 
\Ext^1(F,T)\ar[r] & 0
}
de la derni\`ere assertion de la proposition \ref{struc_ext_3} et du corollaire
\ref{coro_struc_ext_3b}.
\end{proof}

\sepprop

\pagebreak[2]\begin{subsub}{\bf Proposition : }\label{tan3}
Si \ \m{\dim(X)\geq 3} et \ \m{\Hom(G^*,T)=\Hom(F^*,\wT)=\nsp}, on a \hfil\break
\m{\Ext^1(\ke,\ke)={\bf T}}.
\end{subsub}

\begin{proof}
On a \ \m{\Ext^1(G^*,F)\simeq\Hom(G^*,T)=\nsp} \ et \
\m{\Ext^1(F^*,G)\simeq\Hom(F^*,\wT)=\nsp}, donc \
\m{\Ext^1(\ke,\ke)=N+M_*=N_*+M=N+N_*} \
d'apr\`es \ref{s_ext_0}. D'autre part on a aussi
 \ \m{N+N_*={\bf T}} \ d'apr\`es la proposition \ref{tan2}, ou la d\'efinition
 de $\bf T$ plus haut.
\end{proof}

\sepprop

Soient $U$, $V$, $Z$ des vari\'et\'es alg\'ebriques irr\'eductibles
r\'eduites, $\F$,
$\G$ des familles de fibr\'es vectoriels 2-lisses
sur $X$ param\'etr\'ees par $U$, $V$ respectivement. 
Soit $\kt$ une famille de faisceaux parfaits de codimension 2 sur $X$ 
param\'etr\'ee par $Z$, plate sur $Z$. 
On suppose que pour tous point ferm\'es $u$, $v$, $z$ de $U$, $V$, $Z$ 
respectivement on a
\[\Ext^i(\F_u,\kt_z)=\Ext^i(\G_v,\widetilde{\kt}_z)=\nsp \ \ \ \
{\rm si} \ i\geq 1 .\]
On suppose aussi que les d\'eformations semi-universelles des faisceaux
$\kt_z$ et $\widetilde{\kt}_z$ sont r\'eduites. Alors les faisceaux 
\[\kh \ = \ p_{U\times Z *}(\HHom(p^\sharp_U(\F),p^\sharp_Z(\kt))), \ \ 
\kk \ = \ p_{U\times Z *}(\HHom(p^\sharp_V(\G),p^\sharp_Z(\widetilde{\kt})))\]
sont localement libres (cf. \ref{modfin} pour les notations).
Soient $\kh_0$ l'ouvert de $\kh$ vu comme vari\'et\'e alg\'ebrique 
correspondant aux morphismes surjectifs dont le noyau est un faisceau 
r\'egulier, et $\kk_0$ l'ouvert analogue de $\kk$. Soient 
\[\pi_\F : \kh_0\lra U\times Z, \ \ \ \
\pi_\G : \kk_0\lra V\times Z \]
les projections. On a des
morphismes canoniques universels surjectifs de faisceaux coh\'erents sur
$\kh_0\times X$ et $\kk_0\times X$ respectivement :
\[\Pi : (p_U\circ\pi_\F)^\sharp(\F)\lra(p_T\circ\pi_\F)^\sharp(\kt),
\ \ \ \
R : (p_V\circ\pi_\G)^\sharp(\G)\lra(p_T\circ\pi_\G)^\sharp(\widetilde{\kt}) .
\]
Soient \ \m{\kf=\ker(\Pi)}, \m{\kg=\ker(R)}. Ce sont des familles plates
de faisceaux r\'eguliers sur $X$. On a une extension universelle sur
$\kh_0\times\kk_0\times X$
\[0\lra (p_V\circ\pi_\G\circ p_{\kk_0})^\sharp(\G^*)\lra\E\lra
(p_U\circ\pi_\F\circ p_{\kh_0})^\sharp(\kf)\lra 0  .\]
Soit $W$ l'ouvert de $\kh_0\times\kk_0$ correspondant aux extensions larges.
On suppose qu'il est non vide. Soit
\[\D \ = \ R^1p_{W*}(\E^*\ot\E) .\]
C'est un faisceau localement libre sur $W$. Soit $\T$ le sous-fibr\'e 
vectoriel de $\D$ d\'efini de la fa\c con suivante : soit $w$ un point
ferm\'e de $W$, qu'on peut voir comme une paire de morphismes
\[\pi : \F_u\lra\kt_z, \ \ \ \ \rho : \G_v\lra\widetilde{\kt}_z\]
o\`u $u$, $v$, $z$ sont les projections de $w$ sur $U$, $V$, $Z$
respectivement. Alors on a
\[ \T_w \ = \ {\bf T}(\pi,\rho) .\]

Le r\'esultat suivant d\'ecoule de \ref{Def_reg}, \ref{Def_Extens} et
\ref{s_ext4} :

\sepprop

\pagebreak[2]\begin{subsub}{\bf Proposition : }\label{tan_def}
Soit $w\in W$. Alors le morphisme de d\'eformation infinit\'esimale de
Koda\" ira-Spencer de $\E$ au point $w$
\[ \omega_w : T_wW\lra\Ext^1(\E_w,\E_w) \]
est \`a valeurs dans $\T_w$. Soient $u$, $v$, $z$ sont les projections de 
$w$ sur $U$, $V$, $Z$ respectivement. Si $\F$ est une d\'eformation 
compl\`ete de $\F_u$, $\G$ une d\'eformation compl\`ete de $\G_v$ et
$\kt$ une d\'eformation compl\`ete de $\kt_z$, alors l'image de
$\omega_w$ est exactement $\T_w$.
\end{subsub}

\sepprop

On note $W(\pi,\rho)$ le noyau de l'application 
\[A_2(\sigma)\oplus A_2(\sigma^*)\lra\Ext^2(T,T)\]
\'egale \`a la restriction de $(\theta,\theta^*)$, o\`u $\theta$ et
$\theta^*$ sont respectivement les applications canoniques
\[ \Ext^1(G^*,F)\lra\Ext^2(\wT,\wT), \ \ \ \
\Ext^1(F^*,G)\lra\Ext^2(T,T)\]
d\'efinies en \ref{on_ext}. Soit
\[\Delta : \Ext^1(\ke,\ke)\lra A_2(\sigma)\oplus A_2(\sigma^*)\]
l'application canonique.
Des propositions \ref{id4b} et \ref{tan2} 
on d\'eduit imm\'ediatement la

\sepprop

\pagebreak[2]\begin{subsub}{\bf Proposition : }\label{tan4}
On a \ \m{{\bf T}=\ker(\Delta)}, \ \m{\Delta(\Ext^1(\ke,\ke))=W(\pi,\rho)}
\hfil\break et \
\m{\Delta(N+N_*)=\ker(\theta)\oplus\ker(\theta^*)}.
\end{subsub}

\end{sub}

\sepsub

\Ssect{Morphismes \`a valeurs dans $\Ext^1(T,T)$}{s_ext_5}

On a un morphisme canonique
\xmat{
\nu : P\flon[r] & \Hom(G,\wT)\flinc[r] & \Ext^1(G,G)
}
qui est la restriction \`a $P$ du morphisme \ $N\to\Ext^1(G,G)$ \ de
\ref{s_ext_0}. Soit
\[\phi_\wT : P\lra\Ext^1(\wT,\wT)\]
le morphisme compos\'e
\xmat{
P\ar[r]^-\nu & \Hom(G,\wT)\flinc[r] & \Ext^1(\wT,\wT)
}
(pour le second morphisme voir le lemme \ref{reg_3x3}). On a de m\^eme un
morphisme canonique
\[\phi_T : P\lra\Ext^1(T,T) .\]

\sepprop

\pagebreak[2]\begin{subsub}{\bf Proposition : }\label{id5}
On a \ \m{\phi_T = \phi_\wT}, compte tenu de l'identification \
\m{\Ext^1(T,T)\simeq\Ext^1(\wT,\wT)}.
\end{subsub}

\begin{proof}
Analogue \`a celles des propositions \ref{id1} et \ref{id4b}.
\end{proof}
\end{sub}

\sepsub

\Ssect{L'action de $\Hom(\ke,\ke)$ sur $\Ext^1(\ke,\ke)$}{s_ext_4}

D'apr\`es la proposition \ref{larg_lemm2}, on a
\[\Hom(\ke,\ke) \ \simeq \ \C I_\ke\oplus\Hom(F,G^*).\]
On peut donc se restreindre \`a \'etudier l'action de $\Hom(F,G^*)$ sur
$\Ext^1(\ke,\ke)$. On a des isomorphismes canoniques
\[\Hom(F,G^*) \ \simeq \ \Hom(F^{**},G^*) \ \simeq \
\Hom(G^{**},F^*) \ \simeq \ \Hom(G,F^*)  .\]
On note
\[\mu_G : \Hom(\ke,\ke)\ot\Ext^1(\ke,\ke)\lra\Ext^1(\ke,\ke), \ \
\mu_D : \Ext^1(\ke,\ke)\ot\Hom(\ke,\ke)\lra\Ext^1(\ke,\ke)\]
les multiplications.

\sepprop

\pagebreak[2]\begin{subsub}{\bf Proposition : }\label{act_hom}
Les restrictions de $\mu_G$ \`a \ $\Hom(F,G^*)\ot (M+N_*)$ \ et de $\mu_D$
\`a \hfil\break $(M_*+N)\ot\Hom(F,G^*)$ sont nulles. Compte tenu des inclusions
\[\Ext^1(\ke,\ke)/(M+N_*)\simeq 
A_2(\sigma)\subset\Ext^1(G^*,F)\simeq\Hom(G^*,T) , \]
\[\Hom(F^{**},T)/\C\pi = \Ext^1(G,F^*)/\C\sigma^*\subset\Ext^1(\ke,\ke) ,\]
\[\Ext^1(\ke,\ke)/(M_*+N)\simeq 
A_2(\sigma^*)\subset\Ext^1(F^*,G)\simeq\Ext^2(T,F^{**}) , \]
\[\Ext^2(T,G^*)/\C\rho = \Ext^1(F,G^*)/\C\sigma^*\subset\Ext^1(\ke,\ke) ,\]
$\mu_G$ est induite par l'application canonique
\[\Hom(F^{**},G^*)\ot\Hom(G^*,T)\lra\Hom(F^{**},T)/\C\pi ,\]
et $\mu_D$ par
\[\Ext^2(T,F^{**})\ot\Hom(F^{**},G^*)\lra\Ext^2(T,G^*)/\C\rho .\]
\end{subsub}

\begin{proof}
On va montrer que la restriction de $\mu_G$ \`a \ $\Hom(F,G^*)\ot(M+N_*)$ \ 
est nulle (l'assertion concernant $\mu_D$ est analogue).

Rappelons que \ \m{N_*=\Ext^1(F,\ke)}. On a un diagramme commutatif
\xmat{
\Hom(F,G^*)\ot\Ext^1(F,\ke)\flinc[d]\\
\Hom(F,G^*)\ot\Ext^1(\ke,\ke)\flinc[r]\ar[d] & 
\Hom(\ke,\ke)\ot\Ext^1(\ke,\ke)\ar[d]\\
\Hom(F,G^*)\ot\Ext^1(G^*,\ke)\ar[r]^\phi & \Ext^1(\ke,\ke)
}
o\`u la colonne de gauche est exacte. Ceci montre que la restriction de $\mu_G$
\`a \ \m{\Hom(F,G^*)\ot N_*} \ est nulle. Il reste donc \`a montrer que la
restriction $\phi'$ de l'application pr\'ec\'edente $\phi$ \`a 
\m{\Hom(F,G^*)\ot\Ext^1(G^*,G^*)} est nulle. Mais on a 
\m{\Hom(F,G^*)=\Hom(F^{**},G^*)}, donc $\phi'$ se factorise par
\m{\Ext^1(F^{**},G^*)}, qui est nul. Donc \m{\phi'=0}. Les autres assertions se
d\'emontrent ais\'ement en utilisant les r\'esultats de \ref{s_ext_0}.
\end{proof}

\sepprop

On en d\'eduit que \m{\Aut(\ke)} agit trivialement sur $\bf T$, mais n'agit pas
trivialement sur \m{\Ext^1(\ke,\ke)} si celui-ci est distinct de $\bf T$, en
particulier lorsque $X$ est une surface.

\end{sub}

\sepsec


\pagebreak[4]\section{Vari\'et\'es de modules d'extensions larges}\label{vmod}

\medskip

\Ssect{Construction des vari\'et\'es de modules}{vmod0}

\pagebreak[2]\begin{subsub}\label{vmod_hyp} Hypoth\`eses. \rm
Soient $\kx$, $\ky$, $\kz$ des ensembles ouverts de faisceaux coh\'erents sur 
$X$, admettant des vari\'et\'es de modules fins {\bf M}, {\bf N}, {\bf Z} 
respectivement (cf.\ref{modfin}), les faisceaux de $\kz$ \'etant parfaits de
codimension 2. Si $X$ est une surface, on suppose qu'il
existe un entier positif $k$ tel que {\bf Z} soit l'ouvert de \m{\Hilb^k(X)}
des sous-sch\'emas constitu\'es de $k$ points distincts. On note $\F$, $\G$,
$\T$ les faisceaux universels sur \m{{\bf M}\times X}, \m{{\bf N}\times X} et
\m{{\bf Z}\times X} respectivement.
On note \m{\widetilde{\T}} la famille de faisceaux d\'eduite de $\T$ v\'erifiant
\m{\widetilde{\T}_z=\widetilde{(\T_z)}} pour tout \m{z\in {\bf Z}} (on peut 
construire
\m{\widetilde{\T}} par exemple en utilisant des r\'esolutions localement libres
locales de $\T$). Si $X$ est une surface on a \ 
\m{\widetilde{\T}\simeq\T\ot p_X^*(\omega_X^{-1})}
(\m{p_X} d\'esignant la projection \m{{\bf Z}\times X\to X}).

On suppose que tous les faisceaux de $\kx$, $\ky$ sont localement libres,
simples et
2-lisses, et que si $A$, $B$ sont deux fibr\'es de $\kx$ (resp. $\ky$) non
isomorphes, alors on a \ \m{\Hom(A,B)=\nsp}.

On suppose aussi que si \m{\dim(X)>2} les faisceaux de $\kz$ sont simples, et
que pour tout \m{(m,n,z)\in{\bf M}\times{\bf N}\times{\bf Z}}, les
propri\'et\'es suivantes sont v\'erifi\'ees :
\begin{enumerate}
\item[(i)] On a, si $i\geq 1$,
$$\Ext^i(\F_m,\T_z) \ = \ \Ext^i(\G_n,\widetilde{\T}_z) \ = \ \nsp$$
\item[(ii)] Pour tout noyau $E$ d'un morphisme surjectif $\F_m\to\T_z$ (resp.
$\G_n\to\widetilde{\T}_z$), on a \ $D(E)=\Ext^1(E,E)$.
\end{enumerate}
La condition (ii) est v\'erifi\'ee si \m{H^0(\F_m\ot\widetilde{\T}_z)=\nsp}
(resp. \m{H^0(\G_n\ot\T_z)=\nsp}) d'apr\`es la proposition \ref{reg_3x3b}.

Soient $p_M$, $p_N$, $p_Z$, $p$ les projections de \ \m{{\bf M}\times{\bf
N}\times{\bf Z}} \ sur {\bf M}, {\bf N}, \m{\bf Z}, et de \m{{\bf M}\times{\bf
N}\times{\bf Z}\times X} sur \ \m{{\bf M}\times{\bf N}\times{\bf Z}} \
respectivement. Soient
$$\kf \ = \ p_*\big(\HHom(p_M^\sharp(\F),p_Z^\sharp(\T))\big) , \ \ \ \ 
\kg \ = \ p_*\big(\HHom(p_N^\sharp(\G),p_Z^\sharp(\widetilde{\T}))\big) ,
\ \ \ \ \Gamma \ = \ \kf\oplus\kg$$
qui sont des faisceaux localement libres, c'est-\`a-dire des fibr\'es vectoriels
sur \ \m{{\bf M}\times{\bf N}\times{\bf Z}}.
Si \ \m{(m,n,z)\in{\bf M}\times{\bf N}\times{\bf Z}} \ on a
$$\kf_{(m,n,z)} \ \simeq \
\Hom(\F_m,\T_z) , \ \ \ \ \kg_{(m,n,z)} \ \simeq \ 
\Hom(\G_n,\widetilde{\T}_z).$$
Soient \m{\kf^{surj}}, \m{\kg^{surj}} les ouverts correspondant morphismes
surjectifs, et \hfil\break \m{\Gamma^{surj}=
\kf^{surj}\times_{{\bf M}\times{\bf N}\times{\bf Z}}\kg^{surj}\subset\Gamma}.
\end{subsub}

\sepprop

\pagebreak[2]\begin{subsub}{\bf Lemme : }\label{vmod1}
Soient $A$, $A'$ des faisceaux de $\kx$ (resp. $\ky$), $T$, $T'$ des faisceaux
de $\bf Z$ et \ \m{\pi:A\to T}, \m{\pi':A'\to T'} \ des morphismes surjectifs.
Alors si \ \m{\ker(\pi)\simeq\ker(\pi')}, on a \m{A=A'}, \m{T=T'}, et il existe
un automorphisme $g$ de $T$ tel que \m{g\circ\pi=\pi'}.
\end{subsub}

\begin{proof}
Puisque \m{\ker(\pi)\simeq\ker(\pi')}, on a \
\m{A\simeq\ker(\pi)^{**}\simeq\ker(\pi')^{**}\simeq A'} ,
donc \m{A=A'}. Puisque $A$ est simple l'isomorphisme induit 
\m{\ker(\pi)^{**}\simeq\ker(\pi')^{**}} est une homoth\'etie et
\m{\ker(\pi)=\ker(\pi')} comme sous-faisceaux de $A$. On en d\'eduit 
$$T\simeq A/\ker(\pi)=A'/\ker(\pi')\simeq T'$$
et le lemme en d\'ecoule imm\'ediatement.
\end{proof}

\sepprop

Soit \m{\Gamma^0\subset\Gamma^{surj}} l'ouvert correspondant aux extensions
larges. Au dessus de \hfil\break
\m{(m,n,z)\in{\bf M}\times{\bf N}\times{\bf Z}}, \m{\Gamma^0}
est l'ensemble des \ 
\m{(\pi,\rho)\in\Hom(\F_m,\T_z)\oplus\Hom(\G_n,\widetilde{\T}_z)} \ tels que
$\pi$ et $\rho$ soient surjectifs et
\begin{enumerate}
\item[-] $h^i(\F_m^*\ot\G_n^*)=0$ \ si $i\geq 1$,
\item[-] $\Ext^2(\ker(\pi),\G_n^*)=\Ext^2(\ker(\rho),\F_m^*)=\nsp$,
\item[-] $h^i(\F_m\ot\G_n)=0$ \ si $i\leq 1$,
\end{enumerate}

Si \m{(\pi,\rho)\in\Gamma^0_{m,n,z}}, on note \m{\ke(\pi,\rho)} l'extension
large correspondante. On a donc des suites exactes
$$0\lra\G_n^*\lra\ke(\pi,\rho)\lra\ker(\pi)\lra 0 ,$$
$$0\lra\F_m^*\lra\ke(\pi,\rho)^*\lra\ker(\rho)\lra 0 .$$

\sepprop

\pagebreak[2]\begin{subsub}{\bf Proposition : }\label{vmod2}
Soient \m{(m,n,z), (m',n',z')\in{\bf M}\times{\bf N}\times{\bf Z}},
\m{(\pi,\rho)\in\Gamma^0_{(m,n,z)}}, \hfil\break
\m{(\pi',\rho')\in\Gamma^0_{(m',n',z')}}.
Alors on a \ \m{\ke(\pi,\rho)\simeq\ke(\pi',\rho')} \ si et seulement si
\m{(m,n,z)=(m',n',z')} et s'il existe \m{\lambda\in\C^*}, \m{g\in\Aut(\T_z)},  
tels que \ \m{\pi'=g\circ\pi}, \m{\lambda\widetilde{g}\circ\rho'=\rho}.
\end{subsub}

\begin{proof}
Soit \ \m{\theta:\ke(\pi,\rho)\to\ke(\pi',\rho')} \ un isomorphisme. On
consid\`ere les suites exactes
$$0\lra\G_n^*\lra\ke(\pi,\rho)\lra\ker(\pi)\lra 0 ,$$
$$0\lra\G_{n'}^*\lra\ke(\pi',\rho')\lra\ker(\pi')\lra 0 .$$
On a \ \m{\Hom(\G_n^*,\ker(\pi'))=\nsp}, donc \ 
\m{\theta(\G_n^*)\subset\G_{n'}^*}. Comme \ \m{\Hom(\G_n^*,\G_{n'}^*)=\nsp} \ si
\m{n\not=n'}, on a \m{n=n'}, et la restriction de $\theta$ \`a \m{\G_n^*} est
une homoth\'etie de rapport \m{\gamma\not=0}. Donc $\theta$ induit un
isomorphisme \m{\ker(\pi)\simeq\ker(\pi')}. D'apr\`es le lemme \ref{vmod1} on a
\m{m=m'}, \m{z=z'}, et il existe \m{g\in\Aut(\T_z)} tel que \m{\pi'=g\circ\pi}.
On a aussi \m{\ker(\pi)=\ker(\pi')} (comme sous-faisceaux de \m{\F_m}) et
l'automorphisme de \m{\F_m} induit par $\frac{1}{\gamma}\theta$ est une 
homoth\'etie de rapport
\m{\lambda\not=0}. On a donc un diagramme commutatif
\xmat{
0\ar[r] & \G_n^*\ar[r]\fleq[d] & 
\ke(\pi,\rho)\ar[r]\ar[d]^{\frac{1}{\gamma}\theta} &
\ker(\pi)\ar[d]^\lambda\ar[r] & 0\\
0\ar[r] & \G_n^*\ar[r] & \ke(\pi',\rho')\ar[r] & \ker(\pi')\ar[r] & 0
}
En dualisant on en d\'eduit le carr\'e commutatif
\xmat{
\G_n\fleq[d]\ar[r]^-\alpha & \EExt^1(\ker(\pi),\ko_X)\\
\G_n\ar[r]^-{\alpha'} & \EExt^1(\ker(\pi'),\ko_X)\ar[u]^\lambda
}
On a un diagramme commutatif
\xmat{
0\ar[r] & \ker(\pi)\ar[r]\fleq[d] & \F_n\ar[r]\fleq[d] & \T_z\ar[d]^g\ar[r] &
0\\
0\ar[r] & \ker(\pi')\ar[r] & \F_n\ar[r] & \T_z\ar[r] & 0
}
En dualisant on en d\'eduit le carr\'e commutatif
\xmat{
\EExt^1(\ker(\pi),\ko_X)\ar[r]^-\partial\fleq[d] &
\EExt^2(\T_z,\ko_X)=\widetilde{\T}_z \\
\EExt^1(\ker(\pi'),\ko_X)\ar[r]^-{\partial'} &
\EExt^2(\T_z,\ko_X)=\widetilde{\T}_z\ar[u]^{\widetilde{g}}
}
En regroupant ces deux carr\'es commutatifs on obtient le carr\'e commutatif
\xmat{
\G_n\fleq[d]\ar[rr]^-\alpha\ar@/^2pc/[rrrr]^\rho
 &  &
\EExt^1(\ker(\pi),\ko_X)\ar[rr]^-\partial
&  & {\widetilde{\T}_z}
 \\
\G_n\ar[rr]^-{\alpha'}\ar@/_2pc/[rrrr]_{\rho'}
& & \EExt^1(\ker(\pi'),\ko_X)\ar[rr]^-{\partial'}\ar[u]^\lambda & &
{\widetilde{\T}_z}\ar[u]^{\lambda\widetilde{g}}
}
On a donc \ \m{\lambda\widetilde{g}\circ\rho'=\rho}. La r\'eciproque est
imm\'ediate.
\end{proof}

\sepprop

\pagebreak[2]\begin{subsub}\label{vmod3b2} Construction des vari\'et\'es de 
modules. \rm On va construire une "vari\'et\'e de modules" pour les extensions
pr\'ec\'edentes. Traitons d'abord le cas le plus simple, c'est \`a dire
\m{\dim(X)>2}. Soient 
\[U \ = \P(\kf^{surj})\times_{{\bf M}\times{\bf N}\times{\bf Z}}\P(\kg^{surj})\]
et \m{\M(\kx,\ky,\kz)} l'ouvert de $U$ correspondant aux extensions larges.

On suppose maintenant que $X$ est une surface.
Soit ${U_k}$ l'ouvert de $X^k$ constitu\'e des \m{(x_1,\ldots,x_k)}
tels que \m{x_i\not=x_j} si \m{1\leq i<j\leq k}. Pour \m{1\leq j\leq k} soient
\m{\lambda_j:{U_k}\to X} la restriction de la $j$-i\`eme projection,
$$f_i:{\bf M}\times{\bf N}\times{U_k}\lra{\bf M}\times X$$
la compos\'ee de la projection \m{{\bf M}\times{\bf N}\times X\to {\bf M}\times
X} et de \m{I_{{\bf M}\times{\bf N}}\times\lambda_i},
$$g_i:{\bf M}\times{\bf N}\times{U_k}\lra{\bf N}\times X$$
la compos\'ee de la projection \m{{\bf M}\times{\bf N}\times X\to {\bf N}\times
X} et de \m{I_{{\bf M}\times{\bf N}}\times\lambda_i}. Soient \
\m{\F_i=f_i^*(\F)}, \m{\G_i=g_i^\sharp(\G)} \ et
$$W \ = \ M^s(\F_1^*,\ldots,\F_k^*,\G_1^*,\ldots,\G_k^*) $$
(cf. \ref{quotalg}). Sur $W$ agit de mani\`ere \'evidente le groupe \m{\Sigma_k}
des permutations de \m{\lbrace 1,\ldots,k\rbrace}. Soit
$$\M(\kx,\ky,\kz) \ = \ W/\Sigma_k .$$
C'est une vari\'et\'e quasiprojective lisse. La projection \m{W\to
{\bf M}\times{\bf N}\times{U_k}} passe au quotient et d\'efinit un
morphisme \m{\M(\kx,\ky,\kz)\to{\bf M}\times{\bf N}\times{\bf Z}}. 

Dans tous les cas, d'apr\`es la
proposition \ref{vmod2} et la proposition \ref{quotalg1} les points ferm\'es de
\m{\M(\kx,\ky,\kz)} s'identifient aux classes d'isomorphisme d'extensions 
larges du type \m{\ke(\pi,\rho)}.
On note ${\rm Larg}(\kx,\ky,\kz)$ l'ensemble des classes d'isomorphisme
d'extensions larges \m{\ke(\pi,\rho)}, qui est donc aussi l'ensemble des points
ferm\'es de \m{\M(\kx,\ky,\kz)}.

Soit $r$ (resp. $s$) le rang des fibr\'es de $\kx$ (resp. $\ky$),
\m{d=\chi(\F_m,\T_z)}, \m{e=\chi(\G_n,\wT_z)} (pour \ \m{(m,n,z)\in
{\bf M}\times{\bf N}\times{\bf Z}}). Si $X$ est une surface on a
\[\dim(\M(\kx,\ky,\kz)) \ = \ \dim({\bf M})+\dim({\bf N})+k(r+s+1) - 1 ,\]
et si \m{\dim(X)>2}
\[\dim(\M(\kx,\ky,\kz)) \ = \ \dim({\bf M})+\dim({\bf N})+\dim({\bf Z})+d+e-2
.\]
\end{subsub}
\end{sub}

\sepsub

\Ssect{Familles pures d'extensions larges}{vmod01}

Soit \m{n=\dim(X)}.
Soit $\ke$ une famille de fibr\'es de \m{{\rm Larg}(\kx,\ky,\kz)} param\'etr\'ee
par une vari\'et\'e alg\'ebrique $S$. C'est donc un fibr\'e vectoriel sur
\m{S\times X}. Soit
$$\ku \ = \ p^\sharp_S(\ke^*)\ot p^\sharp_{\bf N}(\G^*)\ot p_X^*(\omega_X)$$
($p_S$, $p_{\bf N}$, $p_X$ d\'esignant les projections de \m{S\times{\bf
N}} sur $S$, {\bf N}, et de \m{S\times{\bf N}\times X} sur $X$ respectivement).
Pour tout point ferm\'e \m{(s,n)} de \m{S\times{\bf N}}, le morphisme
$$R^np_{S\times {\bf N}*}(\ku)\ot_{\ko_{(s,n)}}\C\lra
H^n(X,\ke_s^*\ot\G_n^*\ot\omega_X)$$
(\m{p_{S\times {\bf N}}} d\'esignant la projection \m{S\times{\bf N}\times X
\to S\times {\bf N}})
est un isomorphisme. Par dualit\'e de Serre on a un isomorphisme
$$H^n(X,\ke_s^*\ot\G_n^*\ot\omega_X) \ \simeq \ \Hom(\G_n^*,\ke_s)^* .$$
Par cons\'equent pour tout point ferm\'e $s$ de $S$ il existe un unique point
ferm\'e $n$ de $\bf N$ tel que \
\m{(s,n)\in{\rm supp}(R^np_{S\times {\bf N} *}(\ku))} ,
et on a \ \m{\dim(\Hom(G^*_{n},\ke_s)) = 1}.

\sepprop

\pagebreak[2]\begin{subsub}{\bf D\'efinition : }\label{vmod3}
On dit que $\ke$ est {\em pure} si il existe un morphisme \hfil\break
 \m{\phi:S\to
{\rm supp}(R^np_{S\times {\bf N} *}(\ku))} \ tel que \ \m{p_S\circ\phi=I_S},
et un fibr\'e en droites $L$ sur $S$ tel que \
\m{R^np_{S\times {\bf N} *}(\ku)\simeq\phi^*(L)}.
\end{subsub}

\sepprop

Dans le cas des fibr\'es instables de rang 2 sur $\P_2$ cette d\'efinition est
\'equivalente \`a celle donn\'ee dans \cite{st}. 
Supposons que $\ke$ soit pure. Alors $\phi$ est une immersion ferm\'ee.
Soit $\kf$ un faisceau
coh\'erent sur \m{S\times{\bf N}}. Alors on a un isomorphisme canonique
\begin{eqnarray*}
p_{S\times{\bf N}*}(p^\sharp_S(\ke)\ot p^\sharp_{\bf N}(\G)\ot
p_{S\times{\bf N}}^*(\kf)) & \simeq &
\HHom(R^np_{S\times{\bf N}*}
(p^\sharp_S(\ke^*)\ot p^\sharp_{\bf N}(\G^*)\ot p_X^*(\omega_X)),\kf)\\
& = & \HHom(\phi^*(L),\kf) .
\end{eqnarray*}
de dualit\'e relative (cf. \cite{kl}).
Soit \m{\alpha : S\to{\bf N}} la seconde composante de $\phi$. Alors, en prenant
\m{\kf=\ko_{\phi(S)}} dans ce qui pr\'ec\`ede, on voit qu'on a un
isomorphisme
\[p_{S*}(\alpha^\sharp(\G)\ot\ke) \ \simeq \ L^* .\]
Il en d\'ecoule que le morphisme canonique
\[\theta : p_S^*p_{S*}(\alpha^\sharp(\G)\ot\ke)\ot\alpha^\sharp(\G^*)\lra\ke\]
est injectif (comme morphisme de faisceaux). En utilisant par exemple le
corollaire 5.7 de \cite{sga1}, expos\'e IV, on voit que \m{\ku=\coker(\theta)}
est une famille plate de faisceaux r\'eguliers. En utilisant les r\'esultats du
chapitre \ref{f_reg} on voit que \m{\ku^*} est localement libre et que
\m{\ku^{**}/\ku} est une famille de faisceaux de $\kz$. On en d\'eduit
ais\'ement la

\sepprop

\pagebreak[2]\begin{subsub}{\bf Proposition : }\label{vmod3b}
La famille $\ke$ est pure si et seulement si $\ke^*$ l'est.
\end{subsub}

\sepprop

\pagebreak[2]\begin{subsub}{\bf Proposition : }\label{vmod4}
Si $\ke$ est pure, il existe un unique morphisme \ 
\m{f_\ke : S\to\M(\kx,\ky,\kz)}
\ tel que pour tout point ferm\'e de $S$, \m{f_\ke(s)} soit le point de 
\m{\M(\kx,\ky,\kz)} correspondant \`a l'extension large \m{\ke_s}.
\end{subsub}

\begin{proof}
Analogue \`a la proposition 2.7 de \cite{st}. On utilise la suite exacte
\[0\lra p_S^*p_{S*}(\alpha^\sharp(\G)\ot\ke)\ot\alpha^\sharp(\G^*)\lra\ke\lra\ku
\lra 0 .\]
\end{proof}

\sepprop

\pagebreak[2]\begin{subsub}{\bf Remarque : }\rm
Les r\'esultats pr\'ec\'edents s'\'etendent sans peine \`a des
cas un peu plus g\'en\'eraux o\`u {\bf M} ou {\bf N} ne sont plus
n\'ecessairement des vari\'et\'es de modules fins, mais des structures de
vari\'et\'es alg\'ebriques sur $\kx$, $\ky$ respectivement, ayant des
propri\'et\'es moins fortes. Par exemple {\bf M} ou {\bf N} peuvent \^etre des
ouverts de vari\'et\'es de modules de fibr\'es stables et 2-lisses.
\end{subsub}

\end{sub}

\sepsub

\Ssect{Fibr\'es universels}{mod_univ}

On utilise les notations de \ref{vmod3b2}. En utilisant des r\'esultats de
\cite{ra} concernant les extensions universelles (voir aussi \cite{se}, app.
III, p. 198) on montre qu'il existe un {\em fibr\'e universel 
d\'efini localement} sur \m{\M(\kx,\ky,\kz)\times X} 
(cf. \ref{modfin}). On n'en donnera la construction que dans le cas o\`u $X$ 
est une surface, qui est le plus difficile.

\sepsubsub

\pagebreak[2]\begin{subsub}{Construction du fibr\'e universel. }\rm On reprend 
les notations de \ref{vmod_hyp} et \ref{vmod3b2}. Soit
\[U \ = \ \P(\F_1^*)\times_{{\bf M}\times{\bf N}\times X_k}\cdots
\times_{{\bf M}\times{\bf N}\times X_k}
\P(\F_k^*) ,\]
et 
\[\pi_i : U\to\P(\F_i^*) , \ \ \ \ \pi_{\bf M} : U\to {\bf M} , \ \ \ \ 
\pi_{\bf N} : U\to {\bf N} ,\]
\[\pi_{X_k} : U\to X_k ,\ \ \ \ p_U : U\times X\to U , \ \ \ \
\pi : U\to{\bf M}\times{\bf N}\times X_k \]
les projections. Soit $\T'$ le faisceau universel sur \m{X_k\times X} (on a donc
\m{\T'/\Sigma_k=\T}). Soient
\[\Phi : \pi_{\bf M}^\sharp(\F)\ot p_U^*(\pi_1^*(\ko_{\P(\F_1^*)}(-1))\ot\cdots
\ot p_U^*(\pi_k^*(\ko_{\P(\F_k^*)}(-1)))\lra\pi_{X_k}^*(\T')\]
le morphisme surjectif \'evident de faisceaux sur \m{U\times X}
et \ \m{\kf_0=\ker(\Phi)}. Ce dernier est une
famille de faisceaux r\'eguliers sur $X$ param\'etr\'ee par $U$. Posons
\[\A \ = \ \EExt^1_{p_U}(\kf_0,\pi_{\bf N}^\sharp(\G^*)), \ \ \ \
\B \ = \ p_{U*}(\HHom(\kf_0,\pi_{\bf N}^\sharp(\G^*)))  .\]
On a des isomorphismes canoniques
\[\A \ \simeq \ p_{U*}(\HHom(\pi_{\bf N}^\sharp(\G),\pi_{X_k}^*(\T'))) ,\ \ \ \
\B \ \simeq \ p_{U*}(\HHom(\pi_{\bf M}^\sharp(\F),\pi_{\bf N}^\sharp(\G^*))) .\]
Soit \ \m{\E=\A\ot\B}. Si \ \m{y=(m,n,(x_i),(\phi_i))\in U} (\m{x_i\in X},
\m{\phi_i\in\P(\F_{x_i}^*)}) on a
\[\E_y \ = \ \Hom(\G_n,\ssom_{i=1}^k\ko_{x_i})\oplus\Hom(\F_m,\G_n^*) .\]
Il en d\'ecoule qu'on a \ \m{\E=\pi^*(\E')}, avec
\[\E' \ = \ p_*(\HHom(p_{\bf N}^\sharp(\G),p_{k*}(\T')))\ot
p_*(\HHom(p_{\bf M}^\sharp(\F),p_{\bf N}^\sharp(\G^*))) \]
($p$, \m{p_{\bf N}}, \m{p_k}, \m{p_{\bf M}} d\'esignant les projections de
\m{{\bf M}\times{\bf N}\times X_k\times X} sur \m{{\bf M}\times{\bf N}\times 
X_k}, et de \m{{\bf M}\times{\bf N}\times X_k} sur $\bf N$, \m{X_k} et $\bf M$
respectivement). Soit \ \m{V\subset{\bf M}\times{\bf N}\times X_k} \ un ouvert
affine. On a alors \ \m{H^i(V,\E')=\nsp} pour \m{i\geq 1}. Comme on a \
\m{R^i\pi_*(\ko_U)=0} \ pour \m{i\geq 1} il en d\'ecoule qu'on a aussi \ 
\m{H^i(\pi^{-1}(V),\E)=\nsp} pour \m{i\geq 1}. D'apr\`es \cite{ra}, lemma 2.4,
il existe une extension universelle sur \ \m{\P(\ka)_{\mid\pi^{-1}(V)}\times X}
:
\[0\lra\pi_U^\sharp(\pi_{\bf N}^\sharp(\G^*))\ot\pi_\P^*(\ko_{\P(\A)}(1))
\lra\ke\lra\pi_U^\sharp(\kf_0)\lra 0 \]
(\m{\pi_U}, \m{\pi_\P} d\'esignant les projections \m{\P(\A)\to U},
\m{\P(\A)\times X\to\P(\A)} respectivement). Mais
\m{\P(\A)_{\mid\pi^{-1}(V)}/\Sigma_k} est un ouvert de \m{\M(\kx,\ky,\kz)} (cf.
\ref{quotalg2}) et \m{\ke/\Sigma_k} est un fibr\'e universel sur \
\m{(\P(\A)_{\mid\pi^{-1}(V)}/\Sigma_k)\times X}.
\end{subsub}

\sepsubsub

\pagebreak[2]\begin{subsub}{Le cas des vari\'et\'es de dimension sup\'erieure 
\`a 2. }\rm Supposons que \m{\dim(X)>2}. Les hypoth\`eses de \ref{vmod_hyp} et 
les r\'esultats de \ref{s_ext_6} entrainent que les fibr\'es universels locaux
obtenus sont des familles compl\`etes. Dans ce cas \m{\M(\kx,\ky,\kz)} est donc
une vari\'et\'e de modules fins. On obtient ainsi de nouvelles vari\'et\'es de
modules fins de fibr\'es vectoriels non simples (on donne dans \cite{dr} des
exemples de telles vari\'et\'es sur \m{\P_2}). 
\end{subsub}
\end{sub}

\sepsub

\Ssect{Exemples sur $\P_3$}{mod_ex_p3}

Dans les exemples suivants on utilise deux types de faisceaux r\'eguliers sur 
\m{\P_3} : les faisceaux d'id\'eaux de droites ou les faisceaux r\'eguliers
construits en \ref{reg_3x3b4} comme noyaux de morphismes surjectifs
\m{E\to\ko_\ell(m)} ($\ell$ \'etant une droite de \m{\P_3}, \m{m>0} et $E$ un
fibr\'e de corr\'elation nulle). Dans ce dernier cas on prendra pour ${\bf Z}$ 
la grassmannienne des droites de $\P_3$ et pour $\T$ le fibr\'e en droites 
universel de degr\'e \m{m}. Pour toute droite $\ell$ de \m{\P_3} on a donc \ 
\m{\T_\ell=\ko_\ell(m)}. Rappelons que la vari\'et\'e de modules fins 
constitu\'ee des fibr\'es de corr\'elation nulle est
isomorphe \`a un ouvert de \m{\P_5}.  

\sepsubsub

\pagebreak[2]\begin{subsub}\label{vmod_r3}Fibr\'es de rang 3\end{subsub}

Soit \m{n>4} un entier. On consid\`ere des extensions du type
\[0\lra E(n)\lra\ke\lra\ki_\ell\lra 0\]
o\`u $\ell$ est une droite de \m{\P_3}, \m{\ki_\ell} son faisceau d'id\'eaux et
$E$ un fibr\'e de corr\'elation nulle. On montre ais\'ement que les 
propri\'et\'es de \ref{dual} et \ref{vmod_hyp} sont v\'erifi\'ees. Ici $\bf M$
est r\'eduit \`a un point (correspondant \`a $\ko$), $\bf Z$ est la
grassmannienne des droites de \m{\P_3} et \m{{\bf N}} est la vari\'et\'e
de modules des fibr\'es de corr\'elation nulle (c'est-\`a-dire la
vari\'et\'e de modules des fibr\'es stables de rang 2 et de classes de Chern
\m{c_1=0}, \m{c_2=1}).
Les fibr\'es $\ke$ sont de rang 3 et de classes de Chern \m{2n}, \m{n^2+2},
\m{2n+2}. La vari\'et\'e de modules \m{\M(\kx,\ky,\kz)} des extensions larges
du type pr\'ec\'edent est une vari\'et\'e de modules fins. On a \ 
\m{\dim(\M(\kx,\ky,\kz))=2n+14.} 
Les fibr\'es $\ke$ sont lisses, mais on a cependant
\[\dim(\Ext^2(\ke,\ke)) \ = \ 2n + 10 .\]
Pour obtenir ce r\'esultat, on part de la formule \ \m{\chi(\ke,\ke)=4n^2-3}
(obtenue gr\^ace \`a \ref{RR_P3} par exemple). On a
\[\dim(\End(\ke)) \ = \ h^0(E(n))+1 \ = \ \frac{n(n+2)(n+4)}{3}+1 ,\]
\[\dim(\Ext^1(\ke,\ke)) \ = \ \dim(\M(\kx,\ky,\kz)) \ = \ 2n+14 ,\]
\[\dim(\Ext^3(\ke,\ke)) \ = \ \dim(\Hom(\ke,\ke(-4))) \ = \ h^0(E(n-4)) \ =
\ \frac{n(n-2)(n-4)}{3} ,\]
d'o\`u on d\'eduit la dimension de \m{\Ext^2(\ke,\ke)}.

\sepsubsub

\pagebreak[2]\begin{subsub}\label{vmod_r4}Fibr\'es de rang 4\end{subsub}

Soient $m$, $n$ des entiers, avec \m{n>{\rm Max}(m,4)}.
Soient $E$ un fibr\'e de corr\'elation nulle, $\ell$ une 
droite de \m{\P_3} et \m{\pi:E\lra\ko_l(m)} un morphisme surjectif. On 
consid\`ere des extensions du type
\[0\lra E'(n)\lra\ke\lra\ker(\pi)\lra 0 ,\]
o\`u $E'$ est un fibr\'e de corr\'elation nulle. On montre ais\'ement que les 
propri\'et\'es de \ref{dual} et \ref{vmod_hyp} sont v\'erifi\'ees. Ici $\bf M$
est la vari\'et\'e de modules des fibr\'es de corr\'elation nulle, $\bf N$ est
isomorphe \`a $\bf M$ et $\bf Z$ est la grassmannienne des droites de \m{\P_3}.
Les fibr\'es $\ke$ sont de rang 4 et de classes de Chern \m{2n}, \m{n^2+3},
\m{4n-2m+2}. La vari\'et\'e de modules \m{\M(\kx,\ky,\kz)} des extensions larges
du type pr\'ec\'edent est une vari\'et\'e de modules fins. On a \ 
\m{\dim(\M(\kx,\ky,\kz))=2n+20} . Les fibr\'es $\ke$ sont lisses, mais on a 
cependant
\[\dim(\Ext^2(\ke,\ke)) \ = \ 2n - 5 \]
(d\'emonstration analogue \`a celle de \ref{vmod_r3}).
\end{sub}

\sepsec


\pagebreak[4]\section{Extensions larges sur les surfaces}
\label{Surf_ext_larg}

On suppose dans ce chapitre que $X$ est une surface.

\sepsub

\Ssect{Extensions larges g\'en\'eriques}{Ext_gen}

On consid\`ere comme dans \ref{Larg_constr} des fibr\'es vectoriels $\F$,
$\G$, et un faisceau parfait $T$ sur $X$. Soient
\[\pi : \F\lra T, \ \ \ \ \rho : \G\lra\wT\]
des morphismes surjectifs, \m{F=\ker(\pi)}, \m{G=\ker(\rho)},
de telle sorte que \ \m{\F=F^{**}}, \m{\G=G^{**}}. On se place dans le cas 
o\`u
comme dans \ref{dual}, $\pi$ et $\rho$ d\'efinissent une extension large
\[\xymatrix{ (L) \ \ \ \ \ \ \ \
0\ar[r] & G^*\ar[r]^i & \ke\ar[r]^p & F\ar[r] & 0 }\]
associ\'ee \`a \ \m{\sigma\in\Ext^1(F,G^*)} \ et l'extension duale
\[\xymatrix{ (L^*) \ \ \ \ \ \ \ \
0\ar[r] & F^*\ar[r]^{{}^tp} & \ke^*\ar[r]^{{}^ti} & G\ar[r] & 0 }\]
associ\'ee \`a \ \m{\sigma^*\in\Ext^1(G,F^*)}. 
D'apr\`es la d\'efinition \ref{larg_def}, $F^*$ et $G^*$ sont 2-lisses.
On utilise les notations de \ref{Larg0}.

D'apr\`es la proposition \ref{id4}, compte tenu des isomorphismes \
\m{\Ext^1(G^*,F)\simeq\Hom(G^*,T)}, \m{\Ext^2(F,F)\simeq H^2(\ko_X)} \ et de
la dualit\'e de Serre, la transpos\'ee de $\sigma\times$
\[H^0(\omega_X)\lra\Hom(G^{**},\wT\ot\omega_X)\]
est simplement la composition avec $\rho$.

Notons aussi que le fait que $X$ est une surface implique que $\xi_2$ et
$\xi_2^*$ sont surjectives (prop. \ref{surjx2}).

On \'etudiera plus particuli\`erement les {\em extensions larges
g\'en\'eriques}. On emploie ce terme lorsque
le faisceau de torsion $T$ est une somme directe de faisceaux
structuraux de points \break distincts :
\[T \ = \ \som_{x\in Z}\C_{x} , \]
o\`u $Z\subset X$ est fini, $\C_x$ d\'esignant le faisceau structural de
$\lbrace x\rbrace$.
Dans ce cas $\pi$ \'equivaut \`a une suite \m{(\pi_x)_{x\in Z}}, avec 
\m{\pi_x\in F_x^*}, \m{\pi_x\not = 0}. On a 
\[\wT \ = \ \som_{x\in Z}\omega_{X,x}^* .\]
On fixe, pour tout \m{x\in Z}, un isomorphisme \ \m{\omega_{X,x}\simeq\C}, ce 
qui permet d'identifier $\wT$ et \m{\som_{x\in Z}\C_{x}}, et
$\rho$ \'equivaut \`a une suite \m{(\rho_x)_{x\in Z}}, avec
\m{\rho_x\in G_x^*}, \m{\rho_x\not = 0}.

On a des isomorphismes
canoniques
\[\Ext^1(G^*,F) \ \simeq \ \som_{x\in Z}G_x^{**} , \ \ \ \
\Ext^1(F^*,G) \ \simeq \ \som_{x\in Z}F_x^{**} . \]

D'apr\`es la proposition \ref{prop_tra}, on a un isomorphisme canonique
\[\Ext^2(T,T) \ \simeq \ \som_{x\in Z}\omega_{X,x}^* \ \simeq \C^Z , \]
et le morphisme trace s'\'ecrit
\[\xymatrix@R=3pt{
\Ext^2(T,T)\ar[rr] & & H^2(\ko_X)\simeq H^0(\omega_X)^*\\
(\alpha_x)_{x\in Z}\fmaps[rr] & & (s\mapsto\sigg_{x\in Z}\alpha_xs(x))
}\]
On en d\'eduit les multiplications par $\sigma$ et \m{\sigma^*}
\[\xymatrix@R=3pt{
\times\sigma:\Ext^1(G^*,F)=\som_{x\in Z}G_x^{**}\ar[rr] & & \Ext^2(G^*,G^*)
\simeq H^0(\omega_X)^*\\
(\phi_x)_{x\in Z}\fmaps[rr] & & (s\mapsto\sigg_{x\in Z}\langle\phi_x,
\rho_x\rangle s(x))
}\]
\[\xymatrix@R=3pt{
\times\sigma^*:\Ext^1(F^*,G)=\som_{x\in Z}F_x^{**}\ar[rr] & & \Ext^2(F^*,F^*)
\simeq H^0(\omega_X)^*\\
(\psi_x)_{x\in Z}\fmaps[rr] & & (s\mapsto\sigg_{x\in Z}\langle\psi_x,
\pi_x\rangle s(x))
}\]
En ce qui concerne les morphismes canoniques \m{\Ext^1(G^*,F)\to\Ext^2(T,T)}
\hfil\break et \m{\Ext^1(F^*,G)\to\Ext^2(T,T)}, on a
\[\xymatrix@R=3pt{
\Ext^1(G^*,F)=\som_{x\in Z}G_x^{**}\ar[rr] & & \Ext^2(T,T)\\
(\phi_x)_{x\in Z}\fmaps[rr] & & 
(\langle\phi_x,\rho_x\rangle)_{x\in Z}
}\]
\[\xymatrix@R=3pt{
\Ext^1(F^*,G)=\som_{x\in Z}F_x^{**}\ar[rr] & & \Ext^2(T,T)\\
(\psi_x)_{x\in Z}\fmaps[rr] & & 
(\langle\psi_x,\pi_x\rangle)_{x\in Z}
}\]
Rappelons que $\Delta$ d\'esigne le morphisme canonique \ \m{\Ext^1(\ke,\ke)\to
A_2(\sigma)\oplus A_2(\sigma^*)} (cf. \ref{s_ext_6}). D'apr\`es la proposition
\ref{tan4}, on a
\[\Delta(\Ext^1(\ke,\ke)) \ = \ \biggl\lbrace((\phi_x),(\psi_x))\in
\som_{x\in Z}G_x^{**}\times\som_{x\in Z}F_x^{**}; \quad\quad\quad\quad\quad
\quad\quad\quad\null \]
\[ \null\quad\quad \quad\quad\langle\phi_x,\rho_x\rangle+
\langle\psi_x,\pi_x\rangle=0 \ {\rm pour\ tout} \ x\in Z\]
\[ \null\qquad\qquad\qquad\qquad \qquad\qquad\qquad\qquad{\rm et}
\sigg_{x\in Z}\langle\phi_x,\rho_x\rangle s(x)=0  \ {\rm pour\ tout} \ 
s\in H^0(\omega_X)\biggr\rbrace .\]

Si $\eta\in\Ext^1(\ke,\ke)$, on notera $s(\eta)$ (resp. $s^*(\eta)$)
le {\em support} de l'image $\zeta$ de $\eta$ dans $\Ext^1(G^*,F)$ (resp. 
$\Ext^1(F^*,G)$), c'est-\`a-dire l'ensemble des $x\in Z$ tels que la
composante selon $x$ de $\zeta$ est non nulle.

\sepsubsub

\pagebreak[2]\begin{subsub} L'action de \m{\Hom(\ke,\ke)} sur 
\m{\Ext^1(\ke,\ke)}. \rm
Soient \m{\eta\in\Ext^1(\ke,\ke)} et \m{(\phi_x)_{x\in Z}} (resp.
\m{(\psi_x)_{x\in Z}}) son image dans \m{A_2(\sigma)} (resp. \m{A_2(\sigma^*)}).
Soit \m{\lambda\in\Hom(F^{**},G^*)}, vu comme \'el\'ement de \m{\Hom(\ke,\ke)}.
Alors on a d'apr\`es la proposition \ref{act_hom}
\[\lambda\eta\in\Hom(F^{**},T)/C\pi=\biggr(\som_{x\in Z}F_x^*\biggl)/\langle
(\pi_x)\rangle ,\]
\[\lambda\eta \ = \ (\phi_x\circ\lambda_x)_{x\in Z} ,\]
\[\eta\lambda\in\Ext^2(T,G^{*})/C\rho=\biggr(\som_{x\in Z}G_x^*\biggl)/\langle
(\rho_x)\rangle ,\]
\[\lambda\eta \ = \ (\lambda_x(\psi_x))_{x\in Z} ,\]
\end{subsub} 
\end{sub}

\sepsub

\Ssect{Le produit \ $\Ext^1(\ke,\ke)\times\Ext^1(\ke,\ke)\lra\Ext^2(\ke,\ke)$ \
et l'application $\omega_2(\ke)$}
{Ext3}

On consid\`ere l'application bilin\'eaire canonique
\[\mu_0 : \Ext^1(\ke,\ke)\times\Ext^1(\ke,\ke)\lra
\Ext^2(\ke,\ke) .\]
On note $\mu$ la compos\'ee
\xmat{
\Ext^1(\ke,\ke)\times\Ext^1(\ke,\ke)\ar[r]^-{\mu_0} & 
\Ext^2(\ke,\ke)\ar[r] & \Ext^2(G^*,F)=\Ext^2(G^*,F^{**})
}
(le morphisme de droite \'etant induit par $(L)$).
D'apr\`es la proposition \ref{larg_lemm2} on a isomorphisme canonique
\[ \Ext^2(\ke,\ke) \ \simeq \ H^2(\ko_X)\oplus\Ext^2(G^*,F^{**}) .\]

Soit $\sigma\in\Ext^1(\ke,\ke)$. D'apr\`es la proposition \ref{def_tripl},
la d\'eformation double de $\ke$ d\'efinie par $\sigma$ s'\'etend en une
d\'eformation triple si et seulement si \ \m{\mu_0(\sigma,\sigma)=0}.
Mais  \m{\mu_0(\sigma,\sigma)} est toujours contenu dans le noyau de la
trace.  Il en d\'ecoule que la d\'eformation double de $\ke$ d\'efinie par 
$\sigma$ s'\'etend en une d\'eformation triple si et seulement si \ 
\m{\mu(\sigma,\sigma)=0}.

\sepprop

\pagebreak[2]\begin{subsub}{\bf Proposition : }\label{Ext3_2}
L'application $\mu$ s'annule sur \ \m{(M+N_*)\times\Ext^1(\ke,\ke)} \ et
\hfil\break
\m{\Ext^1(\ke,\ke)\times(M_*+N)}. L'application bilin\'eaire induite
\[\Ext^1(\ke,\ke)/(M+N_*)\times\Ext^1(\ke,\ke)/(M_*+N)
\lra\Ext^2(G^*,F^{**})\]
est isomorphe \`a la restriction \`a \m{A_2(\sigma)\times A_2(\sigma^*)} \ 
de l'application bilin\'eaire canonique
\[\Hom(G^*,T)\times\Ext^2(T,F^{**})\lra\Ext^2(G^*,F^{**}) .\]
\end{subsub}

\begin{proof}
 La proposition \ref{Ext3_2} se d\'emontre \`a
l'aide du diagramme commutatif suivant :
\xmat{
\Ext^1(\ke,\ke)\times\Ext^1(\ke,\ke)\flon[d]\ar[r]^-{\mu} & \Ext^2(G^*,F^{**})
\fleq[d] \\
\Ext^1(G^*,\ke)\times\Ext^1(\ke,F^{**})\ar[r] & \Ext^2(G^*,F^{**})\fleq[dd] \\
\Ext^1(G^*,\ke)\times\Ext^1(F,F^{**})\ar[u]^\simeq\ar[d] & \\
\Ext^1(G^*,F)\times\Ext^1(F,F^{**})\ar[r] &  \Ext^2(G^*,F^{**})\fleq[dd] \\
\Hom(G^*,T)\times\Ext^1(F,F^{**})\ar[u]^\simeq\ar[d] & \\
\Hom(G^*,T)\times\Ext^2(T,F^{**})\ar[r] & \Ext^2(G^*,F^{**}) }
qui d\'ecoule de l'associativit\'e des Ext et des r\'esultats du chapitre 
\ref{Surf_ext_larg}.
\end{proof}

\sepprop

Le produit
\[\Hom(G^*,T)\times\Ext^2(T,F^{**})\lra\Ext^2(G^*,F^{**})\]
s'identifie (\`a l'aide de la dualit\'e de Serre) \`a l'application canonique
\[\Hom(G^*,T)\times\Hom(F^{**},T\ot\omega_X)^*\lra
\Hom(F^{**},G^*\ot\omega_X)^* .\]
Si l'extension large est g\'en\'erique, cette application est la somme
directe des applications
\[\xymatrix@R=4pt{
G_x^{**}\times (F_x^{**}\ot\omega_{Xx}^{-1})\ar[rr] & & 
\Hom(F^{**}\ot\omega_X^{-1},G^*)^*\\
(\gamma,\phi)\fmaps[rr] & & 
(\alpha\mapsto\langle\alpha_x(\phi),\gamma\rangle)
}\]
$x$ parcourant $Z$. Cette formule g\'en\'eralise \cite{st2}, theorem (2.8).
On en d\'eduit des propri\'et\'es du module formel de $\ke$ (cf. 
\ref{Mod_Extens}) :

\sepprop

\pagebreak[2]\begin{subsub}{\bf Proposition : }\label{Ext3_3}
On suppose que 
\[H^1(F^*\ot G^*\ot\omega_X\ot\ki_Z) \ = \ H^0(\omega_X) \ = \ \nsp\]
(\m{\ki_Z} d\'esignant le faisceau d'id\'eaux de $Z$). 
Soient $\eta,\eta'\in\Ext^1(\ke,\ke)$. Alors on a \ \m{\mu(\eta,\eta')=0} \ 
si et seulement si
\[s(\eta)\cap s^*(\eta') \ = \ \emptyset .\]
Compte tenu de l'isomorphisme \
\m{\Ext^2(\ke,\ke)^*\simeq\Hom(F^{**}\ot\omega_X^{-1},G^*)}, on a
\[\ker(\omega_2(\ke)) \ = \ \Hom(F^{**}\ot\omega_X^{-1},G^*\ot\ki_Z) .\]
\end{subsub}

\sepprop

Pour tout \m{x\in Z}, soient \m{\pi_{x,0},\ldots,\pi_{x,r}} une base de
\m{F_x^*}, avec \m{\pi_{x,0}=\pi_x}, et \m{\rho_{x,0},\ldots,\rho_{x,s}} une 
base de \m{G_x^*}, avec \m{\rho_{x,0}=\rho_x}. Supposons comme dans la
proposition \ref{Ext3_3} que \m{H^0(\omega_X)=\nsp}. Alors 
\[\coker(\Delta) \ \subset \ \som_{x\in Z}(G_x^{**}\times F_x^{**})\]
est d\'efini par les \'equations \ \m{\rho_{x,0}+\pi_{x,0}=0}. Soient
\m{u_1,\ldots,u_N\in\Ext^1(\ke,\ke)^*} tels que \m{u_1,\ldots,u_N} et les
\m{\rho_{x,0},\ldots,\rho_{x,s},\pi_{x,1},\ldots,\pi_{x,r}}, \m{x\in Z},
constituent une base de \m{\Ext^1(\ke,\ke)^*}. Alors on d\'eduit de ce qui
pr\'ec\`ede le

\sepprop

\pagebreak[2]\begin{subsub}{\bf Corollaire : }\label{Ext3_4}
On suppose que 
\[H^1(F^*\ot G^*\ot\omega_X\ot\ki_Z) \ = \ H^0(\omega_X) \ = \ \nsp .\]
Soit $A$ le module formel de $\ke$. Soit
\[R \ = \ \C\bigl[u_1,\ldots,u_N,(\rho_{x,0},\ldots,\rho_{x,s},
\pi_{x,1},\ldots,\pi_{x,r})_{x\in Z}\bigr]\]
et \ \m{m_R\subset R} \ l'id\'eal maximal engendr\'e par les variables. Alors 
on a
\[A/m_A^3 \ \simeq \ R/J  ,\]
$J$ d\'esignant l'id\'eal engendr\'e par \m{m_R^3} et les \ 
\m{\rho_{x,i}\pi_{x,j}}, \m{i+j>0}, \m{\rho_{x,0}^2}.
\end{subsub}
\end{sub}

\sepsub

\Ssect{D\'eformations des extensions larges}{def_larg_0}

On suppose dans cette partie que le groupe de Picard de $X$ est isomorphe \`a 
$\Z$, le g\'en\'erateur ample $h$ \'etant identifi\'e \`a 1. On peut donc voir 
la premi\`ere classe de Chern d'un faisceau coh\'erent sur $X$ comme un entier.

Soient \m{r_0,r_1,a_0,a_1,b_0,b_1} des entiers, avec \m{r_0\geq 1}, \m{r_1\geq
1}. Pour \m{i=0,1} soit \m{{\bf M}_i} la vari\'et\'e de modules des faisceaux
semi-stables de rang \m{r_i} et de classes de Chern \m{a_i}, \m{b_i} sur $X$
(on suppose que \m{{\bf M}_i} est non vide). On s'int\'eresse \`a des extensions
larges du type
$$0\lra G^*(d)\lra\ke\lra F\lra 0 $$
o\`u \m{G^*} (resp. $F$) est semi-stable de rang \m{r_0} (resp. \m{r_1}) et de
classes de Chern \m{a_0,b_0} (resp. \m{a_1,b_1}). On note $r$, \m{c_1}, \m{c_2}
le rang et les classes de Chern de $\ke$. On a
$$r=r_0+r_1 , \ \ \ \ c_1=a_0+a_1+r_0d, \ \ \ \
c_2=\frac{r_0(r_0-1)}{2}d^2+((a_0+a_1)r_0-a_0)d+a_0a_1+b_0+b_1 .$$
Il existe toujours de telles extensions larges si \m{d\gg 0}. Dans la figure
suivante est repr\'esent\'e le polygone de Harder-Narasimhan \m{P_0} de $\ke$.

\bigskip

\hskip 2cm
\includegraphics{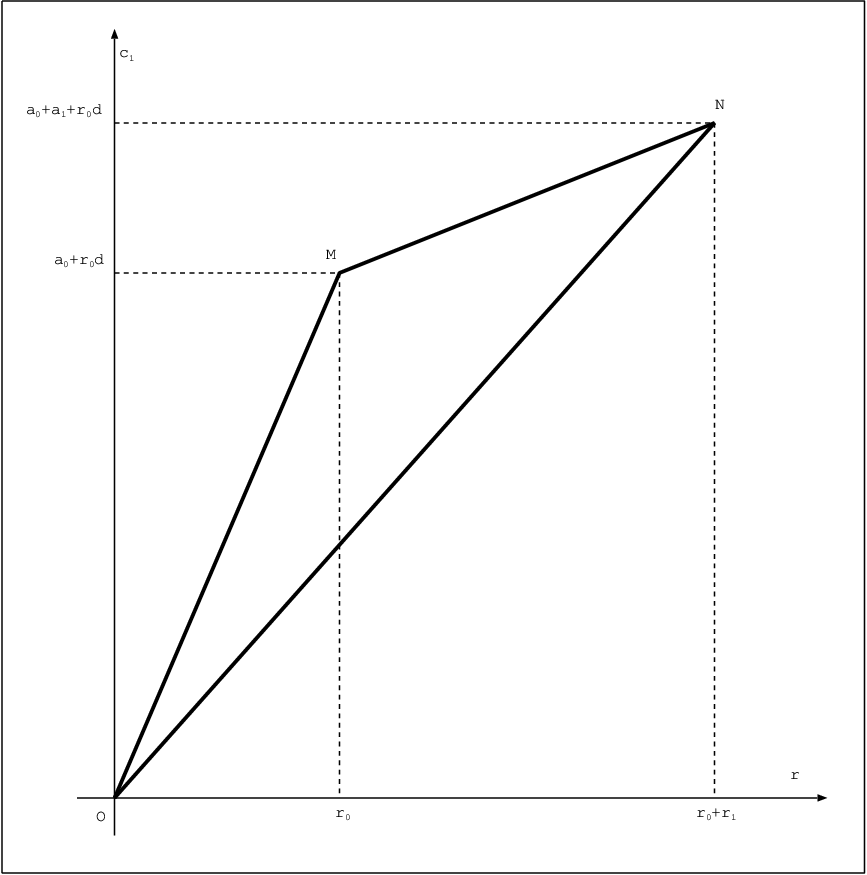}

\medskip

\centerline{{\ttx Figure 4} - Polygone de Harder-Narasimhan des extensions 
larges}

\bigskip

\pagebreak[2]\begin{subsub}{\bf Th\'eor\`eme : }\label{ext_def}
Si \m{d\gg 0} les d\'eformations de $\ke$ sont des extensions larges du m\^eme
type.
\end{subsub}

\begin{proof}
Soient $E$ un faisceau coh\'erent sans torsion sur $X$, de rang $r$ et de 
classes de Chern \m{c_1}, \m{c_2}, et
$$0=E_0\subset E_1\subset\cdots\subset E_n=E$$
sa filtration de Harder-Narasimhan. Pour \m{1\leq i\leq n} soient
$$r_i=rg(E_i/E_{i-1}), \ \ \ \ \alpha_i=c_1(E_i/E_{i-1}), \ \ \ \ 
\beta_i=c_2(E_i/E_{i-1}), \ \ \ \ \Delta_i=\Delta(E_i/E_{i-1}) .$$ 
On a donc
$$\beta_i \ = \ r_i\Delta_i+(\frac{1}{2}-\frac{1}{2r_i})\alpha_i^2 .$$
On a
$$1+c_1h+c_2h^2 \ = \ \prod_{i=1}^n(1+\alpha_ih+\beta_ih^2) ,$$
donc
$$c_1 \ = \ \sigg_{i=1}^n\alpha_i ,$$
\begin{eqnarray*}
c_2 & = & \sigg_{1\leq i<j\leq n}\alpha_i\alpha_j + \sigg_{i=1}^n\beta_i\\
    & = & \frac{1}{2}\bigl(\sigg_{i=1}^n\alpha_i\bigr)^2-
\frac{1}{2}\sigg_{i=1}^n\alpha_i^2+\sigg_{i=1}^nr_i\Delta_i+
\sigg_{i=1}^n(\frac{1}{2}-\frac{1}{2r_i})\alpha_i^2\\
    & = & \frac{1}{2}c_1^2+\sigg_{i=1}^nr_i\Delta_i-
\frac{1}{2}\sigg_{i=1}^n\frac{\alpha_i^2}{r_i} .\\
\end{eqnarray*}
On en d\'eduit que
$$\sigg_{i=1}^n\frac{\alpha_i^2}{r_i} \ = \ r_0d^2+2a_0d+
2\sigg_{i=1}^nr_i\Delta_i+a_0^2+a_1^2-2b_0-2b_1 .$$
Pour tout $i$ on a \m{\Delta_i\geq 0} ({\em in\'egalit\'e de 
Bogomolov}, cf. \cite{bo}, \cite{gi}, \cite{hu_le}). Il existe donc une 
constante $C$ (ind\'ependante de $d$ et de $E$) telle que l'on ait
$$\sigg_{i=1}^n\frac{\alpha_i^2}{r_i} \ \geq \ r_0d^2+2a_0d+C .$$
On utilise les notations de \ref{polHN}. Si \m{P\in\kp(r,c_1)} rappelons qu'on 
note aussi $P$ la fonction \m{[0,r]\to\R} associ\'ee, et \m{P'} sa d\'eriv\'ee 
(d\'efinie en dehors des sommets de $P$). Si \m{P=P(E)} l'in\'egalit\'e 
pr\'ec\'edente s'\'ecrit
$$\int^r_0P'(x)^2dx \ \geq \ r_0d^2+2a_0d+C .$$
On pose
$$m(d,r_0,r_1,a_0,a_1) \ = \ \supp_{P\in\kp(r,c_1),
P<P_0}\bigl(\int^r_0P'(x)^2dx\bigr) .$$
D'apr\`es la proposition \ref{concav1} on a
$$m(d,r_0,r_1,a_0,a_1) \ \leq \ \int^r_0P'_0(x)^2dx .$$
On suppose d'abord prouv\'e le r\'esultat suivant :

\sepprop

\pagebreak[2]\begin{subsub}{\bf Proposition : }\label{ext_def1}
On a \
$\dsp\lim_{d\to\infty}\bigl(r_0d^2+2a_0d-m(d,r_0,r_1,a_0,a_1)\bigr)=\infty .$
\end{subsub}

\sepprop

On va en d\'eduire le th\'eor\`eme \ref{ext_def}. Il faut montrer que les
polygones de Harder-Narasimhan des d\'eformations de $\ke$ sont \'egaux \`a
\m{P_0}. D'apr\'es la proposition
\ref{polHN1} les d\'eformations de $\ke$ ont un polygone de Harder-Narasimhan
\m{P\leq P_0}. On peut m\^eme supposer que \m{P<P_0} car \m{P_0} n'a que trois
sommets. On a vu que
$$\int^r_0P'(x)^2dx \ \geq \ r_0d^2+2a_0d+C ,$$
mais ceci contredit la proposition \ref{ext_def1}. Le th\'eor\`eme \ref{ext_def}
est donc prouv\'e. 

D\'emontrons maintenant la proposition \ref{ext_def1}. Consid\`erons la figure
4. Soient \m{Q_0\in OM}, \m{Q_1\in MN} des points \`a coordonn\'ees enti\`eres.
On suppose que \m{(Q_0,Q_1)\not=(O,N)}, et \m{Q_0}, \m{Q_1} distincts de $M$.
Il existe un nombre fini, ind\'ependant de $d$, de tels points. Les points
\m{Q_0}, \m{Q_1} sont enti\`erement d\'etermin\'es par leur abscisse \m{s_0},
\m{s_1} respectivement. On note \m{\kp_{d,Q_0,Q_1}} le sous-ensemble de
\m{\kp(r,c_1)} constitu\'e des polygones \m{P<P_0} contenant \m{Q_0} et \m{Q_1}.

\bigskip

\hskip 2cm
\includegraphics{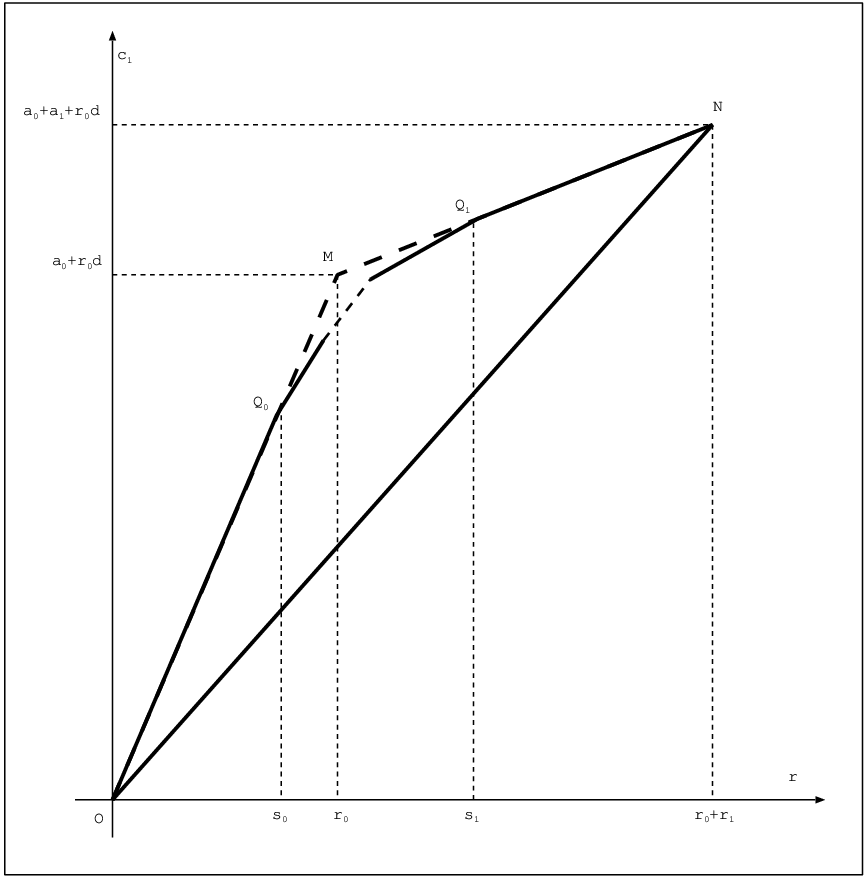}

\medskip

\centerline{{\ttx Figure 5} - Polygones de $\kp_{d,Q_0,Q_1}$}

\bigskip

Si \m{P\in\kp_{d,Q_0,Q_1}} on a donc \m{OQ_0\subset P} et \m{Q_1N\subset P}.
Soient \m{\kp_d\subset\kp(r,c_1)} le sous-ensemble constitu\'e des polygones
\m{P<P_0} et \m{\kp'_d\subset\kp_d} le compl\'ement de l'union des 
\m{\kp_{d,Q_0,Q_1}} (avec \m{(Q_0,Q_1)\not=(O,N)}).

On d\'emontre la proposition \ref{ext_def1} par r\'ecurrence sur $r$. Le premier
cas est \m{r=2}. Dans ce cas on a \m{r_0=r_1=1}, et
$$m(d,1,1,a_0,a_1) \ = \ (d+a_0-1)^2+(a_1+1)^2 \ = \ d^2+2(a_0-1)d+(a_0-1)^2+
(a_1+1)^2 .$$
Cette valeur est obtenue pour le polygone maximal de \m{\kp_d}, dont le sommet
du milieu est \m{(1,a_0+d-1)}. On a donc
$$r_0d^2+2a_0d-m(d,1,1,a_0,a_1) \ = \ 2d-(a_0-1)^2-(a_1+1)^2 ,$$
d'o\`u la proposition \ref{ext_def1} dans ce cas.

On suppose maintenant que la proposition est prouv\'ee si \m{r<R}, et que
\m{r=R>2}. Soient \m{Q_0\in OM}, \m{Q_1\in MN} comme pr\'ec\'edemment, et
\m{a'_0=a_0s_O/r_0}. Soient \m{P\in\kp_{d,Q_0,Q_1}} et \m{\overline{P}} la 
restriction de $P$ \`a \m{[s_0,s_1]}. On a
$$\int^r_0P'(x)^2dx \ = \ s_0d^2+2a'_Od+\frac{{a'_0}^2}{s_0}+
(r_1-s_1)\frac{a_1^2}{r_1}+\int^{s_1}_{s_0}\overline{P}'(x)^2dx .$$
Donc
$$r_0d^2+2a_0d-\int^r_0P'(x)^2dx \ = \ (r_0-s_0)d^2+2(a_0-a'_0)d-
\int^{s_1}_{s_0}\overline{P}'(x)^2dx-\frac{{a'_0}^2}{s_0}-
(r_1-s_1)\frac{a_1^2}{r_1} .$$
Soit \m{\rho>0}. L'hypoth\`ese de r\'ecurrence appliqu\'ee au cas \m{r=s_1-s_0}
montre que
$$(r_0-s_0)d^2+2(a_0-a'_0)d-\int^{s_1}_{s_0}\overline{P}'(x)^2dx \ \geq \ \rho$$
pour \m{d\gg 0}. On a donc
$$(1) \ \ \ \ \ \ \ \ 
r_0d^2+2a_0d-\int^r_0P'(x)^2dx \  \geq \ \rho-\frac{{a'_0}^2}{s_0}-
(r_1-s_1)\frac{a_1^2}{r_1}$$
pour \m{d\gg 0}.

On consid\`ere maintenant les polygones de \m{\kp'_d}. La pente de leur dernier
cot\'e est sup\'erieure \`a \m{a_1/r_1}. Il existe m\^eme un nombre rationnel
\m{\alpha>a_1/r_1} tel que cette pente soit sup\'erieure ou \'egale \`a
$\alpha$, quel que soit $d$. Soient $M'$ le point d'abscisse \m{r_0} de la
droite de pente $\alpha$ passant par $N$, et \m{\beta=MM'}, qui est
ind\'ependant de $d$. Alors tout polygone de \m{\kp'_d} est inf\'erieur ou
\'egal au polygone \m{P_1} suivant :

\bigskip

\hskip 2cm
\includegraphics{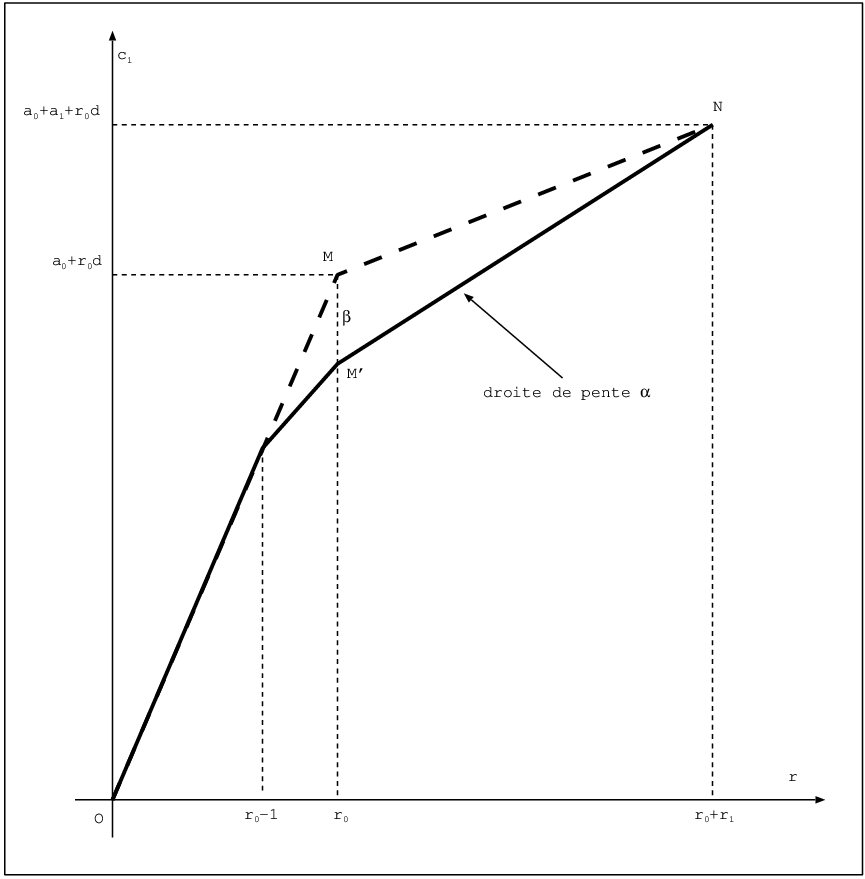}

\medskip

\centerline{{\ttx Figure 6} - Polygone $P_1$}

\bigskip

Soit \m{P\in\kp'_d}. D'apr\`es la proposition \ref{concav1} on a donc
$$\int^r_0P'(x)^2dx \ \leq \ \int^r_0P'_1(x)^2dx .$$
On a
$$\int^r_0P'_1(x)^2dx \ = \ r_0^2d+2(a_0-\beta)d+\frac{a_0^2}{r_0^2}(r_0-1)+
(\frac{a_0}{r_0}-\beta)^2+\alpha^2r_1 ,$$
d'o\`u
$$r_0d^2+2a_0d-\int^r_0P'(x)^2dx \  \geq \ 2\beta d-\frac{a_0^2}{r_0^2}(r_0-1)-
(\frac{a_0}{r_0}-\beta)^2-\alpha^2r_1 ,$$
ce qui, avec l'in\'egalit\'e \m{(1)}, d\'emontre la proposition \ref{ext_def1}.

\end{proof}

\end{sub}

\end{document}